\documentclass[12pt,a4paper,reqno]{amsart}
\usepackage{color}
\usepackage{graphics,epic}
\usepackage{amsmath,amssymb, amsthm}
\usepackage[all,2cell]{xy}

\newtheorem{theorem}{Theorem}[section]
\newtheorem*{theorem*}{Theorem}

\newtheorem{lemma}[theorem]{Lemma}
\newtheorem{proposition}[theorem]{Proposition}
\newtheorem{corollary}[theorem]{Corollary}

\newtheorem*{conjecture*}{Conjecture}

\newtheorem{remark}[theorem]{Remark}

\newcommand{\n}{\mathfrak{n}}

\renewcommand{\hat}[1]{\widehat{#1}}

\newcommand{\id}{{\rm id}}

\newcommand{\Hom}{{\rm Hom}\,}

\newcommand{\Rep}{{\rm Rep}\,}
\newcommand{\Aut}{{\rm Aut}\,}

\newcommand{\Z}{\mathbb{Z}}

\newcommand{\C}{\mathbb{C}}

\def\wt{{\rm wt}}
\def\C{{\mathbb C}}

\def\R{{\mathbb R}}

\def\Z{{\mathbb Z}}

\def\1{{\bf 1}}

\def \af{{\rm af}}
\def \mult{{\rm mult}}

\def \Hom{{\rm Hom}}

\def \pf{\noindent {\bf Proof: \,}}
\def\theequation{5.\arabic{equation}}
\def \h{\mathfrak{h}}
\def \w{\omega}
\def \g{\mathfrak{g}}

\begin{document}

\title[Twisted Verlinde formula for VOAs]{Twisted Verlinde formula for vertex operator algebras}

\author{ Chongying Dong}
\address{Chongying Dong, Department of Mathematics, University of California, Santa Cruz, CA 95064 USA}
\thanks{C. Dong was partially supported by the Simons Foundation 634104}
\email{dong@ucsc.edu}

\author{Xingjun Lin}
\address{Xingjun Lin,  School of Mathematics and Statistics, Wuhan University, Wuhan 430072,  China.}
\thanks{X. Lin was supported by China NSF grant
 12171371}
\email{linxingjun88@126.com}
\begin{abstract}
For a rational and $C_2$-cofinite vertex operator algebra $V$ with an automorphism group $G$ of prime order,  the fusion rules for twisted $V$-modules are studied,  a twisted Verlinde formula which relates fusion rules for $g$-twisted modules to the $S$-matrix in the orbifold theory is established. As an application of the twisted Verlinde formula, a twisted analogue of the Kac-Walton formula is proved,  which gives fusion rules between twisted modules of affine vertex operator algebras in terms of  Clebsch-Gordan coefficients associated to the corresponding finite dimensional simple Lie algebras.
\end{abstract}
\maketitle
\section{Introduction }
\def\theequation{1.\arabic{equation}}
\setcounter{equation}{0}

The Verlinde formula for fusion rules was proposed by E. Verlinde in \cite{V} in the framework of conformal field theory. It can also be formulated in the framework of vertex operator algebras, and has been proved  by Huang \cite{H1}. By the Verlinde formula, one can computes fusion rules of vertex operator algebras in terms of modular transformation matrices of trace functions. Explicitly, let $V$ be a rational and $C_2$-cofinite vertex operator algebra. Then  $V$-module category $\mathcal{C}_V$ is a modular tensor category under certain additional conditions \cite{H}. The fusion product coefficients of  $\mathcal{C}_V$ are called fusion rules of $V$. On the other hand, the conformal block of $V$ spanned by  trace functions on irreducible $V$-modules affords  a representation $\rho_V$ of the modular group $SL_2(\Z)$ \cite{Z}. For $S=\left(\begin{array}{cc}
0 & -1\\
1 & 0
\end{array}\right)$, the matrix $\rho_V(S)$ is called $S$-matrix of $V$.  By the Verlinde formula,  fusion rules of $V$ can be expressed in terms of $S$-matrix of $V$.

For a rational and $C_2$-cofinite vertex operator algebra $V$ and a finite automorphism group $G$ of $V$, it is important to study twisted modules of $V$ \cite{DVVV}, \cite{DLM2}, \cite{H4}. For twisted modules of $V$, one can define fusion rules between twisted modules of $V$ \cite{Xu}, \cite{H3}, \cite{DLXY}. On the other hand,  the conformal block of $V$ spanned by  trace functions on irreducible twisted $V$-modules affords a representation $\rho_{G,V}$ of the modular group $SL_2(\Z)$ \cite{DLM3}. The matrix $\rho_{G,V}(S)$ is called $S$-matrix in the orbifold theory. It is natural to conjecture that fusion rules between twisted modules  of $V$ can be expressed in terms of $S$-matrix in the orbifold theory. In case that $G$ is of order $2$, a twisted Verlinde formula has been discovered in  the framework of conformal field theory \cite{BFS}. In this paper, under the assumption that $G$ is a group of prime order,  we establish a twisted Verlinde formula, which relates  fusion rules between twisted modules  of $V$ to the $S$-matrix in the orbifold theory (see Theorem \ref{main}).

Our work is motivated partly by the work \cite{DeM} of Deshpande and Mukhopadhyay. They established a categorical Verlinde formula, which computes the fusion coefficients for  $G$-crossed modular fusion categories as defined by Turaev \cite{T}.  For a finite group $G$ and a $G$-crossed modular fusion category $\mathcal C$, one may define a categorical $\gamma$-crossed $S$-matrix of $\mathcal C$ for any $\gamma\in G$. By the categorical Verlinde formula, the fusion coefficients of $\mathcal C$ are expressed in terms of the categorical
$\gamma$-crossed $S$-matrix of $\mathcal C$. For a rational and $C_2$-cofinite vertex operator algebra $V$ and a finite automorphism group $G$ of $V$, it is expected that the category $\mathcal C_{G,V}$ of $g$-twisted $V$-modules, $g\in G$, forms a $G$-crossed  modular fusion category. Furthermore, motivated by results for untwisted modules \cite{DLN}, it is expected that the $S$-matrix in the orbifold theory and the categorical
$\gamma$-crossed $S$-matrix of $\mathcal C_{G,V}$ are the same up to a scalar.

The fusion rules of affine vertex operator algebras have been determined in \cite{Wa1, Wa2, K} by using the Verlinde formula, and in \cite{FZ} by using the bimodule theory.    In the framework of conformal field theory, an algorithm for an efficient calculation of fusion rules of twisted representations of untwisted affine Lie algebras has been proposed in \cite{QRS}. In this paper, based on the twisted Verlinde formula, we prove a twisted analogue of the Kac-Walton formula  (Theorem \ref{mainfusion}), which gives fusion rules between twisted modules of affine vertex operator algebras in terms of  Clebsch-Gordan coefficients associated to the corresponding finite dimensional simple Lie algebras. Also see \cite{HK}  for related work. To apply the twisted Verlinde formula, one has  to determine the $S$-matrix in the orbifold theory. For affine vertex operator algebras, we show that the $S$-matrix in the orbifold theory is given by the Kac-Peterson formula (see Theorem \ref{modular1}). The key point in our proof is to use the theory of orbit Lie algebras established in \cite{FSS, FRS}.

The paper is organized as follows: In Section 2, we recall some
facts about twisted modules of vertex operator algebras. In Section 3, we recall some
facts about fusion rules between twisted modules. In Section 4, we recall some
facts about modular invariance properties in orbifold theory of vertex operator algebras. In Section 5, we prove a twisted Verlinde formula for vertex operator algebras. In Section 6, as an application of the twisted Verlinde formula,  we prove a twisted analogue of the Kac-Walton formula.
\section{Basics}
\def\theequation{2.\arabic{equation}}
\setcounter{equation}{0}

In this section, we recall some basic facts about twisted modules of vertex operator algebras. Let $\left(V,Y,\mathbf{1},\omega\right)$ be a vertex operator algebra in the sense of
\cite{FLM} and \cite{FHL} (cf. \cite{Bo}, \cite{LL}).
 An \emph{automorphism} of a vertex operator
algebra $V$ is a linear isomorphism $g$ of $V$ such that $g\left(\omega\right)=\omega$
and $gY\left(v,z\right)g^{-1}=Y\left(g(v),z\right)$ for any $v\in V$ (cf. \cite{FLM}).
Denote by $\mbox{Aut}\left(V\right)$ the group of all automorphisms
of $V$.

 We next recall from \cite{FLM}, \cite{DLM2} the definition of a $g$-twisted $V$-module for  a finite order automorphism $g$ of $V$. Let $g$ be a finite order automorphism of $V$ of order $T$. Then $V$ has the following decomposition
\begin{align}
V=\oplus_{r=0}^{T-1}V^r,
\end{align}
where $V^{r}=\left\{ v\in V\ |\  g(v)=e^{-2\pi i r/T}v\right\} $ for $r\in \Z$.
Note that for $r,s\in \Z$, $V^r=V^s$ if $r\equiv s\ ({\rm mod}\; T)$.
A \emph{weak $g$-twisted $V$-module} is
a vector space $M$ with a linear map
\begin{align*}
Y_{M}(\cdot,z):\   & V\to\left(\text{End}M\right)[[z^{1/T},z^{-1/T}]]\\
 & v\mapsto Y_{M}\left(v,z\right)=\sum_{n\in \frac{1}{T}\mathbb{Z}}v_{n}z^{-n-1}\ \left(v_{n}\in\mbox{End}M\right),
\end{align*}
which satisfies the following conditions: For all $u\in V^{r}$,
$v\in V$, $w\in M$ with $0\le r\le T-1$,

\[
Y_{M}(u,z)=\sum_{n\in \frac{r}{T}+\mathbb{Z}}u_{n}z^{-n-1},
\]

\[
u_{n}w=0\  \   \mbox{ for $n$ sufficiently large},
\]

\[
Y_{M}(\mathbf{1},z)=\id_{M},
\]

\[
z_{0}^{-1}\text{\ensuremath{\delta}}\left(\frac{z_{1}-z_{2}}{z_{0}}\right)Y_{M}(u,z_{1})Y_{M}(v,z_{2})
-z_{0}^{-1}\delta\left(\frac{z_{2}-z_{1}}{-z_{0}}\right)Y_{M}(v,z_{2})Y_{M}(u,z_{1})
\]

\begin{equation}
=z_{2}^{-1}\left(\frac{z_{1}-z_{0}}{z_{2}}\right)^{-\frac{r}{T}}\delta\left(\frac{z_{1}-z_{0}}{z_{2}}\right)Y_{M}\left(Y\left(u,z_{0}\right)v,z_{2}\right),\label{Jacobi for twisted V-module}
\end{equation}
where $\delta\left(z\right)=\sum_{n\in\mathbb{Z}}z^{n}$.

An \emph{admissible $g$-twisted $V$-module}
is a weak $g$-twisted module with a $\frac{1}{T}\mathbb{Z}_{+}$-grading $M=\oplus_{n\in\frac{1}{T}\mathbb{Z}_{+}}M(n)$
such that $u_{m}M\left(n\right)\subset M\left(\mbox{wt}u-m-1+n\right)$
for homogeneous $u\in V$ and $m,n\in\frac{1}{T}\mathbb{Z}.$
A $g$-\emph{twisted $V$-module} is a weak $g$-twisted $V$-module
$M$ which carries a $\mathbb{C}$-grading induced by the spectrum
of $L(0)$, where $L(0)$ is the component operator of $Y(\omega,z)=\sum_{n\in\mathbb{Z}}L(n)z^{-n-2}.$
That is, we have $M=\bigoplus_{\lambda\in\mathbb{C}}M_{\lambda},$
where $M_{\lambda}=\left\{ w\in M\ |\  L(0)w=\lambda w\right\} $. Moreover, it is required that
$\dim M_{\lambda}<\infty$ for all $\lambda$  and for any fixed $\lambda_0,$ $M_{\frac{n}{T}+\lambda_0}=0$
for all small enough integers $n.$

In case $g=1$, we recover the notions of weak, ordinary and admissible
$V$-modules (see \cite{DLM2}).
A vertex operator algebra $V$ is said to be \emph{$g$-rational}
if the admissible $g$-twisted module category is semisimple. In particular, $V$
is said to be \emph{rational} if $V$ is $1$-rational.

If $V$ is a $g$-rational vertex
operator algebra, it is proved in \cite{DLM2} that there are only finitely many irreducible admissible
$g$-twisted $V$-modules up to isomorphism and any irreducible
admissible $g$-twisted $V$-module is ordinary. If $M=\oplus_{n\in\frac{1}{T}\mathbb{Z}_{+}}M(n)$ is an irreducible admissible $g$-twisted $V$-module,
then there is a complex number
$\lambda_{M}$ such that $L(0)|_{M(n)}=\lambda_{M}+n$ for all $n$ (cf. \cite{DLM2}).
As a convention, we assume $M(0)\ne0$, and $\lambda_{M}$ is called the {\em conformal weight} of $M.$

A vertex operator algebra $V$ is said to be {\em $C_{2}$-cofinite}
if $V/C_{2}(V)$ is finite dimensional, where $C_{2}(V)={\rm span}\{ u_{-2}v\ |\ u,v\in V\}.$ A vertex operator algebra $V=\oplus_{n\in\mathbb{Z}}V_{n}$ is said
to be of {\em CFT type} if $V_{n}=0$ for all negative integers $n$ and $V_{0}=\mathbb{C}{\bf 1}$.
The following result has been proved in \cite{ADJR}.
\begin{proposition}\label{g-rational}
 Assume that $V$ is rational and $C_{2}$-cofinite. Then  $V$ is $g$-rational
 for any finite automorphism $g$.
\end{proposition}

 Let $M=\bigoplus_{n\in\frac{1}{T}\mathbb{Z}_{+}}M(n)$
be an admissible $g$-twisted $V$-module.  Set
\[
M'=\bigoplus_{n\in\frac{1}{T}\mathbb{Z}_{+}}M\left(n\right)^{*},
\]
the {\em restricted dual}, where $M(n)^{*}=\text{Hom}_{\mathbb{C}}(M(n),\mathbb{C})$.
For $v\in V$, define a vertex operator $Y_{M'}(v,z)$ on $M'$  via
\begin{eqnarray*}
\langle Y_{M'}(v,z)f,u\rangle= \langle f,Y_{M}(e^{zL(1)}(-z^{-2})^{L(0)}v,z^{-1})u\rangle,
\end{eqnarray*}
where $\langle f,w\rangle=f(w)$ is the natural paring $M'\times M\to\mathbb{C}.$
On the other hand, if $M=\oplus_{\lambda\in \C}M_{\lambda}$ is a $g$-twisted $V$-module,
we define $M'=\oplus_{\lambda\in \C}M_{\lambda}^{*}$ and define $Y_{M'}(v,z)$ for $v\in V$ in the same way.
The following result has essentially established in \cite{FHL, Xu}.
\begin{proposition} If $(M,Y_{M})$ is an admissible $g$-twisted $V$-module, then
 $(M',Y_{M'})$ carries the structure of
an admissible $g^{-1}$-twisted $V$-module. On the other hand,  if $(M,Y_{M})$ is a $g$-twisted $V$-module, then
 $(M',Y_{M'})$ carries the structure of
a $g^{-1}$-twisted $V$-module. Moreover, $M$ is irreducible if and only if $M'$ is irreducible.
\end{proposition}

In particular, if $(M,Y_{M})$ is a $V$-module, then $(M',Y_{M'})$ is also a $V$-module. A $V$-module $M$ is said to be self-dual if $M$ and $M'$ are isomorphic. A vertex operator algebra $V$ is said to be {\em self-dual} if $V$ and $V'$ are isomorphic $V$-modules.

For any subgroup $G\le\mbox{Aut}\left(V\right)$, then the set of $G$-fixed points
$$V^{G}:=\left\{ v\in V\ |\ g\left(v\right)=v \  \mbox{ for  }g\in G\right\} $$
is a vertex operator subalgebra.
The following result has been established in \cite{CM,M2}.
\begin{theorem}\label{CM}
Assume that $V$ is a simple, rational, $C_{2}$-cofinite and self-dual vertex operator
algebra of CFT type. Then for any solvable subgroup $G$ of ${\rm Aut}(V)$, $V^{G}$ is a simple, rational, $C_{2}$-cofinite and
self-dual vertex operator algebra of CFT type.
\end{theorem}

\section{Fusion rules between twisted modules}\label{sfusion}
\def\theequation{3.\arabic{equation}}
\setcounter{equation}{0}
In this section, we recall from \cite{DLXY} some facts about fusion rules between twisted modules. Throughout this section, we assume that $V$ is a simple, rational, $C_2$-cofinite and self-dual vertex operator algebra of CFT type
and $G$ is a finite abelian automorphism group of $V$. By Theorem \ref{CM}, $V^{G}$ is a simple, rational, $C_{2}$-cofinite and
self-dual vertex operator algebra of CFT type.

First, we recall some facts about modular tensor categories from \cite{BK}, \cite{KO}. Let $\mathcal{C}$ be a modular tensor category defined as in \cite{BK}, $\1_{\mathcal{C}}$ be the unit object in $\mathcal{C}$. An {\em algebra} in $\mathcal{C}$ is an object $A\in \mathcal{C}$ along with morphisms $\mu: A\otimes A\to A$ and $\iota_A:\1_{\mathcal{C}}\hookrightarrow A$ such that the following conditions hold:

(i) {\em Associativity}. Compositions $\mu\circ (\mu\otimes \id)\circ a,\ \mu\circ (\id\otimes \mu):A\otimes (A\otimes A)\to A$ are equal, where $a$ denotes the associativity isomorphism $a: A\otimes (A\otimes A)\to (A\otimes A)\otimes A$;

(ii) {\em Unit}. Composition $\mu \circ (\iota_A\otimes A):A=\1_{\mathcal{C}}\otimes A \to A$ is equal to $\id_A$;

(iii) {\em Uniqueness of unit}. $\dim \Hom_{\mathcal{C}}(1_{\mathcal{C}}, A)=1$.

This completes the definition. We will denote the  algebra just defined by $(A, \mu, \iota_A)$ or briefly by $A$.  An algebra $A$ is called {\em commutative} if $\mu \circ c_{A,A}:A\otimes A\to A$ is equal to $\mu$, where $c$ denotes the braiding of $\mathcal{C}$.

 For  an  algebra $A$  in $\mathcal{C}$, we define the category $\Rep (A)$ as follows: objects are pairs $(M,\mu_M)$, where $M\in \mathcal{C}$ and $\mu_M: A\otimes M\to M$ is a morphism in $\mathcal{C}$ such that:
\begin{align*}
& \mu_M\circ (\mu\otimes \id)\circ a=\mu_M\circ (\id \otimes \mu_M):A\otimes (A\otimes  M)\to M;\\
& \mu_M\circ(\iota_A\otimes \id)=\id: \1\otimes M\to M.
\end{align*}
The morphisms are defined by
\begin{align*}
&\Hom_{\Rep (A)}((M_1,\mu_{M_1}), (M_2,\mu_{M_2}))\\
&\ \ \ \ \ \ \ \ \ \ \ =\{\phi\in \Hom_{\mathcal{C}}(M_1, M_2)|\mu_{M_2}\circ(\id\otimes \phi)=\phi\circ\mu_{M_1}:A\otimes M_1\to M_2\}.
\end{align*}
The following result has been established in Theorem 1.5 of \cite{KO}.
\begin{theorem}
$\Rep (A)$ is a tensor category with unit object $A$.
\end{theorem}

We also need the following results, which have been established in Theorem 1.6 of \cite{KO}.
\begin{theorem}\label{induction}
Define functor $\mathcal{F}: \mathcal{C}\to \Rep (A)$ by $\mathcal{F}(V)=A\otimes V$, $\mu_{\mathcal{F}(V)}=\mu\otimes \id$. Then\\
(1) The functor $\mathcal{F}$ is exact.\\
(2) $\mathcal{F}$ is a tensor functor.\\
(3) Define functor $G: \Rep (A)\to \mathcal{C}$ by $G(V, \mu_V)=V$. Then for any $X\in \Rep (A)$,
$$\Hom_{\Rep (A)}(F(V), X)=\Hom_{\mathcal{C}}(V, G(X)).$$
\end{theorem}

Denote by ${\mathcal C}_{V}$ the category of ordinary $V$-modules and  by ${\mathcal C}_{V^G}$ the category of ordinary $V^G$-modules.
From \cite{H}, both ${\mathcal C}_{V}$ and ${\mathcal C}_{V^{G}}$ are modular tensor
categories.  Furthermore, from  \cite{KO} and \cite{CKM},  $V$
is a commutative associative algebra in ${\mathcal C}_{V^{G}}$ as
$V=\oplus_{\chi\in {\rm irr}(G)}V^{\chi}$, where ${\rm irr}(G)$ is the set of irreducible characters of $G$ and $V^{\chi}$ are irreducible $V^{G}$-modules (cf. \cite{DJX}), i.e., simple objects in
${\mathcal C}_{V^{G}}.$ Therefore, we can consider the category ${\rm Rep}(V)$, which coincides with that defined in Definition 3.1 of \cite{DLXY}.
As a consequence, there is a categorical
tensor product functor $\boxtimes_{V}$ in the category of ${\rm Rep}(V)$, which is associative. Moreover, ${\rm Rep}(V)$ is a fusion category. In particular,
${\rm Rep}(V)$ is a semisimple category with finitely many inequivalent simple objects.
 The following results, which describe the objects in  ${\rm Rep}(V)$, have been established in \cite{DLXY} (see also \cite{Ki1,Ki2}).
\begin{proposition}\label{object}
(1) If $W$ is a $g$-twisted $V$-module with $g\in G$, then $W$ is an object of ${\rm Rep}(V).$ Furthermore, if
$W_i$ is a $g_i$-twisted $V$-module with $g_i\in G$ for $i=1,2$, then $W_1$ and $W_2$ are equivalent objects
in ${\rm Rep}(V)$ if and only if $g_1=g_2$ and $W_1\simeq W_2$ as twisted $V$-modules.\\
(2) If $W$ is a simple object in ${\rm Rep}(V),$ then $W$ is an irreducible $g$-twisted $V$-module  for some
$g\in G.$
\end{proposition}

We are now ready to define fusion rules between twisted modules. Let $g_1,g_2, g_3\in G$. For any $g_i$-twisted module $M^i$, $i=1,2,3$, as they are objects in ${\rm Rep}(V)$
by Proposition \ref{object},  $M^1\boxtimes_{V}M^2$ exists in ${\rm Rep}(V).$ Since ${\rm Rep}(V)$ is a fusion category, $M^1\boxtimes_{V}M^2$ is completely reducible. The {\em fusion rule} $N_{M^1, M^2}^{M^3}$ is defined to be the multiplicity of $M^3$ in $M^1\boxtimes_{V}M^2$.

We next recall from \cite{DLXY} some properties about fusion rules between twisted modules. Let $g_1,g_2,g_3$ be mutually commuting
automorphisms of $V$ of orders $T_{1}$, $T_2$, $T_3$, respectively.
 In this case, $V$ decomposes into
the direct sum of common eigenspaces for $g_1$ and $g_2$:
\begin{align*}
V=\bigoplus_{0\le j_{1}<T_1,\ 0\le j_{2}<T_2}V^{\left(j_{1},j_{2}\right)},
\end{align*}
where for $j_1,j_2\in \Z$, \emph{
\begin{align}\label{Vij}
V^{\left(j_{1},j_{2}\right)}=\left\{ v\in V\ |\  g_{s}(v)=e^{-2\pi ij_{s}/T_{s}},s=1,2\right\} .
\end{align}}
For any complex number $\alpha$, we define
$$(-1)^{\alpha}=e^{\alpha \pi i}.$$
 We now define intertwining operators among weak $g_{s}$-twisted modules
$(M_{s},Y_{M_{s}})$ for $s=1,2,3$.
An {\em intertwining operator of type $\binom{M_3}{M_{1}\ M_{2}}$}
is a linear map
\begin{align*}
\mathcal{Y}(\cdot,z): \  \  &  M_{1}\to\left(\text{Hom}(M_{2},M_{3})\right)\{ z\}\nonumber\\
&w\mapsto\mathcal{Y}\left(w,z\right)=\sum_{n\in\mathbb{C}}w_{n}z^{-n-1}
\end{align*}
such that for any $w^{1}\in M_{1},\ w^{2}\in M_{2}$ and for any fixed $c\in\mathbb{C}$,
\[
w_{c+n}^{1}w^{2}=0\   \   \mbox{  for $n\in\mathbb{Q}$ sufficiently large},
\]
\begin{align}
 & z_{0}^{-1}\left(\frac{z_{1}-z_{2}}{z_{0}}\right)^{j_{1}/T_{1}}\delta\left(\frac{z_{1}-z_{2}}{z_{0}}\right)
 Y_{M_{3}}(u,z_{1})\mathcal{Y}(w,z_{2})\nonumber \\
 &\   \  \   \  \
  -z_{0}^{-1}\left(\frac{z_{2}-z_{1}}{-z_{0}}\right)^{j_{1}/T_{1}}\delta\left(\frac{z_{2}-z_{1}}{-z_{0}}\right)
  \mathcal{Y}(w,z_{2})Y_{M_{2}}(u,z_{1})\nonumber \\
 =\  \  &z_{2}^{-1}\left(\frac{z_{1}-z_{0}}{z_{2}}\right)^{-j_{2}/T_{2}}\delta\left(\frac{z_{1}-z_{0}}{z_{2}}\right)
 \mathcal{Y}\left(Y_{M_{1}}(u,z_{0})w,z_{2}\right)\label{Twisted Intertwining}
\end{align}
 on $M_{2}$ for $u\in V^{(j_{1},j_{2})}$ with $j_1,j_2\in \Z$ and
$w\in M_{1}$,  and
\[
\frac{d}{dz}\mathcal{Y}(w,z)=\mathcal{Y}(L(-1)w,z).
\]
All intertwining operators of type $\binom{M_3}{M_{1}\ M_{2}}$
form a vector space, which we denote by \emph{$I_{V}\binom{M_3}{M_{1}\ M_{2}}$.} The following result has been established in Remark 2.19 and Theorem 3.6 of \cite{DLXY}.
\begin{theorem}
Let $g_1,g_2, g_3\in G$ and  $M^i$, $i=1,2,3$, be $g_i$-twisted modules. Then $N_{M_{1}, M_{2}}^{M_{3}}=\dim I_{V}\binom{M_3}{M_{1}\ M_{2}}$.
\end{theorem}

As a result, the following result has been essentially proved in \cite{Xu} (cf. Remark 2.14 of \cite{DLXY}).
\begin{proposition}\label{vanish}
If there are weak $g_{s}$-twisted modules
$(M_{s},Y_{M_{s}})$ for $s=1,2,3$ such that $N_{M_{1},\ M_{2}}^{M_{3}}>0$,
then $g_{3}=g_{1}g_{2}$.
\end{proposition}

By the similar proof as that of Proposition 2.2.2 of \cite{G}, we have the following property of fusion rules (see also Corollary 5.2 of \cite{H3}).
\begin{proposition}\label{symmetric}
Let $g_1,g_2$ be commuting finite order automorphisms of a vertex operator algebra $V$ and
let $M_i$ be a $g_i$-twisted $V$-module for $i=1,2,3$ with $g_3=g_1g_2$.
Then $N_{M_{1}, M_{2}}^{M_{3}}=N_{M_{2}\circ g_1, M_{1}}^{M_{3}}=N_{M_{2}, M_{1}\circ g_2^{-1}}^{M_{3}}.$
\end{proposition}


\section{Modular invariance in orbifold theory}
\def\theequation{4.\arabic{equation}}
\setcounter{equation}{0}
In this section we recall some results about modular invariance in orbifold theory from \cite{Z,DLM3}. Throughout this section, we assume that $V$ is a simple, rational, $C_2$-cofinite and self-dual vertex operator algebra of CFT type. By Proposition \ref{g-rational}, $V$ is $g$-rational
 for any finite automorphism $g$.

 Let $G$ be a finite automorphism group $V$. For
$g,h\in G$ and
a weak $g$-twisted $V$-module $\left(M,Y_{M}\right)$, there is a weak $h^{-1}gh$-twisted
$V$-module $\left(M\circ h,Y_{M\circ h}\right)$, where $M\circ h\cong M$
as vector spaces and $Y_{M\circ h}\left(v,z\right)=Y_{M}\left(h(v),z\right)$
for $v\in V$. This defines a right action of $G$
on the set of weak twisted $V$-modules and on isomorphism classes of weak twisted
$V$-modules.
$M$ is called \emph{$h$-stable} if $M$ and $M\circ h$ are isomorphic.

Assume that $g,h$ commute. Then $h$ acts on the set of $g$-twisted $V$-modules.
Denote by $\mathcal{M}(g)$
the equivalence classes of irreducible $g$-twisted $V$-modules and
\[\mathcal{M}(g,h)=\left\{ M\in\mathcal{M}(g)\mid M\circ h\cong M\right\} .\]
Both $\mathcal{M}(g)$ and $\mathcal{M}(g,h)$
are finite sets since $V$ is $g$-rational for all $g$.

Let $M$ be an irreducible $g$-twisted $V$-module and $G_{M}$ be
the subgroup of $G$ consisting of $h\in G$ such that $M\circ h$ and
$M$ are isomorphic. By the Schur's Lemma there is a projective representation
$\phi$ of $G_{M}$ on $M$ such that
\[
\phi\left(h\right)Y\left(u,z\right)\phi\left(h\right)^{-1}=Y\left(h(u),z\right)
\]
for $h\in G_{M}$. If $h=1$ we take $\phi\left(1\right)=1$. We will need the following result (cf. Page 144 of \cite{DXY}).
\begin{lemma}\label{stablizer}
For any admissible $g$-twisted-module
$M$, $g$ acts naturally on $M$ such that $g|_{M(n)}=e^{2\pi in}$ for $n\in\frac{1}{T}\Z.$ In particular, $g$ lies in $G_{M}$.
\end{lemma}

Let $\mathbb{H}$
be the complex upper half-plane. Here and below we set $q=e^{2\pi i\tau}$,
where $\tau\in\mathbb{H}$. Let $P\left(G\right)$ be the set of the ordered commutating pairs
in $G$. For $\left(g,h\right)\in P\left(G\right)$ and $M\in\mathcal{M}(g,h),$ set
\begin{equation}
Z_{M}\left(v,\left(g,h\right),\tau\right)=\text{tr}_{M}o\left(v\right)\phi\left(h\right)q^{L\left(0\right)-c/24}=q^{\lambda-c/24}\sum_{n\in\frac{1}{T}\mathbb{Z}_{+}}\text{tr}_{M_{\lambda+n}}o\left(v\right)\phi\left(h\right)q^{n},\label{def of trace function}
\end{equation}
where $o\left(v\right)=v_{\text{wt}v-1}$ for homogeneous $v\in V$. Then $Z_{M}\left(v,\left(g,h\right),\tau\right)$ is a holomorphic
function on $\mathbb{H}$ \cite{Z, DLM3}. We write $Z_{M}\left(v,\tau\right)=Z_{M}\left(v,\left(g,1\right),\tau\right)$
for short.

Recall that there is another vertex operator algebra $\left(V,Y\left[\ \cdot, z \ \right],\boldsymbol{1},\tilde{\omega}\right)$
associated to $V$ (see \cite{Z}). Here $\tilde{\omega}=\omega-c/24$
and for homogeneous $v\in V$,
\[
Y\left[v,z\right]=Y\left(v,e^{z}-1\right)e^{z\cdot\text{wt}v}=\sum_{n\in\mathbb{Z}}v\left[n\right]z^{n-1}.
\]
We write
\[
Y\left[\tilde{\omega},z\right]=\sum_{n\in\mathbb{Z}}L\left[n\right]z^{-n-2}.
\]
The weight of a homogeneous $v\in V$ in the second vertex operator
algebra is denoted by $\text{wt}\left[v\right].$

Let $W$ be the vector space spanned by functions $$\{Z_{M}\left(v,\left(g,h\right),\tau\right)|\left(g,h\right)\in P\left(G\right), M\in\mathcal{M}(g,h)\}.$$ Then it is proved in  \cite{DLM3} that
the dimension of $W$ is equal to $\sum_{\left(g,h\right)\in P\left(G\right)}\left|\mathcal{M}(g,h)\right|$.
Now we define an action of the modular group  $\Gamma=SL_{2}\left(\mathbb{Z}\right)$ on $W$ such
that
\[
Z_{M}|_{\gamma}\left(v,\left(g,h\right),\tau\right)=\left(c\tau+d\right)^{-\text{wt}\left[v\right]}Z_{M}\left(v,\left(g,h\right),\gamma\tau\right),
\]
where $\gamma\tau=\frac{a\tau+b}{c\tau+d},$ if $\gamma=\left(\begin{array}{cc}
a & b\\
c & d
\end{array}\right)\in\Gamma=SL\left(2,\mathbb{Z}\right).$
Then the following result has been established in \cite{DLM3,Z}.

\begin{theorem}\label{modular}
Let $V$, $G$ and $W$
be as before. Then there is a representation $\rho:\Gamma\to GL\left(W\right)$ such
that for $\left(g,h\right)\in P\left(G\right)$, $\gamma=\left(\begin{array}{cc}
a & b\\
c & d
\end{array}\right)\in\Gamma$ and $M\in\mathcal{M}(g,h),$
\[
Z_{M}|_{\gamma}\left(v,\left(g,h\right),\tau\right)=\sum_{N\in\mathcal{M}\left(g^{a}h^{c},g^{b}h^{d}\right)}\gamma_{M,N}Z_{N}\left(v,\left(g^{a}h^{c},g^{b}h^{d}\right),\tau\right),
\]
where $\rho\left(\gamma\right)=\left(\gamma_{M,N}\right)$. That is,
\[
Z_{M}\left(v,\left(g,h\right),\gamma\tau\right)=\left(c\tau+d\right)^{\wt\left[v\right]}\sum_{N\in\mathcal{M}\left(g^{a}h^{c},g^{b}h^{d}\right)}\gamma_{M,N}Z_{N}\left(v,\left(g^{a}h^{c},g^{b}h^{d}\right),\tau\right).
\]

\end{theorem}

Since the modular group $\Gamma$ is generated by $S=\left(\begin{array}{cc}
0 & -1\\
1 & 0
\end{array}\right)$ and $T=\left(\begin{array}{cc}
1 & 1\\
0 & 1
\end{array}\right)$, the representation $\rho$ is uniquely determined by $\rho\left(S\right)$
and $\rho\left(T\right).$ The matrix $\rho\left(S\right)$ is called
the \emph{$S$-matrix} of the orbifold theory. Consider special
cases of the $S$-transformation:

\[
Z_{M}\left(v,\left(g,1\right), -\frac{1}{\tau}\right)=\tau^{\text{wt}\left[v\right]}\sum_{N\in\mathcal{M}\left(1,g^{-1}\right)}S_{M, N}Z_{N}\left(v,\left(1,g^{-1}\right),\tau\right)
\]
for $M\in\mathcal{M}\left(g\right)$ and

\[
Z_{N}\left(v,\left(1,g\right),-\frac{1}{\tau}\right)=\tau^{\text{wt}\left[v\right]}\sum_{M\in\mathcal{M}\left(g\right)}S_{N, M}Z_{M}\left(v,\tau\right)
\]
for $N\in\mathcal{M}\left(1\right)$. The matrix $S=\left(S_{M,N}\right)_{M,N\in\mathcal{M}\left(1\right)}$ is
called the \emph{$S$-matrix of  $V$.} The following results have been established in \cite{H1}, \cite{DLN}.
\begin{theorem}\label{verlinde}
Let $V$ be a simple, rational, $C_2$-cofinite and self-dual vertex operator algebra of CFT type, $M^0=V, M^1, \cdots, M^s$ be all inequivalent irreducible $V$-modules. Then \\
(1) $S$ is symmetric and unitary.\\
(2) The Verlinde formula holds
$$N_{M^i, M^j}^{M^k}=\sum_{0\leq r\leq s}\frac{S_{M^i, M^r}S_{M^j, M^r}\overline{S_{M^k, M^r}}}{S_{V, M^r}}.$$
\end{theorem}
\section{Twisted Verlinde formula for vertex operator algebras}
\def\theequation{5.\arabic{equation}}
\setcounter{equation}{0}
In this section, we will prove twisted Verlinde formula for vertex operator algebras. Throughout this section, we assume that $V$ is a simple, rational, $C_2$-cofinite and self-dual vertex operator algebra of CFT type
and $\sigma$ is an automorphism of $V$ of prime order $p$. Let $G=\langle \sigma\rangle$ be the subgroup of $\Aut (V)$ generated by $\sigma$. By Theorem \ref{CM}, $V^{G}$ is a simple, rational, $C_{2}$-cofinite and
self-dual vertex operator algebra of CFT type. In addition,  by Proposition \ref{g-rational}, $V$ is $g$-rational
 for any $g\in G$.
 Our goal in this section is to prove the following twisted Verlinde formula for vertex operator algebras.
 \begin{theorem}\label{main}
 Let $V$ be a simple, rational, $C_2$-cofinite and self-dual vertex operator algebra of CFT type,
$\sigma$ be an automorphism of $V$ of prime order $p$, and $G=\langle \sigma\rangle$ be the subgroup of $\Aut (V)$ generated by $\sigma$. Let $g_1, g_2, g_3\in G$ be automorphisms of $V$ such that $g_1\neq 1$ or $g_2\neq 1$, and $g_3=g_1g_2$. Then for any $M^i\in \mathcal{M}(g_i)$, $i=1,2,3$,
$$N_{M^1, M^2}^{M^3}=\sum_{W\in \mathcal{M}(1,\sigma)}\frac{S_{M^1, W}S_{M^2, W}\overline{S_{M^3,W}}}{S_{V, W}}.$$
 \end{theorem}
To prove  Theorem \ref{main}, we need several results established in the following subsections. We will prove Theorem \ref{main} at the end of this section.
\subsection{$S$-matrix of  $V^G$} In this subsection, we recall from \cite{DRX1} some facts about $S$-matrix of $V^G$. We also prove some facts  about $S$-matrix of $V^G$ which will be used in the proof of Theorem \ref{main}.

 First, we recall from \cite{DRX} some facts about irreducible modules for $V^G$. Let $M$ be an irreducible $g$-twisted $V$-module. Recall that $G_{M}$ is
the subgroup of $G$ consisting of $h\in G$ such that $M\circ h$ and
$M$ are isomorphic. Since the order of $G$ is a prime number. Then we have $G_M=G$ or $G_M=\{1\}$. Let $M$ be an irreducible $g$-twisted $V$-module such that $G_M=G$. Then there is a projective representation $\phi$ of $G$ on $M$. Since $G$ is an abelian group, it follows that $\phi$ is a representation of $G$ (cf. Proposition 5.3 of \cite{DM}). For $0\leq s\leq p-1$, let $\Lambda_s$ be the character of $G$ such that $\Lambda_s(\sigma)=e^{\frac{-2\pi i s}{p}}$. Then $\{\Lambda_s|0\leq s\leq p-1\}$ is the set of all irreducible characters of $G$. As a result, $M$ has the following decomposition
$$M=\bigoplus_{0\leq s\leq p-1} M_{\Lambda_s},$$
where $M_{\Lambda_s}=\{w\in M|\phi(\sigma)w=\Lambda_s(\sigma)w\}$. The following results have been established in \cite{DY}, \cite{MT} (see also \cite{DRX}).
\begin{theorem}\label{irreducible1}
Let $M$ be an irreducible $g$-twisted $V$-module such that $G_M=G$, $N$ be an irreducible $h$-twisted $V$-module such that $G_N=G$. Then\\
(1) $M_{\Lambda_s}$ is nonzero for any $0\leq s\leq p-1$.\\
(2) Each $M_{\Lambda_s}$ is an irreducible $V^G$-module.\\
(3) $M_{\Lambda_s}$ and $M_{\Lambda_t}$ are equivalent  $V^G$-module if and only if $s=t$.\\
(4) If $M$ and $N$ are inequivalent, $M_{\Lambda_s}$ and $N_{\Lambda_t}$ are inequivalent  $V^G$-module for any $s, t$.
\end{theorem}
For any $g\in G$ such that $g\neq 1$ and any irreducible $g$-twisted $V$-module $M$,   we have $G_M=G$ by Lemma \ref{stablizer}. Therefore, if $M$ is an irreducible $g$-twisted $V$-module $M$ such that $G_M=\{1\}$,  $g$ must be $1$. Note that $G$ acts on the set $\mathcal{M}(1)\backslash \mathcal{M}(1,\sigma)$. Decompose $\mathcal{M}(1)\backslash \mathcal{M}(1,\sigma)$ into a disjoint union of orbits
$$\mathcal{M}(1)\backslash \mathcal{M}(1,\sigma)=\cup_{j\in J}\mathcal{O}_j.$$
For each orbit $\mathcal{O}_j$, we fix a representative $M^j$. Then it is proved in \cite{DY} that $M^j$ is an irreducible $V^G$-module. Moreover, $M^j\circ h$ and $M^j$ are isomorphic $V^G$-modules (cf. Page 146 of \cite{DRX1}). Set $$\mathcal{S}=\cup_{g\in G, g\neq 1}\mathcal{M}(g)\cup\mathcal{M}(1,\sigma),$$ the following result follows from Theorem 3.3 of \cite{DRX}.
\begin{theorem}\label{irreducible2}
The set $ \{M_{\Lambda_s}|0\leq s\leq p-1, M\in \mathcal{S}\}\cup \{M^j|j\in J\}$ gives a complete list of inequivalent irreducible  $V^G$-modules.
\end{theorem}

We now recall from \cite{DRX1} some facts about the $S$-matrix of $V^G$. For  an irreducible $g$-twisted $V$-module $M$, set $\mathcal{O}_M=\{M\circ h|h\in G\}$. Let $g_1, g_2\in G$, $M^1$ be an irreducible $g_1$-twisted $V$-module and $M^2$ be an irreducible $g_2$-twisted $V$-module. We define $C_{M^1,M^2}$ to be an least subset of $G$ such that
$$\{M^2\circ \psi|\psi\in C_{M^1,M^2}\}=\mathcal{O}_{M^2}\cap(\cup_{h\in G_{M^1}}\mathcal{M}(h, g_1^{-1})).$$
By a direct computation, we have the following results.
\begin{lemma}\label{orbit}
(1) For $M^1, M^2\in \mathcal{S}$, $C_{M^1,M^2}=\{1\}.$\\
(2) For $M^j, j\in J$, and $M^2\in \mathcal{M}(1,\sigma)$,  $C_{M^j,M^2}=\{1\}$.\\
(3) For $M^j, j\in J$, and $M^2\in\cup_{g\in G, g\neq 1}\mathcal{M}(g)$,
$C_{M^j,M^2}=\emptyset$.\\
(4) For $M^i, M^j, i, j\in J$, $C_{M^i,M^j}=\{1\}$.
\end{lemma}
Then the following results follow from Lemma \ref{orbit} and Corollary 5.4 of \cite{DRX1}.
\begin{theorem}\label{s-matrix}
Let $V$ be a simple, rational, $C_2$-cofinite and self-dual vertex operator algebra of CFT type,
$\sigma$ be an automorphism of $V$ of prime order $p$. Then the entries of $S$-matrix of $V^G$ are as follows:\\
(1) Let $M^1\in \mathcal{S}$ be an irreducible $g_1$-twisted $V$-module, $M^2\in \mathcal{S}$ be an irreducible $g_2$-twisted $V$-module. Then $S_{M^1_{\Lambda_s}, M^2_{\Lambda_t}}=\frac{1}{p}S_{M^1, M^2}\overline{\Lambda_s(g_2)}\Lambda_t(g_1^{-1}).$\\
(2) For any $M^j, j\in J$, and  $M^2\in \mathcal{M}(1,\sigma)$,
$S_{M^j, M^2_{\Lambda_t}}=S_{M^1, M^2}.$\\
(3) For any $M^j, j\in J$, and  $M^2\in\cup_{g\in G, g\neq 1}\mathcal{M}(g)$,
$S_{M^j, M^2_{\Lambda_t}}=0.$
\end{theorem}

As a consequence, we have the following results.
\begin{proposition}\label{s-matrix2}
(1) Let $M^1\in \mathcal{S}$ be an irreducible $g_1$-twisted $V$-module, $M^2\in \mathcal{S}$ be an irreducible $g_2$-twisted $V$-module. Then $S_{M^1_{\Lambda_s}, M^2_{\Lambda_t}}=\Lambda_t(g_1^{-1})S_{M^1_{\Lambda_s}, M^2_{\Lambda_0}}$.\\
(2) For any $M^j, j\in J$, and  $M^2\in \mathcal{S}$, we have $S_{M^j, M^2_{\Lambda_t}}=S_{M^j, M^2_{\Lambda_0}}$.
\end{proposition}
\subsection{Relation between fusion rules for $V$ and $V^G$}
In this subsection we establish relation between fusion rules for $V$ and $V^G$. First, we have the following result about fusion rules.
\begin{proposition}\label{gsymmetric}
Let $g_1,g_2, g_3\in G$ and  $M_i$, $i=1,2,3$, be irreducible $g_i$-twisted modules. Then for any $h\in G$, $N_{M_{1}\circ h, M_{2}\circ h}^{M_{3}\circ h}= N_{M_{1}, M_{2}}^{M_{3}}$.
\end{proposition}
\pf Define a linear map
\begin{align*}
\Phi: I_{V}\binom{M_3}{M_{1}\ M_{2}}&\to I_{V}\binom{M_3\circ h}{M_{1}\circ h\ M_{2}\circ h}\\
\mathcal{Y}(\cdot,z)&\mapsto \mathcal{Y}^h(\cdot,z),
\end{align*}
where $\mathcal{Y}^h(\cdot,z)$ is a linear map defined by
\begin{align*}
\mathcal{Y}^h(\cdot,z): \  \  &  M_{1}\circ h\to\left(\text{Hom}(M_{2}\circ h,M_{3}\circ h)\right)\{ z\}\nonumber\\
&w\mapsto\mathcal{Y}\left(w,z\right)=\sum_{n\in\mathbb{C}}w_{n}z^{-n-1}.
\end{align*}
We need to prove that $\mathcal{Y}^h(\cdot,z)$  is an intertwining operator of type $\binom{M_3\circ h}{M_{1}\circ h\ M_{2}\circ h}$. It is sufficient to the following identity
\begin{align*}
 & z_{0}^{-1}\left(\frac{z_{1}-z_{2}}{z_{0}}\right)^{j_{1}/T_{1}}\delta\left(\frac{z_{1}-z_{2}}{z_{0}}\right)
 Y_{M_{3}\circ h}(u,z_{1})\mathcal{Y}^h(w,z_{2})\nonumber \\
 &\   \  \   \  \
  -z_{0}^{-1}\left(\frac{z_{2}-z_{1}}{-z_{0}}\right)^{j_{1}/T_{1}}\delta\left(\frac{z_{2}-z_{1}}{-z_{0}}\right)
  \mathcal{Y}^h(w,z_{2})Y_{M_{2}\circ h}(u,z_{1})\nonumber \\
 =\  \  &z_{2}^{-1}\left(\frac{z_{1}-z_{0}}{z_{2}}\right)^{-j_{2}/T_{2}}\delta\left(\frac{z_{1}-z_{0}}{z_{2}}\right)
 \mathcal{Y}^h\left(Y_{M_{1}\circ h}(u,z_{0})w,z_{2}\right)
\end{align*}
holds for $u\in V^{(j_{1},j_{2})}$ with $j_1,j_2\in \Z$ and $w\in M_{1}\circ h$. By the formula (\ref{Twisted Intertwining}), we have
\begin{align*}
 & z_{0}^{-1}\left(\frac{z_{1}-z_{2}}{z_{0}}\right)^{j_{1}/T_{1}}\delta\left(\frac{z_{1}-z_{2}}{z_{0}}\right)
 Y_{M_{3}\circ h}(u,z_{1})\mathcal{Y}^h(w,z_{2})\nonumber \\
 &\   \  \   \  \
  -z_{0}^{-1}\left(\frac{z_{2}-z_{1}}{-z_{0}}\right)^{j_{1}/T_{1}}\delta\left(\frac{z_{2}-z_{1}}{-z_{0}}\right)
  \mathcal{Y}^h(w,z_{2})Y_{M_{2}\circ h}(u,z_{1})\nonumber \\
 =\  \
 & z_{0}^{-1}\left(\frac{z_{1}-z_{2}}{z_{0}}\right)^{j_{1}/T_{1}}\delta\left(\frac{z_{1}-z_{2}}{z_{0}}\right)
 Y_{M_{3}}(h(u),z_{1})\mathcal{Y}(w,z_{2})\nonumber \\
 &\   \  \   \  \
  -z_{0}^{-1}\left(\frac{z_{2}-z_{1}}{-z_{0}}\right)^{j_{1}/T_{1}}\delta\left(\frac{z_{2}-z_{1}}{-z_{0}}\right)
  \mathcal{Y}(w,z_{2})Y_{M_{2}}(h(u),z_{1})\nonumber \\
 =\  \ &z_{2}^{-1}\left(\frac{z_{1}-z_{0}}{z_{2}}\right)^{-j_{2}/T_{2}}\delta\left(\frac{z_{1}-z_{0}}{z_{2}}\right)
 \mathcal{Y}\left(Y_{M_{1}}(h(u),z_{0})w,z_{2}\right)\\
 =\  \ &z_{2}^{-1}\left(\frac{z_{1}-z_{0}}{z_{2}}\right)^{-j_{2}/T_{2}}\delta\left(\frac{z_{1}-z_{0}}{z_{2}}\right)
 \mathcal{Y}^h\left(Y_{M_{1}\circ h}(u,z_{0})w,z_{2}\right).
\end{align*}
Therefore, $\Phi$ is an injective map, this implies that $N_{M_{1}\circ h, M_{2}\circ h}^{M_{3}\circ h}\geq N_{M_{1}, M_{2}}^{M_{3}}$. As a consequence, $N_{(M_{1}\circ h)\circ h^{-1}, (M_{2}\circ h)\circ h^{-1}}^{(M_{3}\circ h)\circ h^{-1}}\geq N_{M_{1}\circ h, M_{2}\circ h}^{M_{3}\circ h}$. This implies that $ N_{M_{1}, M_{2}}^{M_{3}}\geq N_{M_{1}\circ h, M_{2}\circ h}^{M_{3}\circ h}$. Therefore, $N_{M_{1}\circ h, M_{2}\circ h}^{M_{3}\circ h}= N_{M_{1}, M_{2}}^{M_{3}}$.
\qed

We also need the following result.
\begin{proposition}\label{induction1}
Let $\mathcal{F}$ be the functor defined in Section \ref{sfusion}. Then \\
(1) For any $M\in \mathcal{S}$, and $\Lambda_i$, we have $\mathcal{F}(M_{\Lambda_i})=M$.\\
(2) For any $j\in J$, we have $\mathcal{F}(M^j)=\oplus_{W\in \mathcal{O}_j}W$.
\end{proposition}
\pf  (1) By (3) of Theorem \ref{induction}, for any $N\in \Rep (V)$, $$\Hom_{\Rep (V)}(\mathcal{F}(M_{\Lambda_i}), N)=\Hom_{\mathcal{C}_{V^G}}(M_{\Lambda_i}, N).$$
By Theorem \ref{irreducible1}, if $\Hom_{\mathcal{C}_{V^G}}(M_{\Lambda_i}, N)\neq 0$, then $N$ must be isomorphic to $M$. Moreover, $$\dim \Hom_{\mathcal{C}_{V^G}}(M_{\Lambda_i}, M)=1.$$
Since $V$ is $g$-rational for any $g\in G$, this implies that $\mathcal{F}(M_{\Lambda_i})=M$.

(2) By (3) of Theorem \ref{induction}, for any $i\in J$ and $N\in \Rep (V)$, $$\Hom_{\Rep (V)}(\mathcal{F}(M^{j}), N)=\Hom_{\mathcal{C}_{V^G}}(M^{j}, N).$$
By Theorem \ref{irreducible1}, if $\Hom_{\mathcal{C}_{V^G}}(M^{j}, N)\neq 0$, then $N$ must be isomorphic to $M^j\circ h$ for some $h\in G$. Moreover, $$\dim \Hom_{\mathcal{C}_{V^G}}(M^{j}, M^j\circ h)=1.$$
This implies that $\mathcal{F}(M^j)=\oplus_{W\in \mathcal{O}_j}W$.
\qed

As a consequence,  fusion rules for $V$ can be expressed by fusion rules for $V^G$.
\begin{theorem}\label{fusion}
(1) If $M^1, M^2, M^3\in \mathcal{S}$, then
$ N_{M^{1}, M^{2}}^{M^{3}}=\sum_{l=0}^{p-1} N_{M^{1}_{\Lambda_0}, M^{2}_{\Lambda_0}}^{M^{3}_{\Lambda_l}}.$\\
(2)  If $M^1, M^2\in \mathcal{S}$ and $M^j, j\in J$, then
$ N_{M^{1}, M^{2}}^{M^{j}}= N_{M^{1}_{\Lambda_0}, M^{2}_{\Lambda_0}}^{M^{j}}.$\\
(3) If $M^2, M^3\in \mathcal{S}$ and $M^j, j\in J$, then we have
$$ N_{M^{j}, M^{2}}^{M^{3}}=\frac{1}{p}\sum_{l=0}^{p-1} N_{M^{j}, M^{2}_{\Lambda_0}}^{M^{3}_{\Lambda_l}}.$$
\end{theorem}
\pf (1) If $M^1, M^2, M^3\in \mathcal{S}$, we have $\mathcal{F}(M^1_{\Lambda_0})=M^1$ and $\mathcal{F}(M^2_{\Lambda_0})=M^2$ by Proposition \ref{induction1}. Therefore, by Theorem \ref{induction},
\begin{align*}
M^1\boxtimes_{V}M^2&=\mathcal{F}(M^1_{\Lambda_0})\boxtimes_{V}\mathcal{F}(M^2_{\Lambda_0})\\
&=\mathcal{F}(M^1_{\Lambda_0}\boxtimes_{V^G}M^2_{\Lambda_0})\\
&=\mathcal{F}(\oplus_{W}N_{M^{1}_{\Lambda_0}, M^{2}_{\Lambda_0}}^{W}W)\\
&=\oplus_{W}N_{M^{1}_{\Lambda_0}, M^{2}_{\Lambda_0}}^{W}\mathcal{F}(W).
\end{align*}
By Proposition \ref{induction1}, this implies that $ N_{M^{1}, M^{2}}^{M^{3}}=\sum_{l=0}^{p-1} N_{M^{1}_{\Lambda_0}, M^{2}_{\Lambda_0}}^{M^{3}_{\Lambda_l}}.$

(2)  If $M^1, M^2\in \mathcal{S}$ and $M^j, j\in J$, we have $M^1\boxtimes_{V}M^2=\oplus_{W}N_{M^{1}_{\Lambda_0}, M^{2}_{\Lambda_0}}^{W}\mathcal{F}(W)$. By Proposition \ref{induction1}, this implies that
$N_{M^{1}, M^{2}}^{M^{j}}= N_{M^{1}_{\Lambda_0}, M^{2}_{\Lambda_0}}^{M^{j}}.$

(3) If $M^2, M^3\in \mathcal{S}$ and $M^j, j\in J$, we have $\mathcal{F}(M^j)=\oplus_{W\in \mathcal{O}_j}W$ and $\mathcal{F}(M^2_{\Lambda_0})=M^2$ by Proposition \ref{induction1}. Therefore, by Theorem \ref{induction},
\begin{align*}
(\oplus_{W\in \mathcal{O}_j}W)\boxtimes_{V}M^2&=\mathcal{F}(M^j)\boxtimes_{V}\mathcal{F}(M^2_{\Lambda_0})\\
&=\mathcal{F}(M^j\boxtimes_{V^G}M^2_{\Lambda_0})\\
&=\mathcal{F}(\oplus_{X}N_{M^j, M^{2}_{\Lambda_0}}^{X}X)\\
&=\oplus_{X}N_{M^j, M^{2}_{\Lambda_0}}^{X}\mathcal{F}(X).
\end{align*}
By Proposition \ref{induction1}, this implies that $ \sum_{W\in \mathcal{O}_j}N_{W, M^{2}}^{M^{3}}=\sum_{l=0}^{p-1} N_{M^j, M^{2}_{\Lambda_0}}^{M^{3}_{\Lambda_l}}.$ Since $M^2, M^3\in \mathcal{S}$, we have $\sum_{W\in \mathcal{O}_j}N_{W, M^{2}}^{M^{3}}=pN_{M^{j}, M^{2}}^{M^{3}}$ by Proposition \ref{gsymmetric}. Hence,  we have
$ N_{M^{j}, M^{2}}^{M^{3}}=\frac{1}{p}\sum_{l=0}^{p-1} N_{M^{j}, M^{2}_{\Lambda_0}}^{M^{3}_{\Lambda_l}}.$
\qed

\vskip.5cm
An an application of Theorem \ref{fusion}, we show that all intertwining operators of $V^G$-modules are restrictions of intertwining operators of twisted $V$-modules.  First, using the proof of  Proposition 11.9 of \cite{DL} gives:
\begin{proposition}\label{restriction}
Let $g_1,g_2, g_3\in G$ and  $M^i$, $i=1,2,3$, be irreducible $g_i$-twisted modules. If $\mathcal{Y}$ is a nonzero intertwining operator of type  $\binom{M^3}{M^{1}\ M^{2}}$, then $\mathcal{Y}(u, z)v\neq 0$ for any $u\in M^1, v\in M^2$.
\end{proposition}

Let $g_1,g_2, g_3\in G$ and  $M^i$, $i=1,2,3$, be irreducible $g_i$-twisted modules, $N^i$ be irreducible $V^G$-submodules of $M^i$ for $i=1,2$. Define a linear map $\Psi$ by
\begin{align*}
\Psi: I_{V}\binom{M^3}{M^{1}\ M^{2}}&\to  I_{V^G}\binom{M^3}{N^{1}\ N^{2}}\\
\mathcal{Y}(\cdot, z)&\mapsto \mathcal{Y}^r(\cdot, z),
\end{align*}
where $\mathcal{Y}^r(\cdot, z)$ is an element in $I_{V^G}\binom{M^3}{N^{1}\ N^{2}}$ such that $\mathcal{Y}^r(u, z)v=\mathcal{Y}(u, z)v$ for any $u\in N^1$ and $v\in N^2$.  As a direct consequence of Proposition \ref{restriction}, we have
\begin{proposition}\label{restriction1}
$\Psi$ is injective.
\end{proposition}

Combining Theorem \ref{fusion} and Proposition
\ref{restriction1}, we have
\begin{theorem}
If $M^1, M^2\in \mathcal{S}$, then $\Psi$ is surjective. In particular, $\Psi$ is an isomorphism.
\end{theorem}
\pf (1) If $M^1, M^2, M^3\in \mathcal{S}$, and $N^1=M^{1}_{\Lambda_s}$, $N^2=M^{2}_{\Lambda_t}$ for any $0\leq s, t\leq p-1$. By the similar argument as that in Theorem \ref{fusion}, we have
$ N_{M^{1}, M^{2}}^{M^{3}}=\sum_{l=0}^{p-1} N_{M^{1}_{\Lambda_s}, M^{2}_{\Lambda_t}}^{M^{3}_{\Lambda_l}}.$ Therefore, $\Psi$ is surjective by Proposition \ref{restriction1}.

(2)  If $M^1, M^2\in \mathcal{S}$ and $M^3=M^j$ for some $ j\in J$, $N^1=M^{1}_{\Lambda_s}$, $N^2=M^{2}_{\Lambda_t}$ for any $0\leq s, t\leq p-1$. By the similar argument as that in Theorem \ref{fusion}, we have
$ N_{M^{1}, M^{2}}^{M^{j}}= N_{M^{1}_{\Lambda_s}, M^{2}_{\Lambda_t}}^{M^{j}}.$ Therefore, $\Psi$ is surjective by Proposition \ref{restriction1}.
\qed

\vskip.25cm
We next consider the case that  $M^2, M^3\in \mathcal{S}$ and $M^1=M^j$ for some $j\in J$. Recall that $M^j$ and $M^j\circ h$ are isomorphic $V^G$-modules for $h\in G$. Then we may define a linear map  $\tilde\Psi$ by
\begin{align*}
\tilde\Psi: \bigoplus_{W\in \mathcal{O}_j}I_{V}\binom{M^3}{W\ M^{2}}&\to  I_{V^G}\binom{M^3}{N^{1}\ N^{2}}\\
(\mathcal{Y}_W(\cdot, z))_{W\in \mathcal{O}_j}&\mapsto \sum_{W\in \mathcal{O}_j}\mathcal{Y}_W^r(\cdot, z),
\end{align*}
where $\mathcal{Y}^r_W(\cdot, z)$ is an element in $I_{V^G}\binom{M^3}{N^{1}\ N^{2}}$ such that $\mathcal{Y}^r_W(u, z)v=\mathcal{Y}_W(u, z)v$ for any $u\in N^1$ and $v\in N^2$. As a direct consequence of Proposition \ref{restriction}, we have
\begin{theorem}
$\tilde\Psi$ is an isomorphism.
\end{theorem}
\pf If $M^2, M^3\in \mathcal{S}$ and $M^1=M^j$ for some $j\in J$, $N^2=M^{2}_{\Lambda_t}$ for any $0\leq t\leq p-1$. By the similar argument as that in Theorem \ref{fusion}, we have
$$ \sum_{W\in \mathcal{O}_j}N_{W, M^{2}}^{M^{3}}=\sum_{l=0}^{p-1} N_{M^j, M^{2}_{\Lambda_t}}^{M^{3}_{\Lambda_l}}.$$
Therefore, $\tilde\Psi$ is an isomorphism by Proposition \ref{restriction}.
\qed
\subsection{Proof of Theorem \ref{main}} In this subsection, we prove the main result in this paper.  Theorem \ref{main} follows from Propositions \ref{main1}, \ref{main2}, \ref{main3} and Corollary \ref{mainc}. In the following, we use $\mathfrak{L}$ to denote the set of inequivalent irreducible $V^G$-modules. We first consider the case that $M^1, M^2, M^3\in \mathcal{S}$.
\begin{proposition}\label{main1}
Let $V$, $G$ be as in Theorem \ref{main}, $g_1, g_2, g_3\in G$ be automorphisms of $V$ such that $g_1\neq 1$ or $g_2\neq 1$, and $g_3=g_1g_2$. Assume that $M^i\in \mathcal{M}(g_i)\cap \mathcal{S}$, $i=1, 2, 3$. Then we have
$$N_{M^1, M^2}^{M^3}=\sum_{W\in \mathcal{M}(1,\sigma)}\frac{S_{M^1, W}S_{M^2, W}\overline{S_{M^3,W}}}{S_{V, W}}.$$
\end{proposition}
\pf By (1) of Theorem \ref{fusion},
$ N_{M^{1}, M^{2}}^{M^{3}}=\sum_{l=0}^{p-1} N_{M^{1}_{\Lambda_0}, M^{2}_{\Lambda_0}}^{M^{3}_{\Lambda_l}}$ if $M^1, M^2, M^3\in \mathcal{S}$. By Theorems \ref{verlinde}, \ref{irreducible2}, \ref{s-matrix} and Proposition \ref{s-matrix2},
\begin{align*}
N_{M^{1}_{\Lambda_0}, M^{2}_{\Lambda_0}}^{M^{3}_{\Lambda_l}}&=\sum_{N\in\mathfrak{L}}\frac{S_{M^{1}_{\Lambda_0}, N}S_{M^{2}_{\Lambda_0}, N}\overline{S_{M^{3}_{\Lambda_l}, N}}}{S_{V^G, N}}\\
&=\sum_{W\in \mathcal{S}}\sum_{i=0}^{p-1}\frac{S_{M^{1}_{\Lambda_0}, W_{\Lambda_i}}S_{M^{2}_{\Lambda_0}, W_{\Lambda_i}}\overline{S_{M^{3}_{\Lambda_l}, W_{\Lambda_i}}}}{S_{V^G, W_{\Lambda_i}}}+\sum_{j\in J}\frac{S_{M^{1}_{\Lambda_0}, M^j}S_{M^{2}_{\Lambda_0},M^j}\overline{S_{M^{3}_{\Lambda_l}, M^j}}}{S_{V^G, M^j}}\\
&=\sum_{W\in \mathcal{S}}\sum_{i=0}^{p-1}\frac{S_{M^{1}_{\Lambda_0}, W_{\Lambda_i}}S_{M^{2}_{\Lambda_0}, W_{\Lambda_i}}\overline{S_{M^{3}_{\Lambda_l}, W_{\Lambda_i}}}}{S_{V^G, W_{\Lambda_i}}}\\
&=\sum_{W\in \mathcal{S}}\sum_{i=0}^{p-1}\frac{\Lambda_i(g_1^{-1})S _{M^{1}_{\Lambda_0}, W_{\Lambda_0}}\Lambda_i(g_2^{-1})S_{M^{2}_{\Lambda_0}, W_{\Lambda_0}}\overline{\Lambda_i(g_3^{-1})S_{M^{3}_{\Lambda_l}, W_{\Lambda_0}}}}{S_{V^G, W_{\Lambda_0}}}\\
&=\sum_{W\in \mathcal{S}}\sum_{i=0}^{p-1}\frac{S _{M^{1}_{\Lambda_0}, W_{\Lambda_0}}S_{M^{2}_{\Lambda_0}, W_{\Lambda_0}}\overline{S_{M^{3}_{\Lambda_l}, W_{\Lambda_0}}}}{S_{V^G, W_{\Lambda_0}}}\\
&=p\sum_{W\in \mathcal{S}}\frac{S _{M^{1}_{\Lambda_0}, W_{\Lambda_0}}S_{M^{2}_{\Lambda_0}, W_{\Lambda_0}}\overline{S_{M^{3}_{\Lambda_l}, W_{\Lambda_0}}}}{S_{V^G, W_{\Lambda_0}}}.
\end{align*}
Therefore,
\begin{align*}
N_{M^{1}, M^{2}}^{M^{3}}&=\sum_{l=0}^{p-1} N_{M^{1}_{\Lambda_0}, M^{2}_{\Lambda_0}}^{M^{3}_{\Lambda_l}}\\
&=\sum_{l=0}^{p-1} p\sum_{W\in \mathcal{S}}\frac{S _{M^{1}_{\Lambda_0}, W_{\Lambda_0}}S_{M^{2}_{\Lambda_0}, W_{\Lambda_0}}\overline{S_{M^{3}_{\Lambda_l}, W_{\Lambda_0}}}}{S_{V^G, W_{\Lambda_0}}}\\
&=p\sum_{W\in \mathcal{S}}\sum_{l=0}^{p-1} \frac{S _{M^{1}_{\Lambda_0}, W_{\Lambda_0}}S_{M^{2}_{\Lambda_0}, W_{\Lambda_0}}\overline{S_{M^{3}_{\Lambda_l}, W_{\Lambda_0}}}}{S_{V^G, W_{\Lambda_0}}}\\
&=p\sum_{W\in \mathcal{S}} \frac{S _{M^{1}_{\Lambda_0}, W_{\Lambda_0}}S_{M^{2}_{\Lambda_0}, W_{\Lambda_0}}}{S_{V^G, W_{\Lambda_0}}}\sum_{l=0}^{p-1}\overline{S_{M^{3}_{\Lambda_l}, W_{\Lambda_0}}}.
\end{align*}
If $W\in \mathcal{M}(g)$ and $g\neq 1$, we have
$S_{M^{3}_{\Lambda_l}, W_{\Lambda_0}}=\frac{1}{p}\overline{\Lambda_l(g)}S_{M^3, W}$ by Theorem \ref{s-matrix}. Hence, in this case,
\begin{align}\label{vanish1}
\sum_{l=0}^{p-1}\overline{S_{M^{3}_{\Lambda_l}, W_{\Lambda_0}}}=\frac{1}{p}\overline{S_{M^3, W}}\sum_{l=0}^{p-1}\Lambda_l(g)=0.
\end{align}
Thus, we have
\begin{align*}
N_{M^{1}, M^{2}}^{M^{3}}&=p\sum_{W\in \mathcal{M}(1, \sigma)} \frac{S _{M^{1}_{\Lambda_0}, W_{\Lambda_0}}S_{M^{2}_{\Lambda_0}, W_{\Lambda_0}}}{S_{V^G, W_{\Lambda_0}}}\sum_{l=0}^{p-1}\overline{S_{M^{3}_{\Lambda_l}, W_{\Lambda_0}}}\\
&=p\sum_{W\in \mathcal{M}(1, \sigma)} \frac{\frac{1}{p}S _{M^{1}, W} \frac{1}{p}S_{M^{2}, W}}{\frac{1}{p}S_{V, W}}\overline{S_{M^{3}, W}}\\
&=\sum_{W\in \mathcal{M}(1,\sigma)}\frac{S_{M^i, W}S_{M^j, W}\overline{S_{M^k,W}}}{S_{V, W}}.
\end{align*}
This completes the proof.
\qed

 We next consider the case that $M^1, M^2\in \mathcal{S}$ and $M^3=M^j$ for some $j\in J$.
\begin{proposition}\label{main2}
Let $V$, $G$ be as in Theorem \ref{main}, $g_1, g_2, g_3\in G$ be automorphisms of $V$ such that $g_1\neq 1$ or $g_2\neq 1$, and $g_3=g_1g_2$. Assume that $M^i\in \mathcal{M}(g_i)\cap \mathcal{S}$, $i=1,2$  and $M^3=M^j$ for some $j\in J$. Then we have
$$N_{M^1, M^2}^{M^3}=\sum_{W\in \mathcal{M}(1,\sigma)}\frac{S_{M^1, W}S_{M^2, W}\overline{S_{M^3,W}}}{S_{V, W}}.$$
\end{proposition}
\pf By (2) of Theorem \ref{fusion},
$ N_{M^{1}, M^{2}}^{M^{j}}= N_{M^{1}_{\Lambda_0}, M^{2}_{\Lambda_0}}^{M^{j}}$ if $M^1, M^2\in \mathcal{S}$ and $M^j, j\in J$. By Theorems \ref{verlinde}, \ref{irreducible2}, \ref{s-matrix} and Proposition \ref{s-matrix2},
\begin{align*}
N_{M^{1}_{\Lambda_0}, M^{2}_{\Lambda_0}}^{M^{j}}&=\sum_{N\in\mathfrak{L}}\frac{S_{M^{1}_{\Lambda_0}, N}S_{M^{2}_{\Lambda_0}, N}\overline{S_{M^{j}, N}}}{S_{V^G, N}}\\
&=\sum_{W\in \mathcal{S}}\sum_{i=0}^{p-1}\frac{S_{M^{1}_{\Lambda_0}, W_{\Lambda_i}}S_{M^{2}_{\Lambda_0}, W_{\Lambda_i}}\overline{S_{M^{j}, W_{\Lambda_i}}}}{S_{V^G, W_{\Lambda_i}}}+\sum_{j_1\in J}\frac{S_{M^{1}_{\Lambda_0}, M^{j_1}}S_{M^{2}_{\Lambda_0},M^{j_1}}\overline{S_{M^{j}, M^{j_1}}}}{S_{V^G, M^{j_1}}}\\
&=\sum_{W\in \mathcal{S}}\sum_{i=0}^{p-1}\frac{S_{M^{1}_{\Lambda_0}, W_{\Lambda_i}}S_{M^{2}_{\Lambda_0}, W_{\Lambda_i}}\overline{S_{M^{j}, W_{\Lambda_i}}}}{S_{V^G, W_{\Lambda_i}}}\\
&=\sum_{W\in \mathcal{S}}\sum_{i=0}^{p-1}\frac{\Lambda_i(g_1^{-1})S _{M^{1}_{\Lambda_0}, W_{\Lambda_0}}\Lambda_i(g_2^{-1})S_{M^{2}_{\Lambda_0}, W_{\Lambda_0}}\overline{S_{M^{j}, W_{\Lambda_0}}}}{S_{V^G, W_{\Lambda_0}}}\\
&=\sum_{W\in \mathcal{S}}\sum_{i=0}^{p-1}\frac{S _{M^{1}_{\Lambda_0}, W_{\Lambda_0}}S_{M^{2}_{\Lambda_0}, W_{\Lambda_0}}\overline{S_{M^{j}, W_{\Lambda_0}}}}{S_{V^G, W_{\Lambda_0}}}\\
&=p\sum_{W\in \mathcal{S}}\frac{S _{M^{1}_{\Lambda_0}, W_{\Lambda_0}}S_{M^{2}_{\Lambda_0}, W_{\Lambda_0}}\overline{S_{M^{j}, W_{\Lambda_0}}}}{S_{V^G, W_{\Lambda_0}}}=p\sum_{W\in \mathcal{M}(1, \sigma)}\frac{S _{M^{1}_{\Lambda_0}, W_{\Lambda_0}}S_{M^{2}_{\Lambda_0}, W_{\Lambda_0}}\overline{S_{M^{j}, W_{\Lambda_0}}}}{S_{V^G, W_{\Lambda_0}}}\\
&=p\sum_{W\in \mathcal{M}(1, \sigma)}\frac{\frac{1}{p}S_{M^{1}, W}\frac{1}{p}S_{M^{2}, W}\overline{S_{M^{j}, W}}}{\frac{1}{p}S_{V, W}}=\sum_{W\in \mathcal{M}(1, \sigma)}\frac{S_{M^{1}, W}S_{M^{2}, W}\overline{S_{M^{j}, W}}}{S_{V, W}}.
\end{align*}

This completes the proof.
\qed

 Finally, we consider the case that $M^2, M^3\in \mathcal{S}$ and $M^1=M^j$ for some $j\in J$.
\begin{proposition}\label{main3}
Let $V$, $G$ be as in Theorem \ref{main}, $g_1, g_2, g_3\in G$ be automorphisms of $V$ such that $g_1\neq 1$ or $g_2\neq 1$, and $g_3=g_1g_2$. Assume that $M^i\in \mathcal{M}(g_i)\cap \mathcal{S}$, $i=2,3$  and $M^1=M^j$ for some $j\in J$. Then we have
$$N_{M^1, M^2}^{M^3}=\sum_{W\in \mathcal{M}(1,\sigma)}\frac{S_{M^1, W}S_{M^2, W}\overline{S_{M^3,W}}}{S_{V, W}}.$$
\end{proposition}
\pf By (3) of Theorem \ref{fusion},
$ N_{M^{j}, M^{2}}^{M^{3}}=\frac{1}{p}\sum_{l=0}^{p-1} N_{M^{j}, M^{2}_{\Lambda_0}}^{M^{3}_{\Lambda_l}}$ if $M^2, M^3\in \mathcal{S}$ and $M^j, j\in J$. By Theorems \ref{verlinde}, \ref{irreducible2}, \ref{s-matrix} and Proposition \ref{s-matrix2},
\begin{align*}
N_{M^{j}, M^{2}_{\Lambda_0}}^{M^{3}_{\Lambda_l}}&=\sum_{N\in\mathfrak{L}}\frac{S_{M^j, N}S_{M^{2}_{\Lambda_0}, N}\overline{S_{M^{3}_{\Lambda_l}, N}}}{S_{V^G, N}}\\
&=\sum_{W\in \mathcal{S}}\sum_{i=0}^{p-1}\frac{S_{M^{j}, W_{\Lambda_i}}S_{M^{2}_{\Lambda_0}, W_{\Lambda_i}}\overline{S_{M^{3}_{\Lambda_l}, W_{\Lambda_i}}}}{S_{V^G, W_{\Lambda_i}}}+\sum_{j_1\in J}\frac{S_{M^{j}, M^{j_1}}S_{M^{2}_{\Lambda_0},M^{j_1}}\overline{S_{M^{3}_{\Lambda_l}, M^{j_1}}}}{S_{V^G, M^{j_1}}}\\
&=\sum_{W\in \mathcal{S}}\sum_{i=0}^{p-1}\frac{S_{M^{j}, W_{\Lambda_i}}S_{M^{2}_{\Lambda_0}, W_{\Lambda_i}}\overline{S_{M^{3}_{\Lambda_l}, W_{\Lambda_i}}}}{S_{V^G, W_{\Lambda_i}}}\\
&=\sum_{W\in \mathcal{S}}\sum_{i=0}^{p-1}\frac{S _{M^{j}, W_{\Lambda_0}}\Lambda_i(g_2^{-1})S_{M^{2}_{\Lambda_0}, W_{\Lambda_0}}\overline{\Lambda_i(g_3^{-1})S_{M^{3}_{\Lambda_l}, W_{\Lambda_0}}}}{S_{V^G, W_{\Lambda_0}}}\\
&=\sum_{W\in \mathcal{S}}\sum_{i=0}^{p-1}\frac{S_{M^{j}, W_{\Lambda_0}}S_{M^{2}_{\Lambda_0}, W_{\Lambda_0}}\overline{S_{M^{3}_{\Lambda_l}, W_{\Lambda_0}}}}{S_{V^G, W_{\Lambda_0}}}=\sum_{W\in \mathcal{S}}p\frac{S_{M^{j}, W_{\Lambda_0}}S_{M^{2}_{\Lambda_0}, W_{\Lambda_0}}\overline{S_{M^{3}_{\Lambda_l}, W_{\Lambda_0}}}}{S_{V^G, W_{\Lambda_0}}}.
\end{align*}
By the formula (\ref{vanish1}) and Theorem \ref{s-matrix}, we have
\begin{align*}
N_{M^{j}, M^{2}}^{M^{3}}
&=\frac{1}{p}\sum_{l=0}^{p-1} N_{M^{j}, M^{2}_{\Lambda_0}}^{M^{3}_{\Lambda_l}}=\sum_{l=0}^{p-1}\sum_{W\in \mathcal{S}}\frac{S_{M^{j}, W_{\Lambda_0}}S_{M^{2}_{\Lambda_0}, W_{\Lambda_0}}\overline{S_{M^{3}_{\Lambda_l}, W_{\Lambda_0}}}}{S_{V^G, W_{\Lambda_0}}}\\
&=\sum_{W\in \mathcal{S}}\frac{S_{M^{j}, W_{\Lambda_0}}S_{M^{2}_{\Lambda_0}, W_{\Lambda_0}}}{S_{V^G, W_{\Lambda_0}}}\sum_{l=0}^{p-1}\overline{S_{M^{3}_{\Lambda_l}, W_{\Lambda_0}}}\\
&=\sum_{W\in \mathcal{M}(1, \sigma)}\frac{S_{M^{j}, W_{\Lambda_0}}S_{M^{2}_{\Lambda_0}, W_{\Lambda_0}}}{S_{V^G, W_{\Lambda_0}}}\sum_{l=0}^{p-1}\overline{S_{M^{3}_{\Lambda_l}, W_{\Lambda_0}}}\\
&=\sum_{W\in \mathcal{M}(1, \sigma)}\frac{S_{M^{j}, W_{\Lambda_0}}S_{M^{2}_{\Lambda_0}, W_{\Lambda_0}}}{S_{V^G, W_{\Lambda_0}}}\overline{S_{M^{3}, W}}\\
&=\sum_{W\in \mathcal{M}(1, \sigma)}\frac{S_{M^{j}, W}\frac{1}{p}S_{M^{2}, W}}{\frac{1}{p}S_{V, W}}\overline{S_{M^{3}, W}}=\sum_{W\in \mathcal{M}(1,\sigma)}\frac{S_{M^1, W}S_{M^2, W}\overline{S_{M^3,W}}}{S_{V, W}}.
\end{align*}
This completes the proof.
\qed

By Proposition \ref{symmetric}, we have the following result.
\begin{corollary}\label{mainc}
Let $V$, $G$ be as in Theorem \ref{main}, $g_1, g_2, g_3\in G$ be automorphisms of $V$ such that $g_1\neq 1$ or $g_2\neq 1$, and $g_3=g_1g_2$. Assume that $M^i\in \mathcal{M}(g_i)\cap \mathcal{S}$, $i=1,3$  and $M^2=M^j$ for some $j\in J$. Then we have
$$N_{M^1, M^2}^{M^3}=\sum_{W\in \mathcal{M}(1,\sigma)}\frac{S_{M^1, W}S_{M^2, W}\overline{S_{M^3,W}}}{S_{V, W}}.$$
\end{corollary}

\section{Fusion rules between twisted modules of affine vertex operator algebras}
\def\theequation{6.\arabic{equation}}
\setcounter{equation}{0}
In this section, we will determine the $S$-matrix in the orbifold theory of affine vertex operator algebras. Furthermore, we shall prove a twisted analogue of the Kac-Walton formula, which gives fusion rules between twisted modules of affine vertex operator algebras in terms of Clebsch-Gordan coefficients associated to the corresponding finite dimensional simple Lie algebras.
\subsection{Definitions and properties about Kac-Moody algebras}
In this subsection, we recall the modular transformations of characters of integrable representations of affine Kac-Moody algebras, which play an important role in determining the $S$-matrix in the orbifold theory of affine vertex operator algebras.
\subsubsection{Definitions about Kac-Moody algebras} Let $A=(a_{i,j})_{i,j=1}^n$ be a generalized Cartan matrix and $(\h, \Pi, \Pi^{\vee})$ be a realization of $A$ as defined in \cite{K}. In particular, $\h$ is a complex vector space, $\Pi=\{\alpha_1, \cdots, \alpha_n\}\subset \h^*$ and $\Pi^\vee=\{\alpha_1^{\vee}, \cdots, \alpha_n^{\vee}\}\subset \h$ are subsets of $\h^*$ and $\h$, respectively, such that both sets $\Pi$ and $\Pi^{\vee}$ are linearly independent.  Let $\g(A)$ be the Kac-Moody algebra as defined in \cite{K}. We let $e_i, f_i \ (i=1, \cdots, n)$  denote the Chevalley generators of $\g(A)$.

It is proved in \cite{K} that there exists a nondegenerate symmetric invariant bilinear form $(,)$ on $\g(A)$. In the following, we will normalize the form $(,)$ such that $(\alpha, \alpha)=2$ for any long root $\alpha$ of $\g(A)$. It is proved in \cite{K} that the bilinear form $(,)$ is nondegenerate on $\h$. Thus, we have an isomorphism $\nu:\h\to \h^*$ defined by
$$\langle \nu(h), h_1\rangle=\nu(h)(h_1)=(h, h_1),\ \ \  h, h_1\in \h,$$
where $\langle \cdot, \cdot\rangle$ denotes the pairing between a vector space and its dual. Via the isomorphism $\nu:\h\to \h^*$,  one may obtain the induced bilinear form $(,)$ on $\h^*$.

We next recall the notion of the Weyl group of $\g(A)$. For each $i=1,\cdots, n$, we define the fundamental reflection $r_i$ of the space $\h^*$ by
$$r_i(\lambda)=\lambda-\langle \lambda, \alpha^{\vee}_i\rangle\alpha_i, \ \ \ \lambda\in \h^*.$$
The subgroup  $W$ of $GL(\h^*)$ generated by all fundamental reflections is called  the {\em Weyl group} of $\g(A)$. Through the isomorphism $\nu:\h\to \h^*$, we have an action of  the Weyl group $W$ on $\h$ such that
 $$r_i(h)=h-\langle \alpha_i, h\rangle \alpha_i^{\vee}, \ \ \ h\in \h.$$
 For $w\in W$, we use $l(w)$ to denote the length of $w$ and set $\epsilon(w)=(-1)^{l(w)}$.
\subsubsection{Properties about affine Kac-Moody algebras} Let $A$ be a generalized Cartan matrix of affine type of order $l+1$, $S(A)$ be its Dynkin diagram. Let $a_0, a_1, \cdots, a_l$ be the numerical labels of $S(A)$ as in \cite{K}. Then $a_0=1$ unless $A$ is of type $A_{2l}^{(2)}$. We let $a_0^{\vee}, a_1^{\vee}, \cdots, a_l^{\vee}$ denote the numerical labels of $S(A^T)$. The number $h^\vee=\sum_{i=0}^la_i^\vee$ is called the {\em dual Coxter number} of the matrix $A$. It is proved in \cite{K} that the center of $\g(A)$ is 1-dimensional and is spanned by
\begin{align}\label{k}
K=\sum_{i=0}^la_i^\vee\alpha_i^\vee.
\end{align}

Fix an element $d\in \h$ which satisfies the following conditions: For $i=1, \cdots, l$,
\begin{align}\label{d}
\langle \alpha_i, d\rangle=0,\ \ \, \langle \alpha_0, d\rangle=1.
\end{align}
Then $\alpha_0^\vee, \cdots, \alpha_l^\vee, d$ form a basis of $\h$. To give a basis of $\h^*$, we define an element $\Lambda_0\in \h^*$ by the following conditions: For $i=0, \cdots, l$,
\begin{align}
\langle \Lambda_0, \alpha_i^\vee\rangle=\delta_{0, i},\ \ \, \langle \Lambda_0, d\rangle=0.
\end{align}
Then $\alpha_0, \cdots, \alpha_l, \Lambda_0$ form a basis of $\h^*$.
Define
\begin{align}
\delta=\sum_{i=0}^la_i\alpha_i
 \end{align}
 and denote by $\bar{\h}^*$ the linear span over $\C$ of $\alpha_1, \cdots, \alpha_l$. Then we have $$\h^*=\bar{\h}^*\oplus(\C\delta+\C\Lambda_0).$$

 Denote by $\bar\g$ the subalgebra of $\g$ generated by $e_i$ and $f_i$ with $i=1, \cdots, l$. Then it is known \cite{K} that $\bar \g$ is a Kac-Moody algebra associated to the matrix $\bar A$ obtained from $A$ by deleting the $0$th row and column. The elements $e_i, f_i (i=1, \cdots, l)$ are the Chevalley generators of $\bar \g$, and $\bar\h=\bar \g\cap \h$ is its Cartan subalgebra. Let $\overline{W}$ be the Weyl group of $\bar \g$. It is proved in \cite{K} that $\overline{W}$ may be identified with the subgroup of the Weyl group of $\g(A)$ generated by $r_1, \cdots, r_l$.

 We next recall the description of the Weyl group of $\g(A)$ given in \cite{K}. For any $\alpha\in \bar{\h}^*$, we define the following endomorphism $t_{\alpha}$ of the vector space $\h^*$: For any $\lambda\in \h^*$,
 $$t_{\alpha}(\lambda)=\lambda+\langle \lambda, K\rangle -((\lambda, \alpha)+\frac{1}{2}(\alpha,\alpha)\langle \lambda, K\rangle)\delta.$$
 Then it is proved in \cite{K} that
 \begin{align}
 t_{\alpha}t_{\beta}=t_{\alpha+\beta}\text{ and  }t_{w(\alpha)}=wt_{\alpha}w^{-1}
  \end{align}
  for any $w\in \overline{W}$. Following \cite{K}, we define the following important lattice $M\subset \bar{\h}_{\R}^*$, where $\bar{\h}_{\R}$ is the linear span over $\R$ of $\alpha_1^\vee, \cdots, \alpha_l^\vee$.
 Let $\theta=\delta-a_0\alpha_0=\sum_{i=1}^la_i\alpha_i$. Then we have $(\theta, \theta)=2a_0$ (see (6.4.1) of \cite{K}). Let $\theta^\vee=\frac{1}{a_0}\nu^{-1}(\theta)$ and $\Z(\overline{W}\cdot\theta^\vee)$ be the lattice in $\bar{\h}_{\R}$ spanned over $\Z$ by the set $\overline{W}\cdot\theta^\vee$. Set $M=\nu(\Z(\overline{W}\cdot\theta^\vee))$. Then the following result has been proved in Proposition 6.5 of \cite{K}.
 \begin{proposition}\label{Weylg}
 Let $T$ be the subgroup of $GL(\h^*)$ generated by $t_{\alpha}$, $\alpha\in M$. Then $W=\overline{W}\ltimes T$.
 \end{proposition}
\subsubsection{Adjacent root systems} For a twisted affine Cartan matrix $A^\dag$ not of type $A_{2l}^{(2)}$, its {\em adjacent Cartan matrix} $A'=(a'_{i,j})_{i,j=0,\cdots, l}$ is defined to  be an affine Cartan matrix, of type $D_{l+1}^{(2)}$, $A_{2l-1}^{(2)}$, $E_6^{(2)}$, or $D_4^{(3)}$ according as the type of $A^\dag$ is $A_{2l-1}^{(2)}$, $D_{l+1}^{(2)}$, $E_6^{(2)}$, or $D_4^{(3)}$. For the Kac-Moody algebra $\g(A^\dag)$, we have notions $$\nu^\dag, a_i^\dag, a_i^{\dag\vee}, \alpha_i^\dag, \alpha_i^{\dag\vee}, \delta^\dag, d^\dag, \Lambda^\dag_0, K^\dag, W^\dag, M^\dag; \text{ etc.} $$ Similarly, for the Kac-Moody algebra $\g(A')$, we have notions
$$\nu',  a_i', a_i'^{\vee}, \alpha'_i, \alpha_i'^{ \vee}, \delta', d', \Lambda'_0, K', W', M'; \text{ etc.}$$ Let $\h^\dag$ and $\h'$ be the Cartan subalgebras of $\g(A^\dag)$ and $\g(A')$, respectively. We define a linear isomorphism $\varphi: \h^{'*}\to \h^{\dag*}$  such that
\begin{align}\label{varphi}
&\varphi(\delta')=\frac{1}{r}\delta^\dag,\ \ \ \varphi(\Lambda_0')=\Lambda_0^\dag,\notag\\
&\varphi(\alpha_i')=\left\{\begin{array}{lll}
\frac{a_i^{\dag}}{a_i^{\dag\vee}}\alpha_i^\dag&\text{if}& A=A_{2l-1}^{(2)}, D_{l+1}^{(2)};\\
\frac{a_{l+1-i}^{\dag}}{a_{l+1-i}^{\dag\vee}}\alpha_{l+1-i}^\dag&\text{if}& A=E_{6}^{(2)}, D_{4}^{(3)};
\end{array}\right.
\end{align}
for $1\leq i\leq l$, where $r$ is the number such that $A^\dag$ belongs to Table Aff $r$ in \cite{K}.

Let $(,)^\dag$ and $(,)'$ be the   bilinear forms of $\h^{\dag*}$ and $\h^{'*}$, respectively. Then we have $(\varphi(\lambda'),\varphi(\mu'))^\dag=\frac{1}{r}(\lambda',\mu')'$ for any $\lambda',\mu'\in \h^{'*}$ (see Page 332 of \cite{W}).
\subsubsection{Characters of integrable highest weight modules} Let $A=(a_{i,j})_{i,j=1}^n$ be a generalized Cartan matrix,  $\n_+$ (resp. $\n_-$) be the subalgebra of $\g(A)$ generated by $e_1, \cdots, e_n$ (resp. $f_1, \cdots, f_n$). Then we have the triangular decomposition $$\g(A)=\n_+\oplus \h\oplus \n_-.$$

For any $\lambda\in \h^*$, we define the Verma module $V_{\g(A)}(\lambda)$ with highest weight $\lambda$ as follows:
$$V_{\g(A)}(\lambda)=U(\g(A))\otimes_{U(\n_+\oplus \h)}\C_{\lambda},$$
where $U(\g(A))$ denotes the universal enveloping algebra of $\g(A)$, and $\C_\lambda$ is the one dimensional module of  $\n_+\oplus \h$ such that $\n_+\cdot 1=0$, $h\cdot 1=\lambda(h)1$ for $h\in \h$. It is known \cite{K} that  $V_{\g(A)}(\lambda)$ has a unique maximal proper submodule $J(\lambda)$. Let  $L_{\g(A)}(\lambda)$ denote the corresponding irreducible quotient of  $V_{\g(A)}(\lambda)$.  It is known \cite{K} that $L_{\g(A)}(\lambda)$ has the following weight decomposition with respect to $\h$
$$L_{\g(A)}(\lambda)=\oplus_{\mu\in \h^*}L_{\g(A)}(\lambda)_\mu,$$
where $L_{\g(A)}(\lambda)_\mu=\{w\in L_{\g(A)}(\lambda)|h\cdot w=\mu(h) w, ~ \forall h\in \h\}$. The set
$$P(\lambda)=\{\mu\in \h^*|L_{\g(A)}(\lambda)_\mu\neq 0\}$$
 is called {\em the set of weights} of $L_{\g(A)}(\lambda)$.

For $\lambda\in \h^*$, we define a function $e^\lambda$ on $\h$ by $e^\lambda(h)=e^{\lambda(h)}$. Following \cite{K}, we define the {\em character } $ch_{L_{\g(A)}(\lambda)}$ of $L_{\g(A)}(\lambda)$ to be the function
$$h\mapsto ch_{L_{\g(A)}(\lambda)}(h)=\sum_{\mu\in \h^*}\dim L_{\g(A)}(\lambda)_\mu e^{\mu(h)}.$$

Given a non-negative integer $k$, let $P^+_k(A)$ denote the set of all dominant integral weights of level $k$ of $\g(A)$. For $\lambda\in P^+_k(A)$, let $\bar\lambda$ denote the projection of $\lambda$ on $\bar \h^*$. We then define $$\widetilde{P^+_k(A)}=\{\lambda\in P^+_k(A)|\lambda=k\Lambda_0+\bar\lambda\}.$$
For $\lambda\in P^+_k(A)$, the {\em normalized character} $\chi_{L_{\g(A)}(\lambda)}$ of $L_{\g(A)}(\lambda)$ is defined to be the function
$$h\mapsto \chi_{L_{\g(A)}(\lambda)}(h)=e^{-m_{\lambda}\delta(h)}\sum_{\mu\in \h^*}\dim L_{\g(A)}(\lambda)_\mu e^{\mu(h)},$$
where $m_\lambda=\frac{(\lambda+\rho, \lambda+\rho)}{2(k+h^\vee)}-\frac{(\rho,\rho)}{2h^\vee}$, and $\rho$ is the Weyl vector of $\g(A)$ determined by
\begin{align}\label{rho}
\rho(\alpha_i^\vee)=1 (i=1,\cdots, n) \text{ and }\rho(d)=0.
\end{align}
Following \cite{K}, we define the following coordinates in $\h^*\cong \h$:$$ h=(\tau, \mathfrak{z}, t)=2\pi i(-\tau\Lambda_0+\mathfrak{z}+t\delta).$$ Then the following important results have been established in \cite{KP}.
\begin{theorem}\label{KPmodular}
(1) Let $A$ be an affine Cartan matrix of type $X_{n}^{(1)}$,  and $\bar\Delta_+$ be the set of positive roots of $\overline{\g(A)}$. Then for any $\lambda\in \widetilde{P^+_k(A)}$,
$$\chi_{L_{\g(A)}(\lambda)}(\frac{1}{-\tau}, \frac{\mathfrak{z}}{\tau}, t-\frac{(\mathfrak{z},\mathfrak{z})}{2\tau} )=\sum_{\mu\in \widetilde{P^+_k(A)}}a_{\lambda, \mu}\chi_{L_{\g(A)}(\mu)}(\tau, \mathfrak{z}, t),$$
where  $$a_{\lambda, \mu}=i^{|\bar\Delta_+|}|M^{*}/(k+h^\vee)M|^{-\frac{1}{2}}\sum_{w\in \overline{W}}\epsilon(w)e^{-\frac{2\pi i}{k+h^\vee}(\bar \lambda+\bar \rho, w(\bar \mu+\bar \rho))}.$$
(2) Let $A^\dag$ be a twisted affine Cartan matrix not of type $A_{2l}^{(2)}$, $A'$ be its adjacent Cartan matrix, $\bar\Delta^\dag_+$ be the set of positive roots of $\overline{\g(A^\dag)}$ and $W^\dag$ be the Weyl group of $\g(A^\dag)$. Then for any $\lambda^\dag\in \widetilde{P^+_k(A^\dag)}$,
$$\chi_{L_{\g(A^\dag)}(\lambda^\dag)}(\frac{1}{-\tau}, \frac{\mathfrak{z}}{\tau}, t-\frac{(\mathfrak{z},\mathfrak{z})^\dag}{2\tau} )=\sum_{\lambda'\in \widetilde{P^+_k(A')}}a_{\lambda^\dag, \lambda'}\chi_{L_{\g(A')}(\lambda')}(\frac{\tau}{r}, \frac{\mathfrak{z}}{r}, t),$$
where  $r$ is the number such that $A^\dag$ belongs to Table Aff $r$ in \cite{K}, and
\begin{align}\label{amatrix}
a_{\lambda^\dag, \lambda'}=i^{|\bar\Delta^\dag_+|}|M^{\dag*}/(k+h^\vee)M^{\dag}|^{-\frac{1}{2}}|M'/M^\dag|^{\frac{1}{2}}\sum_{w\in \overline{ W}^\dag}\epsilon(w)e^{-\frac{2\pi i}{k+h^\vee}(w(\bar \lambda^\dag+\bar \rho^\dag), \varphi(\bar \lambda'+\bar \rho'))^\dag}.
\end{align}
\end{theorem}
\subsection{Definitions and properties about affine vertex operator algebras}\label{affineVOA} In this subsection, we recall some facts about affine vertex operator algebras.
Let $A$ be an affine Cartan matrix of type $X_n^{(1)}$, $\g(A)$ be the Kac-Moody algebra associated to $A$. Then it is proved in \cite{FZ} that $L_{\g(A)}(k\Lambda_0)$ has a vertex operator algebra structure. Let $\1$ denote a highest weight vector of $L_{\g(A)}(k\Lambda_0)$. We will use $x$ to denote the vector $x\cdot 1 $ for $x\in \overline{\g(A)}$. Let $\{u^i|1\leq i\leq \dim \g\}$ be an orthonormal basis of $\overline{\g(A)}$ with respect to $(,)$. Then  it is proved in \cite{FZ} that
\begin{align*}
\omega=\frac{1}{2(k+h^{\vee})} \sum_{i=1}^{\dim \g} u^i_{-1}u^i_{-1}\1
\end{align*}
 is a Virasoro vector of  $L_{\g(A)}(k\Lambda_0)$, where $h^\vee$ denotes the dual Coxeter number of $\g(A)$. Moreover, the following result has been proved in \cite{DLM2}, \cite{FZ}.
\begin{theorem}
Let $k$ be a positive integer. Then $L_{\g(A)}(k\Lambda_0)$  is a rational and $C_2$-cofinite  vertex operator algebra. Furthermore, $\{L_{\g(A)}(\lambda)|\lambda\in \widetilde{P^+_k(A)}\}$ is the complete set of irreducible modules of $L_{\g(A)}(k\Lambda_0)$.
\end{theorem}

By using the Verlinde formula, Kac and Walton obtained in \cite{Wa1, Wa2, K}  the following Kac-Walton formula, which gives fusion rules between untwisted modules of affine vertex operator algebras in terms of  Clebsch-Gordan coefficients associated to the corresponding finite dimensional
 simple Lie algebras (see also (5.1) of \cite{W1}):
\begin{theorem}
Let $k$ be a positive integer, $A$ be an affine Cartan matrix of type $X_n^{(1)}$. For any $\lambda_1, \lambda_2, \lambda_3\in \widetilde{P^+_{k}(A)}$, we have
\begin{align}\label{Kac-Walton}
N_{L_{\g(A)}(\lambda_1), L_{\g(A)}(\lambda_2)}^{ L_{\g(A)}( \lambda_3)}=\sum_{\substack{\mu\in P^+(\bar A),\\\mu+\bar\rho=w(\bar \lambda_3+\bar\rho)~({\rm mod}~ (k+h^\vee)M)\\(\text{ for some } w\in \overline{W})}}\epsilon(w)\mult_{\bar\lambda_1\otimes \bar\lambda_2}(\mu),
\end{align}
where $\mult_{\bar\lambda_1\otimes \bar\lambda_2}(\mu)=\dim \Hom_{\overline{\g(A)}}(L_{\overline{\g(A)}}(\bar\lambda_1)\otimes L_{\overline{\g(A)}}(\bar\lambda_2), L_{\overline{\g(A)}}(\mu))$ denotes the multiplicity of $L_{\overline{\g(A)}}(\mu)$ in the tensor product of $L_{\overline{\g(A)}}(\bar\lambda_1)\otimes L_{\overline{\g(A)}}(\bar\lambda_2)$.
\end{theorem}

We next recall from \cite{Li1} some  facts about automorphisms of affine vertex operator algebras. Let $\sigma$ be an automorphism of the Lie algebra $\overline{\g(A)}$. It is proved in \cite{Li1} that $\sigma$ induces an automorphism $\tilde\sigma$ of $L_{\g(A)}(k\Lambda_0)$. On the other hand,  automorphism groups of finite dimensional simple Lie algebras have been determined in \cite{K}.
\begin{theorem}
Let $\g$ be a finite dimensional simple Lie algebra,  $\h$ be a Cartan subalgebra of $\g$ and $\Pi=\{\alpha_1, \cdots, \alpha_n\}$ be a set of simple roots. Let $\sigma$ be an automorphism of $\g$ of order $T$. Then $\sigma$ is conjugate to an automorphism of $\g$  in the form
$\mu\exp\left(ad(\frac{2\pi i}{T}h)\right)$, $h\in \h^0$, where $\mu$ is a diagram automorphism preserving $\h$ and $\Pi$, $\h^0$ is the fixed point set of  $\mu$ in $\h$, and $\langle\alpha_i, h\rangle\in \Z (i=1,\cdots, n)$.
\end{theorem}

We now let $\sigma$ be a diagram automorphism of $\overline{\g(A)}$ of order $r$.  Then $\overline{\g(A)}$ has the following decomposition
 $$\overline{\g(A)}=\bigoplus_{j=0}^{r-1}\overline{\g(A)}_j,$$
 where $\overline{\g(A)}_j=\{x\in \overline{\g(A)}|\sigma(x)=e^{2\pi i\frac{j}{r}}x\}$. Consider the Lie algebra
 $$\hat{L}(\overline{\g(A)}, \sigma)=\bigoplus_{j=0}^{r-1} \overline{\g(A)}_j\otimes t^j\C[t^{r}, t^{-r}]\bigoplus \C K \bigoplus \C d$$ with the following Lie brackets
 \begin{align*}
&[x(m), y(n)]=[x, y](m+n)+( x, y) m\delta_{m+n,0}K,\\
&[K, \hat{L}(\overline{\g(A)}, \sigma)]=0, \ \ \ [d, x(m)]=mx(m),
\end{align*}
for $x, y\in \overline{\g(A)}$ and $m,n \in \Z$, where $x(n)$ denotes $x\otimes t^n$. Then it is proved in \cite{K} that $\hat{L}(\overline{\g(A)}, \sigma)$ is isomorphic to
 the twisted affine Kac-Moody algebra $\g(A^\dag)$ of $(\overline{\g(A)}, \sigma)$. We use $K^\dag$ and $d^\dag$ to denote the elements of $\g(A^\dag)$ defined by the formulas (\ref{k}) and (\ref{d}), respectively. Then it is known \cite{K} that $$K^\dag=rK\text{ and }d^\dag=d.$$
  The following result has been proved in Proposition 5.6 of \cite{Li1}.
\begin{proposition}
Let $k$ be a positive integer. Then $\{L_{\g(A^\dag)}(\lambda^\dag)|\lambda^\dag\in \widetilde{P^+_k(A^\dag)}\}$ is the complete set of irreducible $\tilde\sigma^{-1}$-twisted modules of $L_{\g(A)}(k\Lambda_0)$.
\end{proposition}
\pf Consider the Lie algebra
 $$\overline{\g(A)}[ \sigma]=\bigoplus_{j=0}^{r-1} \overline{\g(A)}_j\otimes t^\frac{j}{r}\C[t, t^{-1}]\bigoplus \C \tilde K \bigoplus \C \tilde d$$ with the following Lie brackets
 \begin{align*}
&[x\otimes t^m, y\otimes t^n]=[x, y]\otimes t^{m+n}+( x, y) m\delta_{m+n,0}\tilde K,\\
&[\tilde K, \overline{\g(A)}[ \sigma]]=0, \ \ \ [\tilde d, x\otimes t^n]=nx\otimes t^n,
\end{align*}
for $x, y\in \overline{\g(A)}$ and $m,n \in \frac{1}{r}\Z$. Let $M$ be a highest weight module of $[\overline{\g(A)}[ \sigma], \overline{\g(A)}[ \sigma]]$  of level $k$. It is proved in Proposition 5.6 of \cite{Li1} that $M$ is a  $\tilde\sigma^{-1}$-twisted module of $L_{\g(A)}(k\Lambda_0)$ if and only if $M$ is integrable. Moreover, the vertex operator is given by
$$Y_M(x,z)=\sum_{n\in \frac{1}{r}\Z}x\otimes t^nz^{-n-1}.$$

On the other hand, there is a Lie algebra isomorphism $$\Phi: \overline{\g(A)}[ \sigma]\to \hat{L}(\overline{\g(A)}, \sigma)$$ such that $\Phi(x\otimes t^n)=x(rn)$, $\Phi(\tilde{K})=rK$ and $\Phi(\tilde{d})=\frac{d}{r}$. Therefore, a highest weight $[\g(A^\dag), \g(A^\dag)]$-module $L_{\g(A^\dag)}(\lambda^\dag)$ of level $k$ is a  $\tilde\sigma^{-1}$-twisted module of $L_{\g(A)}(k\Lambda_0)$ if and only if $L_{\g(A^\dag)}(\lambda^\dag)$ is integrable. Moreover, the vertex operator is given by
$$Y_{L_{\g(A^\dag)}(\lambda^\dag)}(x,z)=\sum_{n\in \frac{1}{r}\Z}x(r n)z^{-n-1}.$$
This implies that $\{L_{\g(A^\dag)}(\lambda^\dag)|\lambda^\dag\in \widetilde{P^+_k(A^\dag)}\}$ is the complete set of irreducible $\tilde\sigma^{-1}$-twisted modules of $L_{\g(A)}(k\Lambda_0)$.\qed
\begin{remark}
Note that the definition of twisted module is not same as that in \cite{Li1}. Following \cite{DLM3}, we interchange the notions of $g$-twisted modules and $g^{-1}$-twisted modules.
\end{remark}

We next determine the action of $L(0)$ on twisted modules of $L_{\g(A)}(k\Lambda_0)$. Let $\h^\dag$ and $\bar\h$ be the Cartan subalgebras of $\g(A^\dag)$ and $\overline{\g(A)}$, respectively. Then we have $\h^\dag=\bar\h^0\oplus \C K\oplus \C d$, where $\bar\h^0$ denotes the fixed point set of $\bar\h$ under the action of $\sigma$. Let $(,)$ and $(,)^\dag$ denote the normalized bilinear forms of $\g(A)$ and $\g(A^\dag)$, respectively. Then it is known \cite{K} that
\begin{align}\label{key6}
(\lambda_1^\dag, \lambda_2^\dag)=\frac{1}{r}(\lambda_1^\dag, \lambda_2^\dag)^\dag
\end{align}
for any $\lambda_1^\dag, \lambda_2^\dag\in \h^{\dag*}$.
Then we have the following result.
\begin{proposition}\label{key3}
Let $k$ be a positive integer, and $\g(A^\dag)$ be the twisted affine Kac-Moody algebra as above. For any $\lambda^\dag\in \widetilde{P^+_k(A^\dag)}$, the operators $L(0)-\frac{c}{24}\id$ and $\frac{m_{\lambda^\dag}}{r}\id-\frac{d^\dag}{r}$ on $L_{\g(A^\dag)}(\lambda^\dag)$ are equal, where $c$ denotes the central charge of $L_{\g(A)}(k\Lambda_0)$.
\end{proposition}
\pf Let $\rho^\dag$ be the Weyl vector of $\g(A^\dag)$ defined by the formula (\ref{rho}). Then it is proved in Lemmas 5.3, 3.3 of \cite{KFP} that $$L(0)+\tilde d=\frac{1}{k+h^\vee}\left(\frac{1}{2}(\bar\lambda^\dag+2\bar \rho^\dag, \bar\lambda^\dag)+kz(\overline{\g(A)}, \sigma)\right)\id,$$
where $z(\overline{\g(A)}, \sigma)=\frac{1}{4}\sum_{j=0}^{r-1}\frac{j}{r}(1-\frac{j}{r})\dim \overline{\g(A)}_j$.

Recall that $m_{\lambda^\dag}=\frac{(\lambda^\dag+\rho^\dag, \lambda^\dag+\rho^\dag)^\dag}{2(k+h^\vee)}-\frac{(\rho^\dag,\rho^\dag)^\dag}{2h^\vee}$. Moreover, it is proved in Proposition 6.1 of \cite{KFP} that
$$\frac{(\bar\rho^\dag,\bar\rho^\dag)}{2h^\vee}=\frac{\dim\overline{\g(A)}}{24}-z(\overline{\g(A)}, \sigma).$$Then we have
\begin{align*}
m_{\lambda^\dag}&=\frac{(\lambda^\dag+\rho^\dag,\lambda^\dag+\rho^\dag)^\dag}{2(k+h^\vee)}-\frac{(\rho^\dag,\rho^\dag)^\dag}{2h^\vee}\\
&=\frac{(\bar\lambda^\dag+\bar\rho^\dag,\bar\lambda^\dag+\bar\rho^\dag)^\dag}{2(k+h^\vee)}-\frac{(\bar\rho^\dag,\bar\rho^\dag)^\dag}{2h^\vee}\\
&=\frac{(\bar\lambda^\dag+2\bar\rho^\dag,\bar\lambda^\dag)^\dag}{2(k+h^\vee)}-\frac{k(\bar\rho^\dag,\bar\rho^\dag)^\dag}{2h^\vee(k+h^\vee)}\\
&=\frac{r(\bar\lambda^\dag+2\bar\rho^\dag,\bar\lambda^\dag)}{2(k+h^\vee)}-\frac{rk}{(k+h^\vee)} \left(\frac{\dim\overline{\g(A)}}{24}-z(\overline{\g(A)}, \sigma)\right)\\
&=\frac{r(\bar\lambda^\dag+2\bar\rho^\dag,\bar\lambda^\dag)}{2(k+h^\vee)}+\frac{rkz(\overline{\g(A)}, \sigma)}{(k+h^\vee)}- \frac{rc}{24}.
\end{align*}
Since $\Phi(\tilde{d})=\frac{d}{r}=\frac{d^\dag}{r}$, this implies that the operators $L(0)-\frac{c}{24}\id$ and $\frac{m_{\lambda^\dag}}{r}\id-\frac{d^\dag}{r}$ on $L_{\g(A^\dag)}(\lambda^\dag)$ are equal.
\qed
\subsection{$S$-matrix in the orbifold theory of  affine vertex operator algebras} In this subsection, we determine the $S$-matrix in the orbifold theory of  affine vertex operator algebras, which is a key step in determining fusion rules between twisted modules of affine vertex operator algebras.  The key point in our proof is to use the theory of orbit Lie algebras established in \cite{FSS}, \cite{FRS}.
\subsubsection{Orbit Lie algebras} We recall from \cite{FRS} some facts about orbit Lie algebras. Let $A=(a_{i,j})_{i,j\in I}$ be a symmetric affine Cartan matrix, $\dot{\sigma}: I\to I$ be a bijection of finite order which keeps the Cartan matrix $A$ fixed, i.e., $a_{\dot{\sigma}(i),\dot{\sigma}(j)}=a_{i,j}$ for all $i,j\in I$. Let $\g(A)$ denote the Kac-Moody algebra associated to $A$. Then $\dot{\sigma}$ induces an automorphism $\sigma$ of $\g(A)$. Let $\h$ denote the Cartan subalgebra of $\g(A)$. Then  $\sigma$ preserves the Cartan subalgebra. We use $\h^0$ to denote the fixed point set of $\h$ under the action of $\sigma$.

Let $r$ be the order of $\dot{\sigma}$ and $N_i$ be the length of the $\dot{\sigma}$-orbit of $i$ in $I$. Following \cite{FRS}, we define the following subsets of $I$:
\begin{align*}
&\hat I=\{i\in I|i\leq \dot{\sigma}^l(i), \forall l\in \Z\},\\
&\check{I}=\{i\in \hat I|\sum_{l=0}^{N_i-1}a_{i, \dot{\sigma}^l(i)}>0 \}.
\end{align*}
For $i\in \hat I$, define
\begin{align*}
s_i=\left\{\begin{array}{ll}
\frac{a_{ii}}{\sum_{l=0}^{N_i-1}a_{i, \dot{\sigma}^l(i)}}&\text{ if } i\in \check I;\\
1&\text{otherwise}.
\end{array}\right.
\end{align*}
Following \cite{FRS}, we define the matrix $\hat A=(\hat a_{ij})_{i,j\in \hat I}$ as follows:
$$\hat a_{ij}=s_j\sum_{l=0}^{N_j-1}a_{i, \dot{\sigma}^l(j)}.$$
Then it is proved in Lemma 2.1 of \cite{FRS} that $\hat A$ is a generalized Cartan matrix. Let $\g(\hat A)$ be the Kac-Moody algebra associated to $\hat A$, $\hat{\h}$ be the Cartan subalgebra of $\g(\hat A)$. Let $\hat e_i, \hat f_i (i\in \hat I)$ denote the Chevalley generators of $\g(\hat A)$.  Then the {\em orbit Lie algebra} associated to $\sigma$ of $\g(A)$ is defined to be the Lie subalgebra of $\g(\hat A)$ generated by $\hat e_i, \hat f_i (i\in \check I)$ and $\hat{\h}$.

Following Subsection 2.5 of \cite{FSS}, we define the following maps for affine Cartan matrices. For $A=A_{2n-1}^{(1)}$, we consider the map $\dot{\sigma}: i\mapsto 2n-i \text{ mod } 2n$. For $A=D_{n+1}^{(1)}$, we consider the map $\dot{\sigma}: n+1\mapsto n, n\mapsto n+1, \text{ and } i\mapsto i$ else. For $A=D_4^{(1)}$, we consider the map $\dot{\sigma}: 0\mapsto 0, 1\mapsto 3\mapsto 4\mapsto 1, 2\mapsto 2$. For $A=E_6^{(1)}$, we consider the map $\dot{\sigma}: 1\mapsto 5, 2 \mapsto 4, 3\mapsto 3, 6\mapsto 6, 0\mapsto 0$. Then the following results have been established in  \cite{FSS}.
\begin{proposition}\label{orbitlie}
Let $A$ be an affine Cartan matrix of type $A_{2n-1}^{(1)}$, $D_{n+1}^{(1)}$, $D_{4}^{(1)}$ or $E_6^{(1)}$, $\dot{\sigma}$ be the map defined as above. Then\\
 (1) The orbit Lie algebra associated to $\sigma$ of $\g(A)$ is equal to $\g(\hat A)$. \\
 (2) The orbit Lie algebras associated to $\sigma$ of $\g(A_{2n-1}^{(1)}), \g(D_{n+1}^{(1)}), \g(D_{4}^{(1)}), \g(E_6^{(1)})$ are isomorphic to $ \g(D_{n+1}^{(2)}), \g(A_{2n-1}^{(2)}), \g(D_{4}^{(3)}), \g(E_6^{(2)})$, respectively.
\end{proposition}

Let $A$ be an affine Cartan matrix of type $A_{2n-1}^{(1)}$, $D_{n+1}^{(1)}$, $D_{4}^{(1)}$ or $E_6^{(1)}$, $A^\dag$ be an affine Cartan matrix of type $A_{2n-1}^{(2)}$, $D_{n+1}^{(2)}$, $D_{4}^{(3)}$ or $E_6^{(2)}$, $A'$ be the adjacent Cartan matrix of $A^\dag$. By Proposition \ref{orbitlie}, the orbit Lie algebra $\g(\hat A)$ associated to $\sigma$ of $\g(A)$ is isomorphic to $\g(A')$. Let $\h^0$ be the fixed point set of $\h$ under the action of $\sigma$. Following \cite{FRS}, we  define a linear map $P_\sigma: \h^0\to \h'$ such that $P_\sigma(\sum_{l=0}^{N_i-1}\alpha_{\dot{\sigma}^l(i)}^\vee)=N_i\alpha_i'^{\vee}$
and $P_\sigma(d)=d'$. Then we have the following results.
\begin{lemma}\label{orbitmap}
(1) $P_\sigma$ is a linear isomorphism.\\
(2) For any $h_1, h_2\in \h^0$, $(P_\sigma(h_1),P_\sigma( h_2))'=(h_1, h_2)$, where $(,)$ and $(,)'$ denote the standard invariant bilinear forms of $\g(A)$ and $\g(A')$, respectively.\\
(3) $P_\sigma(K)=K'$, where $K$ and $K'$ denote the canonical central elements of $\g(A)$ and $\g(A')$, respectively.
\end{lemma}
\pf (1) Note that $\dim \h^0=\dim \h'$ and $P_\sigma$ is surjective. Then $P_\sigma$ is a linear isomorphism.

(2) By Lemma 2.2 of \cite{FRS}, $(P_\sigma(h_1),P_\sigma( h_2))'=(h_1, h_2)$ holds for $h_1, h_2\in \h^0\cap [\g(A), \g(A)]$. It is enough to prove that
$(P_\sigma(d),P_\sigma( h_2))'=(d, h_2)$ holds for any $h_2\in \h^0$. Note that $(P_\sigma(d),P_\sigma(d))'=(d',d')'=0=(d, d)$ and $(P_\sigma(d),P_\sigma( \alpha_0^\vee))'=(d', \alpha_0'^{ \vee})'=1=(d, \alpha_0^\vee)$. Finally, $(P_\sigma(d),P_\sigma(\sum_{l=0}^{N_i-1}\alpha_{\dot{\sigma}^l(i)}^\vee))'=N_i(d', \alpha_i'^{\vee})'=0=(d, \sum_{l=0}^{N_i-1}\alpha_{\dot{\sigma}^l(i)}^\vee)$ for $i\neq 0$. Thus,  $(P_\sigma(h_1),P_\sigma( h_2))'=(h_1, h_2)$ holds for any $h_1, h_2\in \h^0$.

(3) Note that $K'$ is determined by $(K', \alpha'^{ \vee}_i)'=0$ and $(K', d')'=1$. On the other hand, $N_i(P_\sigma(K), \alpha'^{ \vee}_i)'=(K, \sum_{l=0}^{N_i-1}\alpha_{\dot{\sigma}^l(i)}^\vee)=0$ and $(P_\sigma(K), d')'=(K, d)=1$. Thus, $P_\sigma(K)=K'$.
\qed

\vskip.25cm
Define a  linear map $\iota^*: \h^*\to (\h^0)^*, \lambda\mapsto \lambda|_{\h^0}$. Since $\sigma: \h\to \h$ is a linear isomorphism,  the dual map $\sigma^*: \h^*\to \h^*$ is also a linear isomorphism. Moreover, we have $\sigma^*(\lambda)(h)=\lambda(\sigma(h))$ and $\sigma^*(\alpha_i)=\alpha_{\dot\sigma^{-1}(i)}$ for any $\lambda\in \h^*$ and $h\in \h$ (see Page 525 of \cite{FRS}). Let $(\h^*)^0$ denote the fixed point set of $\h^*$ under the action of $\sigma^*$. The elements in $(\h^*)^0$ are called {\em symmetric weights}. Note that $(,)|_{\h^0}$ is nondegenerate, then $\iota^*: (\h^*)^0\to (\h^0)^*, \lambda\mapsto \lambda|_{\h^0}$ is a linear isomorphism (see Page 525 of \cite{FRS}).

By Lemma \ref{orbitmap}, we have a linear isomorphism $P_\sigma: \h^0\to \h'$. Therefore, we have a linear isomorphism $\tilde{P}_\sigma: (\h')^*\to (\h^0)^*$ such that $\tilde{P}_\sigma^*(\lambda)(h)=\lambda(P_\sigma(h))$ for any $\lambda\in (\h')^*$ and $h\in \h^0$. Let $P_\sigma^*=(\iota^*)^{-1}\circ\tilde{P}_\sigma^*$. Then $P_\sigma^*: (\h')^*\to (\h^*)^0$ is a linear isomorphism. Moreover, we have the following results.
\begin{lemma}\label{orbitdual}
(1) $P_\sigma^*(\Lambda_0')=\Lambda_0$ and $P_\sigma^*(\delta')=\delta$.\\
(2) For any $\lambda_1, \lambda_2\in (\h')^*$,  $(P_\sigma^*(\lambda_1), P_\sigma^*(\lambda_2))=(\lambda_1, \lambda_2)'$.\\
(3) $P_\sigma^*(\rho')=\rho$ and $P_\sigma^*(\bar\rho')=\bar\rho$.
\end{lemma}
\pf (1) $\Lambda_0$ is determined by
$\langle \Lambda_0, \alpha_i^\vee\rangle=\delta_{0,i}$ and $\langle\Lambda_0, d\rangle=0$. Note that $$\langle P_\sigma^*(\Lambda_0'), \alpha_0^\vee\rangle=\langle (\iota^*)^{-1}\circ\tilde{P}_\sigma^*(\Lambda_0'), \alpha_0^\vee\rangle=\langle \tilde{P}_\sigma^*(\Lambda_0'), \alpha_0^\vee\rangle=\langle \Lambda_0', P_\sigma(\alpha_0^\vee)\rangle=\langle \Lambda_0', \alpha_0'^{ \vee}\rangle=1$$ and
\begin{align*}
\langle P_\sigma^*(\Lambda_0'), \alpha_i^\vee\rangle&=\frac{1}{N_i}\langle P_\sigma^*(\Lambda_0'), \sum_{l=0}^{N_i-1}\alpha_{\dot{\sigma}^l(i)}^\vee\rangle=\frac{1}{N_i}\langle (\iota^*)^{-1}\circ\tilde{P}_\sigma^*(\Lambda_0'), \sum_{l=0}^{N_i-1}\alpha_{\dot{\sigma}^l(i)}^\vee\rangle\\
&=\frac{1}{N_i}\langle \tilde{P}_\sigma^*(\Lambda_0'), \sum_{l=0}^{N_i-1}\alpha_{\dot{\sigma}^l(i)}^\vee\rangle=\frac{1}{N_i}\langle \Lambda_0', P_\sigma(\sum_{l=0}^{N_i-1}\alpha_{\dot{\sigma}^l(i)}^\vee)\rangle\\
&=\langle \Lambda_0', \alpha_i'^{\vee}\rangle=0.
\end{align*}
Moreover, $$\langle P_\sigma^*(\Lambda_0'), d\rangle=\langle (\iota^*)^{-1}\circ\tilde{P}_\sigma^*(\Lambda_0'), d\rangle=\langle \tilde{P}_\sigma^*(\Lambda_0'), d\rangle=\langle \Lambda_0', P_\sigma(d)\rangle=\langle \Lambda_0', d'\rangle=0.$$
Therefore, $P_\sigma^*(\Lambda_0')=\Lambda_0$.

We next prove that $P_\sigma^*(\delta')=\delta$. Note that $\delta$ is determined by $\langle \delta, \alpha_i^\vee\rangle=0$ and $\langle\delta, d\rangle=1$. By the similar arguments as above, we have $\langle P_\sigma^*(\delta'), \alpha_i^\vee\rangle=\langle \delta', \alpha_i'^{\vee}\rangle=0$ and $\langle P_\sigma^*(\delta'), d\rangle=\langle \delta', d'\rangle=1$. Therefore, $P_\sigma^*(\delta')=\delta$.

(2) Since $P_\sigma^*$ is a linear isomorphism, it follows from  Lemma 2.3 of \cite{FRS} that   $(P_\sigma^*(\lambda_1), P_\sigma^*(\lambda_2))=(\lambda_1, \lambda_2)'$ holds for any $\lambda_1, \lambda_2\in (\h')^*$.

(3) By Lemma 3.2 of \cite{FRS}, $P_\sigma^*(\rho')=\rho$. On the other hand, by (6.2.8) of \cite{K}, $\rho'=\bar \rho'+h^\vee\Lambda_0'$ and $\rho=\bar \rho+h^\vee\Lambda_0$, it follows from (1) that $P_\sigma^*(\bar\rho')=\bar\rho$.
\qed

\vskip.25cm
As a corollary, we have
\begin{corollary}\label{key5}
For any $\lambda'\in \widetilde{P^+_k(A')}$, we have $m_{\lambda'}=m_{P_\sigma^*(\lambda')}$.
\end{corollary}
\pf Recall that $m_{\lambda'}=\frac{(\lambda'+\rho', \lambda'+\rho')'}{2(k+h^\vee)}-\frac{(\rho',\rho')'}{2h^\vee}$,  where $\rho'$ is the Weyl vector of $\g(A')$. By Lemma \ref{orbitdual}, we have
\begin{align*}
m_{\lambda'}&=\frac{(P_\sigma^*(\lambda'+\rho'), P_\sigma^*(\lambda'+\rho'))}{2(k+h^\vee)}-\frac{(P_\sigma^*(\rho'),P_\sigma^*(\rho'))}{2h^\vee}\\
&=\frac{(P_\sigma^*(\lambda')+\rho, P_\sigma^*(\lambda')+\rho)}{2(k+h^\vee)}-\frac{(\rho,\rho)}{2h^\vee}=m_{P_\sigma^*(\lambda')}.
\end{align*}
This completes the proof.
\qed

\vskip.25cm
Set $(\widetilde{P^+_k(A)})^\sigma=\widetilde{P^+_k(A)}\cap (\h^*)^0$, then we have the following result.
\begin{proposition}\label{key4}
For any $k\in \Z_{\geq 0}$, $P_\sigma^*(\widetilde{P^+_k(A')})=(\widetilde{P^+_k(A)})^\sigma$.
\end{proposition}
\pf We first prove that $P_\sigma^*(\widetilde{P^+_k(A')})\subset(\widetilde{P^+_k(A)})^\sigma$. By the similar arguments as in the proof of Lemma \ref{orbitdual},  we have $\langle P_\sigma^*(\lambda'), \alpha_i^\vee\rangle=\langle\lambda',\alpha_i'{^\vee}\rangle\in \Z_{\geq 0}$ for any $\lambda'\in \widetilde{P^+_k(A')}$. Then it follows from Lemmas \ref{orbitmap}, \ref{orbitdual} that $P_\sigma^*(\widetilde{P^+_k(A')})\subset(\widetilde{P^+_k(A)})^\sigma$.

We next prove that $P_\sigma^*(\widetilde{P^+_k(A')})\supset(\widetilde{P^+_k(A)})^\sigma$. For any $\lambda\in (\widetilde{P^+_k(A)})^\sigma$, we have
\begin{align*}
\langle P_\sigma^{*-1}(\lambda), \alpha_i'^{\vee}\rangle&=\langle P_\sigma^{*-1}(\lambda), \frac{1}{N_i}P_\sigma(\sum_{l=0}^{N_i-1}\alpha_{\dot{\sigma}^l(i)}^\vee)\rangle
=\langle \tilde{P}_\sigma^{*-1}\circ\iota^*(\lambda), \frac{1}{N_i}P_\sigma(\sum_{l=0}^{N_i-1}\alpha_{\dot{\sigma}^l(i)}^\vee)\rangle\\
&=\langle \iota^*(\lambda), \frac{1}{N_i}\sum_{l=0}^{N_i-1}\alpha_{\dot{\sigma}^l(i)}^\vee\rangle=\langle \lambda, \frac{1}{N_i}\sum_{l=0}^{N_i-1}\alpha_{\dot{\sigma}^l(i)}^\vee\rangle
=\langle \lambda, \alpha_{i}^\vee\rangle\in \Z_{\geq 0}.
\end{align*}
Then it follows from Lemmas \ref{orbitmap}, \ref{orbitdual} that $\widetilde{P^+_k(A')}\supset P_\sigma^{*-1}((\widetilde{P^+_k(A)})^\sigma)$. Therefore, $P_\sigma^*(\widetilde{P^+_k(A')})\supset(\widetilde{P^+_k(A)})^\sigma$.
\qed
\vskip.25cm

For any $\lambda\in (\widetilde{P^+_k(A)})^\sigma$, it is proved in \cite{FRS} that there is a linear isomorphism $\phi(\sigma):L_{\g(A)}(\lambda) \to L_{\g(A)}(\lambda)$ such that $$\phi(\sigma)x\phi(\sigma)^{-1}=\sigma(x)$$ for any $x\in \g(A)$. Then the following result has been established in Theorem 3.1 of \cite{FRS}.
\begin{theorem}\label{key2}
For any $\mu'\in P(L_{\g(A')}(\lambda'))$, $\dim L_{\g(A')}(\lambda')_{\mu'}={\rm tr}_{L_{\g(A)}(P_\sigma^*(\lambda'))_{P_\sigma^*(\mu')}}\phi(\sigma)$.
\end{theorem}

Let $\h^\dag$ denote the Cartan subalgebra of $\g(A^\dag)$. By the results in Subsection 8.3 of \cite{K}, there is a linear isomorphism $\iota: \h^\dag\to \h^0$ such that
\begin{align}\label{iota}
&\iota(K^\dag)=rK,\ \ \ \iota(d^\dag)=d,\notag\\
&\iota(\alpha_i^{\dag\vee})=\sum_{j=0}^{N_i-1}\alpha_{\dot\sigma^l(i)}^\vee,\text{ \ \ \  if } A=A_{2l-1}^{(2)}, D_{l+1}^{(2)}\text{ and } 1\leq i\leq l;\notag\\
&\iota(\alpha_1^{\dag\vee})=\alpha_1^\vee+\alpha_5^\vee, \iota(\alpha_2^{\dag\vee})=\alpha_2^\vee+\alpha_4^\vee, \iota(\alpha_3^{\dag\vee})=\alpha_3^\vee, \iota(\alpha_4^{\dag\vee})=\alpha_6^\vee, \text{\ \ \ if } A=E_{6}^{(2)};\notag\\
&\iota(\alpha_1^{\dag\vee})=\alpha_1^\vee+\alpha_3^\vee+\alpha_4^\vee, \iota(\alpha_2^{\dag\vee})=\alpha_2^\vee, \text{\ \ \ if } A= D_{4}^{(3)}.
\end{align}
Recall that we have linear isomorphisms $\nu:\h\to \h^*$ and $\nu^\dag:\h^\dag\to \h^{\dag*}$. Then we have the following results.
\begin{proposition}\label{key1}
(1) $\nu:\h^0\to (\h^*)^0$ is a linear isomorphism.\\
(2) $\nu\circ \iota\circ \nu^{\dag-1}=P_{\sigma}^*\circ\varphi^{-1}$, where $\nu^{\dag-1}$ denotes $(\nu^{\dag})^{-1}$ and $\varphi: \h^{'*}\to \h^{\dag*}$ is the linear isomorphism defined by the formula (\ref{varphi}).
\end{proposition}
\pf (1) For any $h\in \h^0$, we have $\sigma(h)=h$. Then for any $h_1\in \h$,
\begin{align*}
\langle \sigma^*(\nu(h)), h_1\rangle&=\langle \nu(h), \sigma(h_1)\rangle\\
&=(h, \sigma(h_1))=(\sigma^{-1}(h), h_1)=(h, h_1)\\
&=\langle\nu(h), h_1\rangle.
\end{align*}
Therefore, $\sigma^*(\nu(h))=\nu(h)$. Hence, $\nu(h)\in (\h^*)^0$. Note that $\dim \h^0=\dim (\h^*)^0$ (see Page 525 of \cite{FRS}). It follows that $\nu:\h^0\to (\h^*)^0$ is a linear isomorphism.

(2) We first consider the case that $A^\dag=A_{2l-1}^{(2)}$ or $ D_{l+1}^{(2)}$. By the formula (\ref{iota}), we have
\begin{align*}
&\nu\circ \iota\circ \nu^{\dag-1}(\alpha_i^\dag)=\nu\circ \iota(\frac{a_i^{\dag\vee}}{a_i^{\dag}}\alpha_i^{\dag\vee})=\frac{a_i^{\dag\vee}}{a_i^{\dag}}\nu(\sum_{j=0}^{N_i-1}\alpha_{\dot\sigma^l(i)}^\vee)=\frac{a_i^{\dag\vee}}{a_i^{\dag}}\sum_{j=0}^{N_i-1}\alpha_{\dot\sigma^l(i)},\\
&\nu\circ \iota\circ \nu^{\dag-1}(\delta^\dag)=\nu\circ \iota(K^\dag)=\nu(2K)=2\delta,\\
&\nu\circ \iota\circ \nu^{\dag-1}(\Lambda_0^\dag)=\nu\circ \iota(d^\dag)=\nu(d)=\Lambda_0.
\end{align*}
On the other hand, by the formula (\ref{varphi})and Lemma 2.3 of \cite{FRS},
\begin{align*}
&P_{\sigma}^*\circ\varphi^{-1}(\alpha_i^\dag)=P_{\sigma}^*(\frac{a_i^{\dag\vee}}{a_i^{\dag}}\alpha_i')=\frac{a_i^{\dag\vee}}{a_i^{\dag}}\sum_{j=0}^{N_i-1}\alpha_{\dot\sigma^l(i)},\\
&P_{\sigma}^*\circ\varphi^{-1}(\delta^\dag)=P_{\sigma}^*(2\delta')=2\delta,\ \ P_{\sigma}^*\circ\varphi^{-1}(\Lambda_0^\dag)=P_{\sigma}^*(\Lambda_0')=\Lambda_0.
\end{align*}
Therefore, $\nu\circ \iota\circ \nu^{\dag-1}=P_{\sigma}^*\circ\varphi^{-1}$ holds for $A^\dag=A_{2l-1}^{(2)}$ or $ D_{l+1}^{(2)}$.

We next consider  the case that $A^\dag=E_{6}^{(2)}$. By the formula (\ref{iota}), we have
\begin{align*}
&\nu\circ \iota\circ \nu^{\dag-1}(\alpha_1^\dag)=\nu\circ \iota(\frac{a_1^{\dag\vee}}{a_1^{\dag}}\alpha_1^{\dag\vee})=\frac{a_1^{\dag\vee}}{a_1^{\dag}}(\alpha_{1}+\alpha_5),\\
&\nu\circ \iota\circ \nu^{\dag-1}(\alpha_2^\dag)=\nu\circ \iota(\frac{a_2^{\dag\vee}}{a_2^{\dag}}\alpha_2^{\dag\vee})=\frac{a_2^{\dag\vee}}{a_2^{\dag}}(\alpha_{2}+\alpha_4),\\
&\nu\circ \iota\circ \nu^{\dag-1}(\alpha_3^\dag)=\nu\circ \iota(\frac{a_3^{\dag\vee}}{a_3^{\dag}}\alpha_3^{\dag\vee})=\frac{a_3^{\dag\vee}}{a_3^{\dag}}\alpha_{3},\\
&\nu\circ \iota\circ \nu^{\dag-1}(\alpha_4^\dag)=\nu\circ \iota(\frac{a_4^{\dag\vee}}{a_4^{\dag}}\alpha_4^{\dag\vee})=\frac{a_4^{\dag\vee}}{a_4^{\dag}}\alpha_{6},\\
&\nu\circ \iota\circ \nu^{\dag-1}(\delta^\dag)=\nu\circ \iota(K^\dag)=\nu(2K)=2\delta,\\
&\nu\circ \iota\circ \nu^{\dag-1}(\Lambda_0^\dag)=\nu\circ \iota(d^\dag)=\nu(d)=\Lambda_0.
\end{align*}
On the other hand, by the formula (\ref{varphi})and Lemma 2.3 of \cite{FRS},
\begin{align*}
&P_{\sigma}^*\circ\varphi^{-1}(\alpha_1^\dag)=P_{\sigma}^*(\frac{a_1^{\dag\vee}}{a_1^{\dag}}\alpha_4')=\frac{a_1^{\dag\vee}}{a_1^{\dag}}(\alpha_{1}+\alpha_5),\\
&P_{\sigma}^*\circ\varphi^{-1}(\alpha_2^\dag)=P_{\sigma}^*(\frac{a_2^{\dag\vee}}{a_2^{\dag}}\alpha_3')=\frac{a_2^{\dag\vee}}{a_2^{\dag}}(\alpha_{2}+\alpha_4),\\
&P_{\sigma}^*\circ\varphi^{-1}(\alpha_3^\dag)=P_{\sigma}^*(\frac{a_3^{\dag\vee}}{a_3^{\dag}}\alpha_2')=\frac{a_3^{\dag\vee}}{a_3^{\dag}}\alpha_3,
\ \ P_{\sigma}^*\circ\varphi^{-1}(\alpha_4^\dag)=P_{\sigma}^*(\frac{a_4^{\dag\vee}}{a_4^{\dag}}\alpha_1')=\frac{a_4^{\dag\vee}}{a_4^{\dag}}\alpha_{6},\\
&P_{\sigma}^*\circ\varphi^{-1}(\delta^\dag)=P_{\sigma}^*(2\delta')=2\delta,\ \ P_{\sigma}^*\circ\varphi^{-1}(\Lambda_0^\dag)=P_{\sigma}^*(\Lambda_0')=\Lambda_0.
\end{align*}
Therefore, $\nu\circ \iota\circ \nu^{\dag-1}=P_{\sigma}^*\circ\varphi^{-1}$ holds for $A^\dag=E_{6}^{(2)}$.

We next consider  the case that $A^\dag= D_{4}^{(3)}$.  By the formula (\ref{iota}), we have
\begin{align*}
&\nu\circ \iota\circ \nu^{\dag-1}(\alpha_1^\dag)=\nu\circ \iota(\frac{a_1^{\dag\vee}}{a_1^{\dag}}\alpha_1^{\dag\vee})=\frac{a_1^{\dag\vee}}{a_1^{\dag}}(\alpha_{1}+\alpha_3+\alpha_4),\\
&\nu\circ \iota\circ \nu^{\dag-1}(\alpha_2^\dag)=\nu\circ \iota(\frac{a_2^{\dag\vee}}{a_2^{\dag}}\alpha_2^{\dag\vee})=\frac{a_2^{\dag\vee}}{a_2^{\dag}}\alpha_{2},\\
&\nu\circ \iota\circ \nu^{\dag-1}(\delta^\dag)=\nu\circ \iota(K^\dag)=\nu(3K)=3\delta,\\
&\nu\circ \iota\circ \nu^{\dag-1}(\Lambda_0^\dag)=\nu\circ \iota(d^\dag)=\nu(d)=\Lambda_0.
\end{align*}
On the other hand, by the formula (\ref{varphi})and Lemma 2.3 of \cite{FRS},
\begin{align*}
&P_{\sigma}^*\circ\varphi^{-1}(\alpha_1^\dag)=P_{\sigma}^*(\frac{a_1^{\dag\vee}}{a_1^{\dag}}\alpha_2')=\frac{a_1^{\dag\vee}}{a_1^{\dag}}(\alpha_{1}+\alpha_3+\alpha_4),\\
&P_{\sigma}^*\circ\varphi^{-1}(\alpha_2^\dag)=P_{\sigma}^*(\frac{a_2^{\dag\vee}}{a_2^{\dag}}\alpha_1')=\frac{a_2^{\dag\vee}}{a_2^{\dag}}\alpha_{2},\\
&P_{\sigma}^*\circ\varphi^{-1}(\delta^\dag)=P_{\sigma}^*(3\delta')=3\delta,\ \ P_{\sigma}^*\circ\varphi^{-1}(\Lambda_0^\dag)=P_{\sigma}^*(\Lambda_0')=\Lambda_0.
\end{align*}
Therefore, $\nu\circ \iota\circ \nu^{\dag-1}=P_{\sigma}^*\circ\varphi^{-1}$ holds for $A^\dag=D_{4}^{(3)}$.
\qed
\subsubsection{Li's operators}\label{Li}
 Let $V$ be a vertex operator algebra, $h\in V_1$ be an element satisfying the following conditions:
$$L(n)h=\delta_{n,0}h,\ \ \ h_nh=\delta_{n,1}\langle h, h\rangle\1 \text{ for } n\in \Z_{\geq 0}.$$
Suppose further that $V$ is finitely generated and $h_0$ acts semisimply on $V$ with ${\rm Spec}~ h_0\subset \frac{1}{T}\Z$ for some positive integer $T$. Define $$\varphi(h)=e^{-2\pi ih_0}.$$ Then $\varphi(h)$ is an automorphism of $V$ such that $\varphi(h)^T=1$.

 Set
$$\Delta(h, z)=z^{h_0}\exp\left(\sum_{n=1}^{\infty}\frac{h_n(-z)^{-n}}{-n}\right).$$
Then the following result has been established in Proposition 5.4 of \cite{Li1}.
\begin{proposition}\label{inner}
Let $(M,Y_{M}(\cdot, z))$ be an irreducible $V$-module and $h$ be as above. Set
$$(M^{(h)}, Y_{M^{(h)}}(\cdot, z))=(M, Y_M(\Delta(h, z)\cdot, z)).$$ Then $(M^{(h)}, Y_{M^{(h)}}(\cdot, z))$ is an irreducible $\varphi(h)$-twisted  $V$-module. Moreover, for any irreducible $\varphi(h)$-twisted  $V$-module $N$, there exists an irreducible $V$-module $M$ such that  $N$ is isomorphic to $M^{(h)}$.
\end{proposition}

 Let $g$ be an automorphism of $V$ of finite order, $(M,Y_{M}(\cdot, z))$ be a $g$-twisted $V$-module. For a homogeneous vector $v\in V^0$, $Y_{M}(v, z)=\sum_{n\in \Z}v_nz^{-n-1} $. We define  $$o_M(v)=v_{\wt v-1}$$ and extend linearly. Then the following results have been proved in \cite{AE}.
\begin{proposition}\label{changeVE}
Let $(M,Y_{M}(\cdot, z))$ be an irreducible $V$-module and $h$ be as above. Then \\
(1) Let $v$ be an element in $V$ such that $h_0v=0$. Then $o_{M^{(h)}}(v)=o_M(\Delta(h, 1)v)$.\\
(2) Let $\hat L(0)$ and $L(0)$ be the operators defined by $Y_{M^{(h)}}(\w,z)=\sum_{n\in \Z}\hat L(n)z^{-n-2}$ and $Y_{M}(\w,z)=\sum_{n\in \Z} L(n)z^{-n-2}$, respectively. Then $\hat  L(0)=L(0)+h_0+\frac{1}{2}\langle h,h \rangle.$
\end{proposition}
\subsubsection{$S$-matrix in the orbifold theory of  affine vertex operator algebras} Let $k$ be a positive integer, $L_{\g(A)}(k\Lambda_0)$ be the affine vertex operator algebra defined in Subsection \ref{affineVOA}. Then $L_{\g(A)}(k\Lambda_0)$  is a rational and $C_2$-cofinite  vertex operator algebra. Let $\sigma$ be a diagram automorphism of $\overline{\g(A)}$, $\tilde\sigma$ be the automorphism of $L_{\g(A)}(k\Lambda_0)$ induced from $\sigma$, and $L_{\g(A)}(\lambda)$ be the irreducible highest weight $\g(A)$-module of highest weight $\lambda$. Then $L_{\g(A)}(\lambda)$ is $\tilde\sigma$-stable if and only if $\lambda\in (\widetilde{P^+_k(A)})^\sigma$ (see Page 526 of \cite{FRS}). By Theorem \ref{modular}, we have the following result.
\begin{theorem}
For any $\lambda^\dag\in \widetilde{P^+_k(A^\dag)}$ and $v\in L_{\g(A)}(k\Lambda_0)$, we have
$$Z_{L_{\g(A^\dag)}(\lambda^\dag)}\left(v,\left(\tilde\sigma^{-1},1\right),\frac{-1}{\tau}\right)=\tau^{\wt\left[v\right]}\sum_{\lambda\in (\widetilde{P^+_k(A)})^\sigma}S_{\lambda^\dag,\lambda}Z_{L_{\g(A)}(\lambda)}\left(v,\left(1,\tilde\sigma\right),\tau\right).$$
\end{theorem}

The aim in this subsection is to determine the matrix $(S_{\lambda^\dag,\lambda})_{\lambda^\dag\in \widetilde{P^+_k(A^\dag)}, \lambda\in (\widetilde{P^+_k(A)})^\sigma}$. Let $h\in \bar\h^0$ be an element satisfying ${\rm Spec}~ h_0\subset \frac{1}{T}\Z$ for some positive integer $T$. For any rational number $\epsilon$, $\varphi(\epsilon h)$ is an automorphism of  $L_{\g(A)}(k\Lambda_0)$ of finite order. Then we have the following result.
\begin{lemma}\label{stable1}
For any $\lambda^\dag\in \widetilde{P^+_k(A^\dag)}$, $L_{\g(A^\dag)}(\lambda^\dag)$ is $\varphi(\epsilon h)$-stable.
\end{lemma}
\pf Note that $h_0$ acts semisimply on $L_{\g(A^\dag)}(\lambda^\dag)$. Then $e^{2\pi i \epsilon h_0}$ is a well-defined operator on $L_{\g(A^\dag)}(\lambda^\dag)$. Moreover, for any $v\in L_{\g(A)}(k\Lambda_0)$, we have $$[h_0, Y_{L_{\g(A^\dag)}(\lambda^\dag)}(v,z)] =Y_{L_{\g(A^\dag)}(\lambda^\dag)}(h_0v,z).$$ This implies that
$e^{2\pi i \epsilon h_0}Y_{L_{\g(A^\dag)}(\lambda^\dag)}(v,z)e^{-2\pi i \epsilon h_0}=Y_{L_{\g(A^\dag)}(\lambda^\dag)}(e^{2\pi i \epsilon h_0}v,z)$. Thus, $L_{\g(A^\dag)}(\lambda^\dag)$ is $\varphi(\epsilon h)$-stable.
\qed

We also need the following result.
\begin{lemma}\label{stable2}
$L_{\g(A)}(\lambda)^{(\epsilon h)}$ is $\tilde\sigma$-stable if and only if $\lambda\in (\widetilde{P^+_k(A)})^\sigma$.
\end{lemma}
\pf For $\lambda\in (\widetilde{P^+_k(A)})^\sigma$, $L_{\g(A)}(\lambda)$ is $\tilde\sigma$-stable. Then there exists a linear isomorphism $\phi(\tilde\sigma):L_{\g(A)}(\lambda)\to L_{\g(A)}(\lambda)$ such that
$$\phi(\tilde\sigma)Y_{L_{\g(A)}(\lambda)}(v,z)\phi(\tilde\sigma)^{-1}=Y_{L_{\g(A)}(\lambda)}(\tilde\sigma(v),z).$$
Since $h\in \bar\h^0$, we have $\tilde\sigma h_n=h_n\tilde\sigma$ holds for any $n\in \Z_{\geq 0}$. This implies that
$$\phi(\tilde\sigma)Y_{L_{\g(A)}(\lambda)}(\Delta(h, z)v,z)\phi(\tilde\sigma)^{-1}=Y_{L_{\g(A)}(\lambda)}(\Delta(h,z)\tilde\sigma(v),z).$$
Therefore, $L_{\g(A)}(\lambda)^{(\epsilon h)}$ is $\tilde\sigma$-stable.
Similarly, if $L_{\g(A)}(\lambda)^{(\epsilon h)}$ is $\tilde\sigma$-stable, then $L_{\g(A)}(\lambda)$ is $\tilde\sigma$-stable. Then $\lambda\in (\widetilde{P^+_k(A)})^\sigma$.
\qed

By Theorem \ref{modular} and Lemmas \ref{stable1}, \ref{stable2}, we have the following result.
\begin{proposition}\label{main6}
For any $\lambda^\dag\in \widetilde{P^+_k(A^\dag)}$, we have
\begin{align*}
 &{\rm tr}_{L_{\g(A^\dag)}(\lambda^\dag)}o_{L_{\g(A^\dag)}(\lambda^\dag)}(v)e^{2\pi i\epsilon h_0}e^{2\pi i\frac{-1}{\tau}(L(0)-\frac{c}{24})}\\
 &=e^{2\pi i \tau\frac{1}{2}\langle \epsilon h,\epsilon h \rangle}\sum_{\lambda\in (\widetilde{P^+_k(A)})^\sigma} S_{L_{\g(A^\dag)}(\lambda^\dag), L_{\g(A)}(\lambda)^{(\epsilon h)}}\tau^{\wt[v]}{\rm tr}_{L_{\g(A)}(\lambda)}o_{L_{\g(A)}(\lambda)}(\Delta(\epsilon h, 1)v)\phi(\tilde \sigma)\\
 &\ \ \ \ \ \ \ \ \ \ \ \ \ \ \ \ \ \ \ \ \ \ \ \ \ \ \ \ \ \ \ \ \ \ \ \ \ \ \ \ \ \ \ \ \ \ \ \ \ \ \ \ \ \ \ \ \ \ \ \ \ \ \ \ \ \ \ \ \ \ \ \ \cdot e^{2\pi i\tau \epsilon h_0}e^{2\pi i\tau(L(0)-\frac{c}{24})}.
\end{align*}
\end{proposition}
\pf By Theorem \ref{modular}, Propositions \ref{inner}, \ref{changeVE} and Lemmas \ref{stable1}, \ref{stable2}, we have
\begin{align*}
&{\rm tr}_{L_{\g(A^\dag)}(\lambda^\dag)}o_{L_{\g(A^\dag)}(\lambda^\dag)}(v)e^{2\pi i\epsilon h_0}e^{2\pi i\frac{-1}{\tau}(L(0)-\frac{c}{24})}=Z_{L_{\g(A^\dag)}(\lambda^\dag)}(v, (\tilde\sigma^{-1}, \varphi(\epsilon h)), \frac{-1}{\tau})\\
 &=\sum_{\lambda\in (\widetilde{P^+_k(A)})^\sigma} S_{L_{\g(A^\dag)}(\lambda^\dag), L_{\g(A)}(\lambda)^{(\epsilon h)}}\tau^{\wt[v]}Z_{L_{\g(A)}(\lambda)^{(\epsilon h)}}(v, (\varphi(\epsilon h),\tilde\sigma ), \tau)\\
 &=\sum_{\lambda\in (\widetilde{P^+_k(A)})^\sigma} S_{L_{\g(A^\dag)}(\lambda^\dag), L_{\g(A)}(\lambda)^{(\epsilon h)}}\tau^{\wt[v]}{\rm tr}_{L_{\g(A)}(\lambda)^{(\epsilon h)}}o_{L_{\g(A)}(\lambda)^{(\epsilon h)}}(v)\phi(\tilde \sigma)e^{2\pi i\tau(\hat L(0)-\frac{c}{24})}\\
 &=\sum_{\lambda\in (\widetilde{P^+_k(A)})^\sigma} S_{L_{\g(A^\dag)}(\lambda^\dag), L_{\g(A)}(\lambda)^{(\epsilon h)}}\tau^{\wt[v]}{\rm tr}_{L_{\g(A)}(\lambda)}o_{L_{\g(A)}(\lambda)}(\Delta(\epsilon h, 1)v)\phi(\tilde \sigma)\\
 &\ \ \ \ \ \ \ \ \ \ \ \ \ \ \ \ \ \ \ \ \ \ \ \ \ \ \ \ \ \ \ \ \ \ \ \ \ \ \ \ \ \ \ \ \ \ \ \ \ \ \ \ \ \ \ \ \ \ \ \ \ \ \ \ \ \ \ \ \ \ \ \ \cdot e^{2\pi i\tau(L(0)+\epsilon h_0+\frac{1}{2}\langle \epsilon h,\epsilon h \rangle-\frac{c}{24})}\\
 &=e^{2\pi i \tau\frac{1}{2}\langle \epsilon h,\epsilon h \rangle}\sum_{\lambda\in (\widetilde{P^+_k(A)})^\sigma} S_{L_{\g(A^\dag)}(\lambda^\dag), L_{\g(A)}(\lambda)^{(\epsilon h)}}\tau^{\wt[v]}{\rm tr}_{L_{\g(A)}(\lambda)}o_{L_{\g(A)}(\lambda)}(\Delta(\epsilon h, 1)v)\phi(\tilde \sigma)\\
 &\ \ \ \ \ \ \ \ \ \ \ \ \ \ \ \ \ \ \ \ \ \ \ \ \ \ \ \ \ \ \ \ \ \ \ \ \ \ \ \ \ \ \ \ \ \ \ \ \ \ \ \ \ \ \ \ \ \ \ \ \ \ \ \ \ \ \ \ \ \ \ \ \cdot e^{2\pi i\tau \epsilon h_0}e^{2\pi i\tau(L(0)-\frac{c}{24})}.
\end{align*}
This completes the proof.
\qed

In case that $v=\1$, we have the following result.
\begin{corollary}\label{main4}
For any $\lambda^\dag\in \widetilde{P^+_k(A^\dag)}$, we have
\begin{align*}
 &{\rm tr}_{L_{\g(A^\dag)}(\lambda^\dag)}e^{2\pi i\epsilon h_0}e^{2\pi i\frac{-1}{\tau}(L(0)-\frac{c}{24})}\\
 &=e^{2\pi i \tau\frac{1}{2}\langle \epsilon h,\epsilon h \rangle}\sum_{\lambda\in (\widetilde{P^+_k(A)})^\sigma} S_{L_{\g(A^\dag)}(\lambda^\dag), L_{\g(A)}(\lambda)^{(\epsilon h)}}{\rm tr}_{L_{\g(A)}(\lambda)}\phi(\tilde \sigma) e^{2\pi i\epsilon\tau h_0}e^{2\pi i\tau(L(0)-\frac{c}{24})}.
\end{align*}
\end{corollary}

On the other hand, by Theorem \ref{KPmodular}, we have the following result.
\begin{theorem}\label{main5}
 Let $A^\dag$ be a twisted affine Cartan matrix not of type $A_{2l}^{(2)}$.
Then for any $z\in \C$ and  $\lambda^\dag\in \widetilde{P^+_k(A^\dag)}$, we have
\begin{align*}
 &{\rm tr}_{L_{\g(A^\dag)}(\lambda^\dag)}e^{2\pi i\frac{z}{\tau} h_0}e^{2\pi i\frac{-1}{\tau}(L(0)-\frac{c}{24})}\\
 &=e^{2\pi i k\frac{( zh,zh )}{2\tau}}\sum_{\lambda'\in \widetilde{P^+_k(A')}} a_{\lambda^\dag, \lambda'}{\rm tr}_{L_{\g(A)}(P_{\sigma}^*(\lambda'))}\phi(\tilde \sigma) e^{2\pi iz h_0}e^{2\pi i\tau(L(0)-\frac{c}{24})},
\end{align*}
where $a_{\lambda^\dag, \lambda'}$ is given by the formula (\ref{amatrix}).
\end{theorem}
\pf Let $A^\dag$ be a twisted affine Cartan matrix not of type $A_{2l}^{(2)}$, $A'$ be its adjacent Cartan matrix. By Theorem \ref{KPmodular},   for any $\lambda^\dag\in \widetilde{P^+_k(A^\dag)}$,
$$\chi_{L_{\g(A^\dag)}(\lambda^\dag)}(\frac{1}{-\tau}, \frac{\mathfrak{z}}{\tau}, t-\frac{(\mathfrak{z},\mathfrak{z})^\dag}{2\tau} )=\sum_{\lambda'\in \widetilde{P^+_k(A')}}a_{\lambda^\dag, \lambda'}\chi_{L_{\g(A')}(\lambda')}(\frac{\tau}{r}, \frac{\mathfrak{z}}{r}, t),$$
where  $r$ is the number such that $A^\dag$ belongs to Table Aff $r$ in \cite{K}. In case that $t=0$, we have
\begin{align*}
&e^{2\pi i\frac{-1}{\tau}m_{\lambda^\dag}}\sum_{\mu\in P(L_{\g(A^\dag)}(\lambda^\dag))}\dim L_{\g(A^\dag)}(\lambda^\dag)_\mu e^{2\pi i(\mu, \frac{1}{\tau}\Lambda_0^\dag+\frac{\mathfrak z}{\tau})^\dag}\\
&=e^{2\pi i k\frac{(\mathfrak z, \mathfrak z)^\dag}{2\tau}}\sum_{\lambda'\in \widetilde{P^+_k(A')}}a_{\lambda^\dag, \lambda'}e^{2\pi i\frac{\tau}{r}m_{\lambda'}}\sum_{\mu'\in P(L_{\g(A')}(\lambda'))}\dim L_{\g(A')}(\lambda')_{\mu'}e^{2\pi i(\mu', \frac{-\tau}{r}\Lambda_0'+\varphi^{-1}(\frac{\mathfrak z}{r}))'}.
\end{align*}
By Theorem \ref{key2}, Propositions \ref{key3},  \ref{key4},  Lemma \ref{orbitdual} and Corollary \ref{key5}, we have
\begin{align*}
&e^{2\pi i\frac{\tau}{r}m_{\lambda'}}\sum_{\mu'\in P(L_{\g(A')}(\lambda'))}\dim L_{\g(A')}(\lambda')_{\mu'}e^{2\pi i(\mu', \frac{-\tau}{r}\Lambda_0'+\varphi^{-1}(\frac{\mathfrak z}{r}))'}\\
&=e^{2\pi i\frac{\tau}{r}m_{\lambda'}}\sum_{\mu'\in P(L_{\g(A')}(\lambda'))}{\rm tr}_{L_{\g(A)}(P_{\sigma}^*(\lambda'))_{P_{\sigma}^*(\mu')}}\phi(\sigma)e^{2\pi i(P_{\sigma}^*(\mu'), P_{\sigma}^*(\frac{-\tau}{r}\Lambda_0'+\varphi^{-1}(\frac{\mathfrak z}{r})))}\\
&=e^{2\pi i\frac{\tau}{r}m_{P_{\sigma}^*(\lambda')}}{\rm tr}_{L_{\g(A)}(P_{\sigma}^*(\lambda'))}\phi(\sigma)e^{2\pi i \nu^{-1}(P_{\sigma}^*(\frac{-\tau}{r}\Lambda_0'+\varphi^{-1}(\frac{\mathfrak z}{r})))}\\
&=e^{2\pi i\frac{\tau}{r}m_{P_{\sigma}^*(\lambda')}}{\rm tr}_{L_{\g(A)}(P_{\sigma}^*(\lambda'))}\phi(\sigma)e^{2\pi i \nu^{-1}(\frac{-\tau}{r}\Lambda_0)}e^{2\pi i \nu^{-1}(P_{\sigma}^*\varphi^{-1}(\frac{\mathfrak z}{r}))}\\
&=e^{2\pi i\frac{\tau}{r}m_{P_{\sigma}^*(\lambda')}}{\rm tr}_{L_{\g(A)}(P_{\sigma}^*(\lambda'))}\phi(\sigma)e^{2\pi i \frac{-\tau}{r}d}e^{2\pi i \nu^{-1}(P_{\sigma}^*\varphi^{-1}(\frac{\mathfrak z}{r}))}\\
&={\rm tr}_{L_{\g(A)}(P_{\sigma}^*(\lambda'))}\phi(\sigma)e^{2\pi i \frac{\tau}{r}(L(0)-\frac{c}{24})}e^{2\pi i \nu^{-1}(P_{\sigma}^*\varphi^{-1}(\frac{\mathfrak z}{r}))}.
\end{align*}
On the other hand, by Proposition \ref{key3}, we have
\begin{align*}
&e^{2\pi i\frac{-1}{\tau}m_{\lambda^\dag}}\sum_{\mu\in P(L_{\g(A^\dag)}(\lambda^\dag))}\dim L_{\g(A^\dag)}(\lambda^\dag)_\mu e^{2\pi i(\mu, \frac{1}{\tau}\Lambda_0^\dag+\frac{\mathfrak z}{\tau})^\dag}\\
&=e^{2\pi i\frac{-1}{\tau}m_{\lambda^\dag}}{\rm tr}_{ L_{\g(A^\dag)}(\lambda^\dag)}e^{2\pi i\nu^{\dag-1}( \frac{1}{\tau}\Lambda_0^\dag+\frac{\mathfrak z}{\tau})}\\
&=e^{2\pi i\frac{-1}{\tau}m_{\lambda^\dag}}{\rm tr}_{ L_{\g(A^\dag)}(\lambda^\dag)}e^{2\pi i\nu^{\dag-1}( \frac{1}{\tau}\Lambda_0^\dag)}e^{2\pi i\nu^{\dag-1}(\frac{\mathfrak z}{\tau})}\\
&=e^{2\pi i\frac{-1}{\tau}m_{\lambda^\dag}}{\rm tr}_{ L_{\g(A^\dag)}(\lambda^\dag)}e^{2\pi i \frac{d^\dag}{\tau}}e^{2\pi i\nu^{\dag-1}(\frac{\mathfrak z}{\tau})}\\
&={\rm tr}_{ L_{\g(A^\dag)}(\lambda^\dag)}e^{2\pi i \frac{-r}{\tau}(L(0)-\frac{c}{24})}e^{2\pi i\nu^{\dag-1}(\frac{\mathfrak z}{\tau})}.
\end{align*}
Therefore, we have
\begin{align*}
&{\rm tr}_{ L_{\g(A^\dag)}(\lambda^\dag)}e^{2\pi i \frac{-r}{\tau}(L(0)-\frac{c}{24})}e^{2\pi i\nu^{\dag-1}(\frac{\mathfrak z}{\tau})}\\
&=e^{2\pi i k\frac{(\mathfrak z, \mathfrak z)^\dag}{2\tau}}\sum_{\lambda'\in \widetilde{P^+_k(A')}}a_{\lambda^\dag, \lambda'}{\rm tr}_{L_{\g(A)}(P_{\sigma}^*(\lambda'))}\phi(\sigma)e^{2\pi i \frac{\tau}{r}(L(0)-\frac{c}{24})}e^{2\pi i \nu^{-1}(P_{\sigma}^*\varphi^{-1}(\frac{\mathfrak z}{r}))}.
\end{align*}
By Proposition \ref{key1} and the formula (\ref{key6}), we have
\begin{align*}
&{\rm tr}_{ L_{\g(A^\dag)}(\lambda^\dag)}e^{2\pi i \frac{-r}{\tau}(L(0)-\frac{c}{24})}e^{2\pi i\nu^{\dag-1}(\frac{\mathfrak z}{\tau})}\\
&=e^{2\pi i k\frac{(\iota\nu^{\dag-1}(\mathfrak z), \iota\nu^{\dag-1}(\mathfrak z))}{2r\tau}}\sum_{\lambda'\in \widetilde{P^+_k(A')}}a_{\lambda^\dag, \lambda'}{\rm tr}_{L_{\g(A)}(P_{\sigma}^*(\lambda'))}\phi(\sigma)e^{2\pi i \frac{\tau}{r}(L(0)-\frac{c}{24})}e^{2\pi i \iota\nu^{\dag-1}(\frac{\mathfrak z}{r})}\\
&=e^{2\pi i kr\frac{(\iota\nu^{\dag-1}(\frac{\mathfrak z}{r})), \iota\nu^{\dag-1}(\frac{\mathfrak z}{r})))}{2\tau}}\sum_{\lambda'\in \widetilde{P^+_k(A')}}a_{\lambda^\dag, \lambda'}{\rm tr}_{L_{\g(A)}(P_{\sigma}^*(\lambda'))}\phi(\sigma)e^{2\pi i \frac{\tau}{r}(L(0)-\frac{c}{24})}e^{2\pi i \iota\nu^{\dag-1}(\frac{\mathfrak z}{r})}.
\end{align*}
Note that $\phi(\sigma)=\phi(\tilde \sigma)$. When $\iota\nu^{\dag-1}(\frac{\mathfrak z}{r})=zh$, we have
\begin{align*}
 &{\rm tr}_{L_{\g(A^\dag)}(\lambda^\dag)}e^{2\pi i\frac{z}{\tau} h_0}e^{2\pi i\frac{-1}{\tau}(L(0)-\frac{c}{24})}\\
 &=e^{2\pi i k\frac{( zh,zh )}{2\tau}}\sum_{\lambda'\in \widetilde{P^+_k(A')}} a_{\lambda^\dag, \lambda'}{\rm tr}_{L_{\g(A)}(P_{\sigma}^*(\lambda'))}\phi(\tilde \sigma) e^{2\pi iz h_0}e^{2\pi i\tau(L(0)-\frac{c}{24})}.
\end{align*}
This completes the proof.
\qed
\vskip.25cm
To determine the  matrix $(S_{\lambda^\dag,\lambda})_{\lambda^\dag\in \widetilde{P^+_k(A^\dag)}, \lambda\in (\widetilde{P^+_k(A)})^\sigma}$, we also need the following result.
\begin{lemma}\label{element}
There exists an element $h\in \bar\h^0$ such that\\
(1) ${\rm Spec}~ h_0\subset \frac{1}{T}\Z$ for some positive integer $T$.\\
(2) There exist $a, b\in \R$ such that the numbers $\epsilon\lambda(h)-h_{\lambda}-\frac{c}{24}$ $(\lambda\in (\widetilde{P^+_k(A)})^\sigma)$ are mutually distinct for any fixed number $\epsilon\in (a,b)$, where $h_{\lambda}=\frac{(\bar\lambda+2\bar\rho, \bar\lambda)}{2(k+h^\vee)}$ and $\rho$ is the Weyl vector of $\g(A)$.
\end{lemma}
\pf The equality of any pair of $\lambda(h)-h_{\lambda}-\frac{c}{24}$ defines a certain hyperplane in $\bar\h^0$. Therefore, there exists an element $h\in \bar\h^0$ such that ${\rm Spec}~ h_0\subset \frac{1}{T}\Z$ for some positive integer $T$ and the numbers $\lambda(h)-h_{\lambda}-\frac{c}{24}$ $(\lambda\in (\widetilde{P^+_k(A)})^\sigma)$ are mutually distinct. Furthermore, we can find numbers $a, b\in \R$ such that the numbers $\epsilon\lambda(h)-h_{\lambda}-\frac{c}{24}$ $(\lambda\in (\widetilde{P^+_k(A)})^\sigma)$ are mutually distinct for any fixed number $\epsilon\in (a,b)$.
\qed

\vskip.25cm
We are now ready to prove the main result in this subsection.
\begin{theorem}\label{modular1}
For any $\lambda^\dag\in \widetilde{P^+_k(A^\dag)}$ and $v\in L_{\g(A)}(k\Lambda_0)$, we have
$$Z_{L_{\g(A^\dag)}(\lambda^\dag)}\left(v,\left(\tilde\sigma^{-1},1\right),\frac{-1}{\tau}\right)=\tau^{\wt\left[v\right]}\sum_{\lambda'\in \widetilde{P^+_k(A')}} a_{\lambda^\dag, \lambda'}Z_{L_{\g(A)}(P_{\sigma}^*(\lambda'))}\left(v,\left(1,\tilde\sigma\right),\tau\right),$$
where $a_{\lambda^\dag, \lambda'}$ is given by the formula (\ref{amatrix}).
\end{theorem}
\pf In the following, we let $h\in \bar\h^0$ be an element satisfying the conditions (1), (2) in Lemma  \ref{element}. By Theorem 6.6.14 of \cite{LL}, the conformal weight of $L_{\g(A)}(P_{\sigma}^*(\lambda'))$ is equal to $h_{P_{\sigma}^*(\lambda')}$. Then for any fixed rational number $\epsilon\in (a, b)$, the functions $${\rm tr}_{L_{\g(A)}(P_{\sigma}^*(\lambda'))}\phi(\tilde \sigma) e^{2\pi i\epsilon\tau h_0}e^{2\pi i\tau(L(0)-\frac{c}{24})} ~(\lambda'\in \widetilde{P^+_k(A')})$$ are linearly independent. Therefore, by Corollary \ref{main4} and Theorem \ref{main5}, we have
$$S_{L_{\g(A^\dag)}(\lambda^\dag), L_{\g(A)}(P_{\sigma}^*(\lambda'))^{(\epsilon h)}}=a_{\lambda^\dag, \lambda'}.$$
By Proposition \ref{main6}, we have
\begin{align*}
 &{\rm tr}_{L_{\g(A^\dag)}(\lambda^\dag)}o_{L_{\g(A^\dag)}(\lambda^\dag)}(v)e^{2\pi i\epsilon h_0}e^{2\pi i\frac{-1}{\tau}(L(0)-\frac{c}{24})}\\
 &=e^{2\pi i \tau\frac{1}{2}\langle \epsilon h,\epsilon h \rangle}\sum_{\lambda'\in \widetilde{P^+_k(A')}} a_{\lambda^\dag, \lambda'}\tau^{\wt[v]}{\rm tr}_{L_{\g(A)}(P_{\sigma}^*(\lambda')}o_{L_{\g(A)}(P_{\sigma}^*(\lambda'))}(\Delta(\epsilon h, 1)v)\phi(\tilde \sigma)\\
 &\ \ \ \ \ \ \ \ \ \ \ \ \ \ \ \ \ \ \ \ \ \ \ \ \ \ \ \ \ \ \ \ \ \ \ \ \ \ \ \ \ \ \ \ \ \ \ \ \ \ \ \ \ \ \ \ \ \ \ \ \ \ \ \ \ \ \ \ \ \ \ \ \cdot e^{2\pi i\tau \epsilon h_0}e^{2\pi i\tau(L(0)-\frac{c}{24})}.
\end{align*}
Therefore,
\begin{align}\label{keyidentity}
 &{\rm tr}_{L_{\g(A^\dag)}(\lambda^\dag)}o_{L_{\g(A^\dag)}(\lambda^\dag)}(v)e^{2\pi i\frac{z}{\tau} h_0}e^{2\pi i\frac{-1}{\tau}(L(0)-\frac{c}{24})}\notag\\
 &=e^{2\pi i \frac{1}{2\tau}\langle z h,z h \rangle}\sum_{\lambda'\in \widetilde{P^+_k(A')}} a_{\lambda^\dag, \lambda'}\tau^{\wt[v]}{\rm tr}_{L_{\g(A)}(P_{\sigma}^*(\lambda')}o_{L_{\g(A)}(P_{\sigma}^*(\lambda'))}(\Delta(\frac{z}{\tau} h, 1)v)\phi(\tilde \sigma)\notag\\
 &\ \ \ \ \ \ \ \ \ \ \ \ \ \ \ \ \ \ \ \ \ \ \ \ \ \ \ \ \ \ \ \ \ \ \ \ \ \ \ \ \ \ \ \ \ \ \ \ \ \ \ \ \ \ \ \ \ \ \ \ \ \ \  \cdot e^{2\pi iz h_0}e^{2\pi i\tau(L(0)-\frac{c}{24})}
\end{align}
holds for $z=\epsilon \tau$. Note that these are holomorphic functions of $z$, and the identity (\ref{keyidentity}) holds for any rational number $\epsilon\in (a,b)$. This implies that
\begin{align*}
 &{\rm tr}_{L_{\g(A^\dag)}(\lambda^\dag)}o_{L_{\g(A^\dag)}(\lambda^\dag)}(v)e^{2\pi i\frac{z}{\tau} h_0}e^{2\pi i\frac{-1}{\tau}(L(0)-\frac{c}{24})}\\
 &=e^{2\pi i \frac{1}{2\tau}\langle z h,z h \rangle}\sum_{\lambda'\in \widetilde{P^+_k(A')}} a_{\lambda^\dag, \lambda'}\tau^{\wt[v]}{\rm tr}_{L_{\g(A)}(P_{\sigma}^*(\lambda')}o_{L_{\g(A)}(P_{\sigma}^*(\lambda'))}(\Delta(\frac{z}{\tau} h, 1)v)\phi(\tilde \sigma)\\
 &\ \ \ \ \ \ \ \ \ \ \ \ \ \ \ \ \ \ \ \ \ \ \ \ \ \ \ \ \ \ \ \ \ \ \ \ \ \ \ \ \ \ \ \ \ \ \ \ \ \ \ \ \ \ \ \ \ \ \ \ \ \ \ \ \ \ \ \ \ \ \ \ \cdot e^{2\pi iz h_0}e^{2\pi i\tau(L(0)-\frac{c}{24})}
\end{align*}
holds for any $z\in \C$. In case that $z=0$, we have
\begin{align*}
 &{\rm tr}_{L_{\g(A^\dag)}(\lambda^\dag)}o_{L_{\g(A^\dag)}(\lambda^\dag)}(v)e^{2\pi i\frac{-1}{\tau}(L(0)-\frac{c}{24})}\\
 &=\sum_{\lambda'\in \widetilde{P^+_k(A')}} a_{\lambda^\dag, \lambda'}\tau^{\wt[v]}{\rm tr}_{L_{\g(A)}(P_{\sigma}^*(\lambda')}o_{L_{\g(A)}(P_{\sigma}^*(\lambda'))}(v)\phi(\tilde \sigma)e^{2\pi i\tau(L(0)-\frac{c}{24})}.
\end{align*}
This completes the proof.
\qed
\subsection{Fusion rules between twisted modules of  affine vertex operator algebras} In this subsection, we determine fusion rules between twisted modules of  affine vertex operator algebras by using the twisted Verlinde formula.  We shall prove a twisted analogue of the Kac-Walton formula, which gives fusion rules between twisted modules of affine vertex operator algebras in terms of Clebsch-Gordan coefficients associated to the corresponding finite dimensional simple Lie algebras.
\subsubsection{Fusion rules between twisted modules of  affine vertex operator algebras: diagram automorphisms} Let $k$ be a positive integer, $L_{\g(A)}(k\Lambda_0)$ be the affine vertex operator algebra defined in Subsection \ref{affineVOA}. Then $L_{\g(A)}(k\Lambda_0)$  is a rational and $C_2$-cofinite  vertex operator algebra. Let $\sigma$ be a diagram automorphism of $\overline{\g(A)}$, $\tilde\sigma$ be the automorphism of $L_{\g(A)}(k\Lambda_0)$ induced from $\sigma$. The aim in this subsection is to determine  fusion rules between $\tilde\sigma$-twisted $L_{\g(A)}(k\Lambda_0)$-modules.

Let $A^\dag$ be a twisted affine Cartan matrix not of type $A_{2l}^{(2)}$, and $P^+(\bar A^\dag)$ be the set of dominant weights of $\overline{\g(A^\dag)}$. Then we have following result.
\begin{proposition}\label{weylgroup}
For any $\bar\lambda\in P^+(\bar A^\dag)$, there exists $w_0\in \overline{W}^\dag$ and $\mu^\dag\in P^+_{k+h^\vee}(A^\dag)$ such that $w_0(\bar \mu^\dag)=\bar\lambda+\bar\rho^\dag$ $({\rm mod }~ (k+h^\vee)M^\dag)$.
\end{proposition}
\pf The proof is similar to that of Proposition 6.6 of \cite{K}. Let $\h^\dag$ be the Cartan subalgebra of $\g(A^\dag)$, and $\h^\dag_{\R}$ be a vector space over $\R$ such that $\h^\dag_{\R}\otimes_{\R}\C\cong \h^\dag$. Set $$\h^{\dag*}_{k+h^\vee}=\{\mu^\dag\in \h^{\dag*}_{\R}|\langle\mu^\dag, K^{\dag}\rangle=k+h^\vee\}.$$Define a linear map
\begin{align*}
\Pi: \h^{\dag*}_{k+h^\vee}&\to \bar\h_{\R}^{\dag*}\\
\mu^\dag&\mapsto \bar\mu^\dag,
\end{align*}
where $\bar\mu^\dag$ denotes the projection of $\mu^\dag$ on $\bar\h_{\R}^{\dag*}$.
For any $w\in W^\dag$, we define
 \begin{align*}
\af(w): \bar\h_{\R}^{\dag*}&\to \bar\h_{\R}^{\dag*}\\
\bar\mu^\dag&\mapsto \af(w)(\bar\mu^\dag),
\end{align*}
 where $\af(w)(\bar\mu^\dag)=\Pi( w(\mu^\dag))$ and $\mu^\dag\in \h^{\dag*}_{k+h^\vee}$ is an element such that $\Pi(\mu^\dag)=\bar\mu^\dag$. Since $w(\delta^\dag)=\delta^\dag$ and $\h^{\dag*}_{k+h^\vee}$ is $W^\dag$-invariant,  $\af(w)$ is well-defined. By the definition of $\af(w)$, we have $\af(w_1)\circ\af(w_2)=\af(w_1w_2)$ for any $w_1, w_2\in W^\dag$. Moreover, we have $\Pi\circ w=\af(w)\circ\Pi$ for any $w\in W^\dag$.

  For any $\bar\lambda\in P^+(\bar A^\dag)$, $\bar\lambda+k\Lambda_0^\dag+\rho^\dag$ is an integral weight of $\g(A^\dag)$ of level $k+h^\vee$. Since the Weyl group $W^\dag$ preserves the set of integral weights of $\g(A^\dag)$, it follows from Lemma 3.32 of \cite{W1} that there exists $\mu^\dag\in P_{k+h^\vee}^+(A^\dag)$ and $w\in W^\dag$ such that $w(\mu^\dag)=\bar\lambda+k\Lambda_0^\dag+\rho^\dag$. Therefore, we have $\Pi(w(\mu^\dag))=\bar\lambda+\bar\rho^\dag$. This implies that $\af(w)(\bar \mu^\dag)=\bar\lambda+\bar\rho^\dag$. By Proposition \ref{Weylg}, there exists $\alpha\in M^\dag$ and $w_0\in \overline{W}^\dag$ such that $w=t_{\alpha}w_0$. As a result, we have $\af(w)(\bar \mu^\dag)=\af(t_{\alpha}w_0)(\bar \mu^\dag)=\af(t_{\alpha})\af(w_0)(\bar \mu^\dag)=\bar\lambda+\bar\rho^\dag$.

  We now show that $\af(w_0)(\bar \mu^\dag)=w_0(\bar \mu^\dag)$. By the definition of $\af(w_0)$, we have
  \begin{align*}
  \af(w_0)(\bar \mu^\dag)&=\Pi(w_0(\bar \mu^\dag+(k+h^\vee)\Lambda_0^\dag))=\Pi(w_0(\bar \mu^\dag)+w_0((k+h^\vee)\Lambda_0^\dag))\\
  &=\Pi(w_0(\bar \mu^\dag)+(k+h^\vee)\Lambda_0^\dag)=w_0(\bar \mu^\dag).
  \end{align*}

  Next we show that  $\af(t_{\alpha})(w_0(\bar \mu^\dag))=w_0(\bar \mu^\dag)~ ({\rm mod}~(k+h^\vee)M^\dag)$. In case that $\alpha=\nu^\dag(\theta^{\dag\vee})$, $t_{\alpha}=r_{\alpha_0^\dag}r_{\theta^\dag}$ by the formula (6.5.2) of \cite{K}. Then we have
  \begin{align*}
  \af(t_{\alpha})(w_0(\bar \mu^\dag))=\af(r_{\alpha_0^\dag}r_{\theta^\dag})(w_0(\bar \mu^\dag))=\af(r_{\alpha_0^\dag})\af(r_{\theta^\dag})(w_0(\bar \mu^\dag))=\af(r_{\alpha_0^\dag})(r_{\theta^\dag}(w_0(\bar \mu^\dag))).
  \end{align*}
  On the other hand,  for any $\Lambda\in \bar\h_{\R}^{\dag*}$, we have
  \begin{align}\label{weyltr}
  &\af(r_{\alpha_0^\dag})(\Lambda)=\Pi(r_{\alpha_0^\dag}(\Lambda+(k+h^\vee)\Lambda_0^\dag))\notag\\
  &=\Pi(\Lambda+(k+h^\vee)\Lambda_0^\dag-\langle \Lambda+(k+h^\vee)\Lambda_0^\dag, \alpha_0^{\dag\vee}\rangle\alpha_0^{\dag})\notag\\
  &=\Pi(\Lambda+(k+h^\vee)\Lambda_0^\dag-\langle \Lambda+(k+h^\vee)\Lambda_0^\dag, K^\dag-\theta^{\dag\vee}\rangle(\delta^\dag-\theta^\dag))\notag\\
  &=\Pi(\Lambda+(k+h^\vee)\Lambda_0^\dag+\langle \Lambda, \theta^{\dag\vee}\rangle(\delta^\dag-\theta^\dag)- \langle(k+h^\vee)\Lambda_0^\dag, K^\dag-\theta^{\dag\vee}\rangle(\delta^\dag-\theta^\dag))\notag\\
  &=\Pi(\Lambda+(k+h^\vee)\Lambda_0^\dag+\langle \Lambda, \theta^{\dag\vee}\rangle(\delta^\dag-\theta^\dag)- (k+h^\vee)(\delta^\dag-\theta^\dag))\notag\\
   &=\Lambda-\langle \Lambda, \theta^{\dag\vee}\rangle\theta^\dag~~~({\rm mod }~ (k+h^\vee)M^\dag)\notag\\
  & =r_{\theta^\dag}(\Lambda)~~~({\rm mod }~ (k+h^\vee)M^\dag).
  \end{align}
 Thus, for $\alpha=\nu^\dag(\theta^{\dag\vee})$, we have
 \begin{align*}
  \af(t_{\alpha})(w_0(\bar \mu^\dag))=\af(r_{\alpha_0^\dag})(r_{\theta^\dag}(w_0(\bar \mu^\dag))) =r_{\theta^\dag}(r_{\theta^\dag}(w_0(\bar \mu^\dag)))=w_0(\bar \mu^\dag)~~~({\rm mod }~ (k+h^\vee)M^\dag).
  \end{align*}
  For any $\alpha\in M^\dag$, we have $\alpha=w(\nu^\dag(\theta^{\dag\vee}))$ for some $w\in \overline{W}^\dag$. Therefore, $$t_{\alpha}=t_{w(\nu^\dag(\theta^{\dag\vee}))}=wt_{\nu^\dag(\theta^{\dag\vee})}w^{-1}.$$
  This implies that $\af(t_{\alpha})(w_0(\bar \mu^\dag))=w_0(\bar \mu^\dag)~ ({\rm mod}~(k+h^\vee)M^\dag)$ holds for any $\alpha\in M^\dag$.
  As a result, $\bar\lambda+\bar\rho^\dag=\af(w)(\bar \mu^\dag)=\af(t_{\alpha})(w_0(\bar \mu^\dag))=w_0(\bar \mu^\dag)~~~({\rm mod }~ (k+h^\vee)M^\dag)$.
This completes the proof.
\qed

\vskip.25cm
As a corollary, we have
\begin{corollary}\label{weyl1}
For any $\bar\lambda\in P^+(\bar A^\dag)$, either $r_{\alpha}(\bar\lambda+\bar\rho^\dag)=\bar\lambda+\bar\rho^\dag$ $({\rm mod }~ (k+h^\vee)M^\dag)$ for some root $\alpha$ of $\overline{\g(A^\dag)}$, or else there exists a unique $w_0\in \overline{W}^\dag$ and unique $\lambda^\dag\in \widetilde{P^+_{k}(A^\dag)}$ such that $w_0(\bar \lambda^\dag+\bar\rho^\dag)=\bar\lambda+\bar\rho^\dag$ $({\rm mod }~ (k+h^\vee)M^\dag)$.
\end{corollary}
\pf Assume that $r_{\alpha}(\bar\lambda+\bar\rho^\dag)\neq \bar\lambda+\bar\rho^\dag$ $({\rm mod }~ (k+h^\vee)M^\dag)$ for any root $\alpha$ of $\overline{\g(A^\dag)}$.
By Lemma 3.32 of \cite{W1}, there exists $\mu^\dag\in P_{k+h^\vee}^+(A^\dag)$ and $w\in W^\dag$ such that $w(\mu^\dag)=\bar\lambda+k\Lambda_0^\dag+\rho^\dag$. We next show that  $\langle \mu^\dag, \alpha_i^{\dag\vee}\rangle>0$.

If $\langle \mu^\dag, \alpha_0^{\dag\vee}\rangle=0$, then we have $r_{\alpha_0^\dag}(\mu^\dag)=\mu^\dag$. On the other hand, by the similar argument as that in the formula (\ref{weyltr}), we have $\overline{r_{\alpha_0^\dag}(\mu^\dag)}=r_{\theta^\dag}(\bar\mu^\dag)~({\rm mod }~ (k+h^\vee)M^\dag)$. In particular, $\bar\mu^\dag=r_{\theta^\dag}(\bar\mu^\dag)~({\rm mod }~ (k+h^\vee)M^\dag)$. By Proposition \ref{weylgroup}, there exists $w_0\in \overline{W}^\dag$  such that $w_0(\bar \mu^\dag)=\bar\lambda+\bar\rho^\dag~({\rm mod }~ (k+h^\vee)M^\dag).$
This implies that $r_{\theta^\dag}w_0^{-1}(\bar\lambda+\bar\rho^\dag)=w_0^{-1}(\bar\lambda+\bar\rho^\dag)~({\rm mod }~ (k+h^\vee)M^\dag)$. In particular, $w_0 r_{\theta^\dag}w_0^{-1}(\bar\lambda+\bar\rho^\dag)=\bar\lambda+\bar\rho^\dag~({\rm mod }~ (k+h^\vee)M^\dag)$. Therefore, $ r_{w_0 (\theta^\dag)}(\bar\lambda+\bar\rho^\dag)=\bar\lambda+\bar\rho^\dag~({\rm mod }~ (k+h^\vee)M^\dag)$, this is a contradiction.

If $\langle \mu^\dag, \alpha_i^{\dag\vee}\rangle=0$ for $i\neq 0$. Then we have $r_{\alpha_i^\dag}(\mu^\dag)=\mu^\dag$. On the other hand,
$\overline{r_{\alpha_i^\dag}(\mu^\dag)}=r_{\alpha_i^\dag}(\bar\mu^\dag)$. In particular, $\bar\mu^\dag=r_{\alpha_i^\dag}(\bar\mu^\dag)$. This implies that $r_{\alpha_i^\dag}w_0^{-1}(\bar\lambda+\bar\rho^\dag)=w_0^{-1}(\bar\lambda+\bar\rho^\dag)~({\rm mod }~ (k+h^\vee)M^\dag)$. In particular, $w_0 r_{\alpha_i^\dag}w_0^{-1}(\bar\lambda+\bar\rho^\dag)=\bar\lambda+\bar\rho^\dag~({\rm mod }~ (k+h^\vee)M^\dag)$. Therefore, $ r_{w_0 (\alpha_i^\dag)}(\bar\lambda+\bar\rho^\dag)=\bar\lambda+\bar\rho^\dag~({\rm mod }~ (k+h^\vee)M^\dag)$, this is a contradiction.

 Therefore, $\langle \mu^\dag, \alpha_i^{\dag\vee}\rangle>0$. Set $\lambda^\dag=\bar\mu^\dag+(k+h^\vee)\Lambda_0-\rho^\dag$. Then we have $\lambda^\dag\in \widetilde{P^+_{k}(A^\dag)}$. By Proposition \ref{weylgroup}, there exists $w_0\in \overline{W}^\dag$  such that $$w_0(\bar \lambda^\dag+\bar\rho^\dag)=w_0(\bar \mu^\dag)=\bar\lambda+\bar\rho^\dag~({\rm mod }~ (k+h^\vee)M^\dag).$$
We now show that there exists a unique $w_0\in \overline{W}^\dag$ and unique $\lambda^\dag\in \widetilde{P^+_{k}(A^\dag)}$ such that $w_0(\bar \lambda^\dag+\bar\rho^\dag)=\bar\lambda+\bar\rho^\dag$ $({\rm mod }~ (k+h^\vee)M^\dag)$. This follows from the fact that the set $\{\bar\rho^\dag+\bar\beta^\dag|\beta^\dag\in \widetilde{P^+_{k}(A^\dag)}\}$ sits in the interior of fundamental alcove with respect to the action of $\overline{W}^\dag\ltimes (k+h^\vee )M^\dag$ on $\bar\h^{\dag*}_{\R}$ (see Page 18 of \cite{HK}). This completes the proof.
\qed

\vskip.25cm
Let $A$ be an affine Cartan matrix of type $A_{2n-1}^{(1)}$, $D_{n+1}^{(1)}$, $D_{4}^{(1)}$ or $E_6^{(1)}$, $A^\dag$ be an affine Cartan matrix of type $A_{2n-1}^{(2)}$, $D_{n+1}^{(2)}$, $D_{4}^{(3)}$ or $E_6^{(2)}$, respectively. Then there is an embedding $\iota:\overline{\g(A^\dag)}\to \overline{\g(A)}$. Let $\lambda_1\in \widetilde{P^+_{k}(A)}$, $\lambda^\dag_2\in \widetilde{P^+_{k}(A^\dag)}$. Then $L_{\overline{\g(A)}}(\bar \lambda_1)$ may be viewed as a $\overline{\g(A^\dag)}$-module. Consider the tensor product $L_{\overline{\g(A)}}(\bar \lambda_1)\otimes L_{\overline{\g(A^\dag)}}(\bar \lambda_2)$ of $\overline{\g(A^\dag)}$-modules $L_{\overline{\g(A)}}(\bar \lambda_1)$ and $L_{\overline{\g(A^\dag)}}(\bar \lambda_2)$, we then have the following decomposition
\begin{align}\label{decomp}
L_{\overline{\g(A)}}(\bar \lambda_1)\otimes L_{\overline{\g(A^\dag)}}(\bar \lambda_2)=\bigoplus_{\mu\in P^+(\bar A^\dag)}\mult_{\bar\lambda_1\otimes \bar\lambda_2}(\mu) L_{\overline{\g(A^\dag)}}(\mu),
 \end{align}
 where $\mult_{\bar\lambda_1\otimes \bar\lambda_2}(\mu)$ denotes the multiplicity of $L_{\overline{\g(A^\dag)}}(\mu)$ in $L_{\overline{\g(A)}}(\bar \lambda_1)\otimes L_{\overline{\g(A^\dag)}}(\bar \lambda_2)$. We are now ready to prove the main result in this subsection.
\begin{theorem}\label{mainfusion}
Let $k$ be a positive integer, $A$ be an affine Cartan matrix of type $A_{2n-1}^{(1)}$, $D_{n+1}^{(1)}$, $D_{4}^{(1)}$ or $E_6^{(1)}$, $A^\dag$ be an affine Cartan matrix of type $A_{2n-1}^{(2)}$, $D_{n+1}^{(2)}$, $D_{4}^{(3)}$ or $E_6^{(2)}$, respectively. Let $\lambda_1\in \widetilde{P^+_{k}(A)}$, $\lambda^\dag_2, \lambda^\dag_3\in \widetilde{P^+_{k}(A^\dag)}$. Then we have
$$N_{L_{\g(A)}(\lambda_1), L_{\g(A^\dag)}(\lambda_2^\dag)}^{ L_{\g(A^\dag)}( \lambda_3^\dag)}=\sum_{\substack{\mu\in P^+(\bar A^\dag),\\\mu+\bar\rho^\dag=w(\bar \lambda_3^\dag+\bar\rho^\dag)~({\rm mod}~ (k+h^\vee)M^\dag)\\(\text{ for some } w\in \overline{W}^\dag)}}\epsilon(w)\mult_{\bar\lambda_1\otimes \bar\lambda_2^\dag}(\mu).$$
\end{theorem}
\pf  By the decomposition (\ref{decomp}), we have
\begin{align*}
ch_{L_{\overline{\g(A)}}(\bar \lambda_1)\otimes L_{\overline{\g(A^\dag)}}(\bar \lambda_2^\dag)}(h)=\sum_{\mu\in P^+(\bar A^\dag)}\mult_{\bar\lambda_1\otimes \bar\lambda_2^\dag}(\mu) ch_{L_{\overline{\g(A^\dag)}}(\mu)}(h),
 \end{align*}
 holds for any $h\in \bar\h^\dag$. This implies that
 \begin{align*}
ch_{L_{\overline{\g(A)}}(\bar \lambda_1)}(\iota(h))ch_{ L_{\overline{\g(A^\dag)}}(\bar \lambda_2^\dag)}(h)=\sum_{\mu\in P^+(\bar A^\dag)}\mult_{\bar\lambda_1\otimes \bar\lambda_2^\dag}(\mu) ch_{L_{\overline{\g(A^\dag)}}(\mu)}(h).
 \end{align*}
 Note that for any  $\lambda'\in \widetilde{P^+_{k}(A')}$, we have $\nu^{\dag-1}(\varphi(\bar\lambda'+\bar\rho'))\in \bar\h^\dag$. Then we have
 \begin{align*}
&ch_{L_{\overline{\g(A)}}(\bar \lambda_1)}(\iota(\frac{-2\pi i\nu^{\dag-1}(\varphi(\bar\lambda'+\bar\rho'))}{k+h^\vee}))ch_{ L_{\overline{\g(A^\dag)}}(\bar \lambda_2^\dag)}(\frac{-2\pi i\nu^{\dag-1}(\varphi(\bar\lambda'+\bar\rho'))}{k+h^\vee})\\
&=\sum_{\mu\in P^+(\bar A^\dag)}\mult_{\bar\lambda_1\otimes \bar\lambda_2^\dag}(\mu) ch_{L_{\overline{\g(A^\dag)}}(\mu)}(\frac{-2\pi i\nu^{\dag-1}(\varphi(\bar\lambda'+\bar\rho'))}{k+h^\vee}).
 \end{align*}

 If $r_{\alpha}(\bar\mu+\bar\rho^\dag)=\bar\mu+\bar\rho^\dag$ $({\rm mod }~ (k+h^\vee)M^\dag)$ for some root $\alpha$ of $\overline{\g(A^\dag)}$, it follows from (13.9.1) of \cite{K} that
 $$ch_{L_{\overline{\g(A^\dag)}}(\mu)}(\frac{-2\pi i\nu^{\dag-1}(\varphi(\bar\lambda'+\bar\rho'))}{k+h^\vee})=ch_{L_{\overline{\g(A^\dag)}}(r_{\alpha}(\mu+\bar\rho^\dag)-\bar\rho^\dag)}(\frac{-2\pi i\nu^{\dag-1}(\varphi(\bar\lambda'+\bar\rho'))}{k+h^\vee}).$$
On the other hand,  $ch_{L_{\overline{\g(A^\dag)}}(\mu)}=\epsilon(w)ch_{L_{\overline{\g(A^\dag)}}(\mu_1)}$ if $\mu=w(\mu_1+\bar\rho^\dag)-\bar\rho^\dag$. This implies that
$$ch_{L_{\overline{\g(A^\dag)}}(\mu)}(\frac{-2\pi i\nu^{\dag-1}(\varphi(\bar\lambda'+\bar\rho'))}{k+h^\vee})=0.$$
 By Corollary \ref{weyl1}, we have
\begin{align*}
&ch_{L_{\overline{\g(A)}}(\bar \lambda_1)}(\iota(\frac{-2\pi i\nu^{\dag-1}(\varphi(\bar\lambda'+\bar\rho'))}{k+h^\vee}))ch_{ L_{\overline{\g(A^\dag)}}(\bar \lambda_2^\dag)}(\frac{-2\pi i\nu^{\dag-1}(\varphi(\bar\lambda'+\bar\rho'))}{k+h^\vee})\\
&=\sum_{\lambda^\dag\in \widetilde{P^+_{k}(A^\dag)}}\sum_{\substack{\mu\in P^+(\bar A^\dag),\\\mu+\bar\rho^\dag=w(\bar \lambda^\dag+\bar\rho^\dag)~({\rm mod}~ (k+h^\vee)M^\dag)\\(\text{ for some } w\in \overline{W}^\dag)}}\mult_{\bar\lambda_1\otimes \bar\lambda_2^\dag}(\mu) ch_{L_{\overline{\g(A^\dag)}}(\mu)}(\frac{-2\pi i\nu^{\dag-1}(\varphi(\bar\lambda'+\bar\rho'))}{k+h^\vee})\\
&=\sum_{\lambda^\dag\in \widetilde{P^+_{k}(A^\dag)}}W_{L_{\g(A)}(\lambda_1), L_{\g(A^\dag)}(\lambda_2^\dag)}^{ L_{\g(A^\dag)}( \lambda^\dag)} ch_{L_{\overline{\g(A^\dag)}}(\bar\lambda^\dag)}(\frac{-2\pi i\nu^{\dag-1}(\varphi(\bar\lambda'+\bar\rho'))}{k+h^\vee}),
 \end{align*}
 where
 $$W_{L_{\g(A)}(\lambda_1), L_{\g(A^\dag)}(\lambda_2^\dag)}^{ L_{\g(A^\dag)}( \lambda^\dag)}=\sum_{\substack{\mu\in P^+(\bar A^\dag),\\\mu+\bar\rho^\dag=w(\bar \lambda^\dag+\bar\rho^\dag)~({\rm mod}~ (k+h^\vee)M^\dag)\\(\text{ for some } w\in \overline{W}^\dag)}}\epsilon(w)\mult_{\bar\lambda_1\otimes \bar\lambda_2^\dag}(\mu).$$
 Combining with the Weyl character formula yields
 \begin{align*}
 &\frac{\sum_{w\in \overline{W}}\epsilon(w)e^{\frac{-2\pi i\langle w(\bar\lambda_1+\bar\rho), \iota(\nu^{\dag-1}(\varphi(\bar\lambda'+\bar\rho')))\rangle}{k+h^\vee}}}{\sum_{w\in \overline{W}}\epsilon(w)e^{\frac{-2\pi i\langle w(\bar\rho), \iota(\nu^{\dag-1}(\varphi(\bar\lambda'+\bar\rho')))\rangle}{k+h^\vee}}}
 \frac{\sum_{w\in \overline{W}^\dag}\epsilon(w)e^{\frac{-2\pi i\langle w(\bar\lambda_2^\dag+\bar\rho^\dag), \nu^{\dag-1}(\varphi(\bar\lambda'+\bar\rho'))\rangle}{k+h^\vee}}}{\sum_{w\in \overline{W}^\dag}\epsilon(w)e^{\frac{-2\pi i\langle w(\bar\rho^\dag), \nu^{\dag-1}(\varphi(\bar\lambda'+\bar\rho'))\rangle}{k+h^\vee}}}\\
 &=\sum_{\lambda^\dag\in \widetilde{P^+_{k}(A^\dag)}}W_{L_{\g(A)}(\lambda_1), L_{\g(A^\dag)}(\lambda_2^\dag)}^{ L_{\g(A^\dag)}( \lambda^\dag)} \frac{\sum_{w\in \overline{W}^\dag}\epsilon(w)e^{\frac{-2\pi i\langle w(\bar\lambda^\dag+\bar\rho^\dag), \nu^{\dag-1}(\varphi(\bar\lambda'+\bar\rho'))\rangle}{k+h^\vee}}}{\sum_{w\in \overline{W}^\dag}\epsilon(w)e^{\frac{-2\pi i\langle w(\bar\rho^\dag), \nu^{\dag-1}(\varphi(\bar\lambda'+\bar\rho'))\rangle}{k+h^\vee}}}.
 \end{align*}
 Using Proposition \ref{key1} and multiplying both sides by
 $$\frac{\sum_{w\in \overline{W}^\dag}\epsilon(w)e^{\frac{-2\pi i\langle w(\bar\rho^\dag), \nu^{\dag-1}(\varphi(\bar\lambda'+\bar\rho'))\rangle}{k+h^\vee}}}{\sum_{w\in \overline{W}}\epsilon(w)e^{\frac{-2\pi i\langle w(\bar\rho), \iota(\nu^{\dag-1}(\varphi(\bar\lambda'+\bar\rho')))\rangle}{k+h^\vee}}},$$
 we have
\begin{align*}
 &\frac{\sum_{w\in \overline{W}}\epsilon(w)e^{\frac{-2\pi i( w(\bar\lambda_1+\bar\rho), P^*_{\sigma}(\bar\lambda'+\bar\rho'))}{k+h^\vee}}}{\sum_{w\in \overline{W}}\epsilon(w)e^{\frac{-2\pi i( w(\bar\rho), P^*_{\sigma}(\bar\lambda'+\bar\rho'))}{k+h^\vee}}}
 \frac{\sum_{w\in \overline{W}^\dag}\epsilon(w)e^{\frac{-2\pi i( w(\bar\lambda_2^\dag+\bar\rho^\dag), \varphi(\bar\lambda'+\bar\rho'))^\dag}{k+h^\vee}}}{\sum_{w\in \overline{W}}\epsilon(w)e^{\frac{-2\pi i( w(\bar\rho), P^*_{\sigma}(\bar\lambda'+\bar\rho'))}{k+h^\vee}}}\\
 &=\sum_{\lambda^\dag\in \widetilde{P^+_{k}(A^\dag)}}W_{L_{\g(A)}(\lambda_1), L_{\g(A^\dag)}(\lambda_2^\dag)}^{ L_{\g(A^\dag)}( \lambda^\dag)}\frac{\sum_{w\in \overline{W}^\dag}\epsilon(w)e^{\frac{-2\pi i( w(\bar\lambda^\dag+\bar\rho^\dag), \varphi(\bar\lambda'+\bar\rho'))^\dag}{k+h^\vee}}}{\sum_{w\in \overline{W}}\epsilon(w)e^{\frac{-2\pi i( w(\bar\rho), P^*_{\sigma}(\bar\lambda'+\bar\rho'))}{k+h^\vee}}}.
 \end{align*}
 By Theorem \ref{KPmodular} and Lemma \ref{orbitdual}, we have
 \begin{align*}
 \frac{a_{\lambda_1, P^*_{\sigma}(\lambda')}}{a_{k\Lambda_0, P^*_{\sigma}(\lambda')}}
 \frac{a_{\lambda_2^\dag, \lambda'}}{a_{k\Lambda_0, P^*_{\sigma}(\lambda')}}=\sum_{\lambda^\dag\in \widetilde{P^+_{k}(A^\dag)}}W_{L_{\g(A)}(\lambda_1), L_{\g(A^\dag)}(\lambda_2^\dag)}^{ L_{\g(A^\dag)}( \lambda^\dag)}\frac{a_{\lambda^\dag, \lambda'}}{a_{k\Lambda_0, P^*_{\sigma}(\lambda')}}.
 \end{align*}
 By Lemma 4.3 of \cite{DKR} and Theorem \ref{modular1}, we have
 \begin{align*}
 &\frac{S_{L_{\g(A)}(\lambda_1), L_{\g(A)}(P^*_{\sigma}(\lambda'))}}{S_{L_{\g(A)}(k\Lambda_0), L_{\g(A)}(P^*_{\sigma}(\lambda'))}}
 \frac{S_{L_{\g(A^\dag)}(\lambda_2^\dag), L_{\g(A)}(P^*_{\sigma}(\lambda'))}}{S_{L_{\g(A)}(k\Lambda_0), L_{\g(A)}(P^*_{\sigma}(\lambda'))}}\\
 &=\sum_{\lambda^\dag\in \widetilde{P^+_{k}(A^\dag)}}W_{L_{\g(A)}(\lambda_1), L_{\g(A^\dag)}(\lambda_2^\dag)}^{ L_{\g(A^\dag)}( \lambda^\dag)}\frac{S_{L_{\g(A^\dag)}(\lambda^\dag), L_{\g(A)}(P^*_{\sigma}(\lambda'))}}{S_{L_{\g(A)}(k\Lambda_0), L_{\g(A)}(P^*_{\sigma}(\lambda'))}}.
 \end{align*}
 Hence,  by Lemma 5.2 of \cite{Ho2} and Proposition 4.6 of \cite{KP}, we have
 \begin{align*}
 &\sum_{\lambda'\in \widetilde{P^+_{k}(A')}}\frac{S_{L_{\g(A)}(\lambda_1), L_{\g(A)}(P^*_{\sigma}(\lambda'))}S_{L_{\g(A^\dag)}(\lambda_2^\dag), L_{\g(A)}(P^*_{\sigma}(\lambda'))}\overline{S_{L_{\g(A^\dag)}(\lambda_3^\dag),L_{\g(A)}(P^*_{\sigma}(\lambda'))}}}{S_{L_{\g(A)}(k\Lambda_0), L_{\g(A)}(P^*_{\sigma}(\lambda'))}}\\
 &=\sum_{\lambda'\in \widetilde{P^+_{k}(A')}}\sum_{\lambda^\dag\in \widetilde{P^+_{k}(A^\dag)}}W_{L_{\g(A)}(\lambda_1), L_{\g(A^\dag)}(\lambda_2^\dag)}^{ L_{\g(A^\dag)}( \lambda^\dag)}S_{L_{\g(A^\dag)}(\lambda^\dag), L_{\g(A)}(P^*_{\sigma}(\lambda'))}\overline{S_{L_{\g(A^\dag)}(\lambda_3^\dag),L_{\g(A)}(P^*_{\sigma}(\lambda'))}}\\
&=\sum_{\lambda^\dag\in \widetilde{P^+_{k}(A^\dag)}}W_{L_{\g(A)}(\lambda_1), L_{\g(A^\dag)}(\lambda_2^\dag)}^{ L_{\g(A^\dag)}( \lambda^\dag)}\sum_{\lambda'\in \widetilde{P^+_{k}(A')}}S_{L_{\g(A^\dag)}(\lambda^\dag), L_{\g(A)}(P^*_{\sigma}(\lambda'))}\overline{S_{L_{\g(A^\dag)}(\lambda_3^\dag),L_{\g(A)}(P^*_{\sigma}(\lambda'))}}\\
&=\sum_{\lambda^\dag\in \widetilde{P^+_{k}(A^\dag)}}W_{L_{\g(A)}(\lambda_1), L_{\g(A^\dag)}(\lambda_2^\dag)}^{ L_{\g(A^\dag)}( \lambda^\dag)}\delta_{\lambda^\dag, \lambda_3^\dag}=W_{L_{\g(A)}(\lambda_1), L_{\g(A^\dag)}(\lambda_2^\dag)}^{ L_{\g(A^\dag)}( \lambda_3^\dag)}.
 \end{align*}
 By Proposition \ref{key4} and Theorem \ref{main}, we have
 $$N_{L_{\g(A)}(\lambda_1), L_{\g(A^\dag)}(\lambda_2^\dag)}^{ L_{\g(A^\dag)}( \lambda_3^\dag)}=\sum_{\substack{\mu\in P^+(\bar A^\dag),\\\mu+\bar\rho^\dag=w(\bar \lambda_3^\dag+\bar\rho^\dag)~({\rm mod}~ (k+h^\vee)M^\dag)\\(\text{ for some } w\in \overline{W}^\dag)}}\epsilon(w)\mult_{\bar\lambda_1\otimes \bar\lambda_2^\dag}(\mu).$$
 This completes the proof.
\qed
\begin{remark}
(1) Theorem \ref{mainfusion} may be viewed as a twisted analogue of the Kac-Walton formula (see the formula (\ref{Kac-Walton})).

(2) By Proposition \ref{symmetric} and Corollary 6.3 of \cite{H3}, fusion product  between two $\tilde\sigma$-twisted $L_{\g(A)}(k\Lambda_0)$-modules can be determined explicitly when $A$ is an affine Cartan matrix of type $A_{2n-1}^{(1)}$, $D_{n+1}^{(1)}$ or $E_6^{(1)}$ when the order of $\tilde\sigma$ is 2. If the order of $\tilde\sigma$ is 3, the treatment is more complicated and will appear in somewhere else.
\end{remark}
\subsubsection{Fusion rules between twisted modules of  affine vertex operator algebras: general automorphisms}
Let $k$ be a positive integer, $L_{\g(A)}(k\Lambda_0)$ be the affine vertex operator algebra defined in Subsection \ref{affineVOA}. Let $\sigma$ be a diagram automorphism of $\overline{\g(A)}$, $\tilde\sigma$ be the automorphism of $L_{\g(A)}(k\Lambda_0)$ induced from $\sigma$. Let $h\in \bar\h^0$ be an element satisfying ${\rm Spec}~ h_0\subset \frac{1}{T}\Z$ for some positive integer $T$. By the discussion in Subsection \ref{Li}, we have an automorphism $\varphi(h)$ of $L_{\g(A)}(k\Lambda_0)$. Since $h\in \bar\h^0$, $\varphi(h)$  commutes with $\tilde\sigma$. The aim in this subsection is to determine  fusion rules between $\tilde\sigma\varphi(h)$-twisted $L_{\g(A)}(k\Lambda_0)$-modules.

 By Proposition 5.4 of \cite{Li1}, $\{L_{\g(A^\dag)}(\lambda^\dag)^{(h)}|\lambda^\dag\in \widetilde{P^+_k(A^\dag)}\}$ is the complete set of irreducible $\tilde\sigma\varphi(h)$-twisted modules of $L_{\g(A)}(k\Lambda_0)$. On the other hand, the following result has been proved in Lemma 2.7 of \cite{DLM1}.
\begin{proposition}
Let $L_{\g(A)}(\lambda_1)$ be an irreducible $L_{\g(A)}(k\Lambda_0)$-module, $L_{\g(A^\dag)}(\lambda_2^\dag)$ and $ L_{\g(A^\dag)}( \lambda_3^\dag)$ be irreducible $\tilde\sigma$-twisted  $L_{\g(A)}(k\Lambda_0)$-modules, $I(\cdot, z)$ be an intertwining operator of type $\binom{L_{\g(A^\dag)}( \lambda_3^\dag)}{L_{\g(A)}(\lambda_1)\ L_{\g(A^\dag)}( \lambda_2^\dag)}$. Then $I(\Delta(h,z)\cdot, z)$ is an intertwining operator of type $\binom{L_{\g(A^\dag)}( \lambda_3^\dag)^{(h)}}{L_{\g(A)}(\lambda_1)\ L_{\g(A^\dag)}( \lambda_2^\dag)^{(h)}}$.
\end{proposition}

Furthermore, by Lemma 2.5 of \cite{DLM1}, we have the following result about fusion rules between $\tilde\sigma\varphi(h)$-twisted $L_{\g(A)}(k\Lambda_0)$-modules.
\begin{proposition}
Let $L_{\g(A)}(\lambda_1)$ be an irreducible $L_{\g(A)}(k\Lambda_0)$-module, $L_{\g(A^\dag)}(\lambda_2^\dag)$ and $ L_{\g(A^\dag)}( \lambda_3^\dag)$ be irreducible $\tilde\sigma$-twisted  $L_{\g(A)}(k\Lambda_0)$-modules. Then we have
$$N_{L_{\g(A)}(\lambda_1), L_{\g(A^\dag)}(\lambda_2^\dag)^{(h)}}^{ L_{\g(A^\dag)}( \lambda_3^\dag)^{(h)}}=N_{L_{\g(A)}(\lambda_1), L_{\g(A^\dag)}(\lambda_2^\dag)}^{ L_{\g(A^\dag)}( \lambda_3^\dag)}.$$
\end{proposition}

\end{document}